\documentclass[11pt]{article}
\date{}
\usepackage{graphicx}
\usepackage{enumerate,amssymb,amsmath,setspace,hyperref}
\usepackage{amssymb,amsfonts,amsthm,amscd}

\usepackage{caption}
\usepackage[]{subcaption}
\usepackage{tikz}
\usetikzlibrary{snakes}

\DeclareCaptionLabelFormat{opening}{}
\numberwithin{equation}{section}

\font\sml=cmr6

\newtheorem{thm}{Theorem}[section]
\newtheorem{lem}[thm]{Lemma}
\newtheorem{prop}[thm]{Proposition}
\newtheorem{cor}[thm]{Corollary}
\newtheorem{conj}[thm]{{Conjecture}}
\newtheorem{example}[thm]{{Example}}
\newtheorem{problem}[thm]{{Problem}}
\newtheorem{question}[thm]{{Question}}
\newtheorem{remark}[thm]{{Remark}}
\newtheorem{definition}[thm]{{Definition}}

\newcommand{\PP}{{\mathbb P}} \newcommand{\RR}{\mathbb{R}}
\newcommand{\QQ}{\mathbb{Q}} \newcommand{\CC}{{\mathbb C}}
\newcommand{\ZZ}{{\mathbb Z}} \newcommand{\NN}{{\mathbb N}}
\newcommand{\FF}{{\mathbb F}}

\newcommand{\acal}{\mathcal{A}}
\newcommand{\ccal}{\mathcal{C}}

\newcommand{\gcal}{\mathcal{G}}\newcommand{\hcal}{\mathcal{H}}

\newcommand{\lcal}{\mathcal{L}}

\newcommand{\ocal}{\mathcal{O}}
\newcommand{\qcal}{\mathcal{Q}}

\newcommand{\vcal}{\mathcal{V}}
\newcommand{\jcal}{\mathcal{J}}

\def\calQ{\qcal}  
 \def\calA{\acal} \def\calL{\lcal}
 \def\calH{\hcal} \def\H{\hcal}
 \def\calV{\vcal} 
\def\calC{\ccal} \def\calO{\ocal}
 \def\calJ{\jcal} \def\calG{\gcal}

\def\a{\alpha} \def\be{\beta}
\def\o{\omega} \def\O{\Omega}
\def\th{\theta} 
\def\vp{\varphi} \def\eps{\epsilon}

\newcommand{\dbar}{\bar\partial}
\newcommand{\ddbar}{\partial\dbar}
\newcommand{\diag}{{\operatorname{diag}}}

\def\Im{{\operatorname{Im}}}
\def\Re{{\operatorname{Re}}}

\def\sing{{\operatorname{sing}}}
\def\max{{\operatorname{max}}}

\def\phg{{\operatorname{phg}}}
\def\tr{\hbox{\rm tr}}
\def\Id{\hbox{\rm Id}}
\def\D{\Delta}

\def\MA{Monge--Amp\`ere } \def\MAno{Monge--Amp\`ere}

\def\K{K\"ahler } 
\def\KE{K\"ahler--Einstein } 
\def\KEE{K\"ahler--Einstein edge } 
\def\KRF{K\"ahler--Ricci flow }

\def\polishl{\char'40l}
\def\Kolodziej{Ko\polishl{}odziej} \def\Blocki{B\polishl{}ocki}
\def\Holder{H\"older }

\def\ovp{{\o_{\vp}}}
\def\on{\omega^n}
\def\ovpn{\omega_{\vp}^n}
\def\ovps{\omega_{\vp_s}}

\def\Ric{\hbox{\rm Ric}\,}

\def\oKE{\o_{\h{\sml KE}}}

\def\vpKE{\vp_{\h{\sml KE}}}
\def\ovpn{\o^n_{\vp}}
\def\gM{g_{\hbox{\sml M}}}

\def\h#1{\hbox{#1}}

\def\text{\textstyle}
\def\dis{\displaystyle}

\def\q{\quad} \def\qq{\qquad}

\def\b#1{\bar{#1}}
\def\MsmD{M\!\setminus\!  D}

\def\PSH{\mathrm{PSH}}

\def\ra{\rightarrow}
\def\pa{\partial}
\def\isom{\cong}
\def\w{\wedge}
\def\i{\sqrt{-1}}

\def\Aut{{\operatorname{Aut}}}
\def\aut{{\operatorname{aut}}}

\def\Ds{\calD_s^{0,\gamma}}
\def\De{\calD_e^{0,\gamma}}
\def\Dw{\calD_w^{0,\gamma}}
\let\s=\sigma
\def\calHo{\calH_\omega}
\def\ovp{{\o_\vp}}
\def\sm{\setminus}
\def\al{\alpha}
\def\be{\beta}

\newcommand{\calD}{\mathcal D}

\def\Aut{{\operatorname{Aut}}}

\def\Ds{\calD_s^{0,\gamma}}
\def\De{\calD_e^{0,\gamma}}
\def\Dw{\calD_w^{0,\gamma}}
\def\tildeDs{{\widetilde\calD}_s^{0,\gamma}}
\def\tildeDe{{\widetilde\calD}_e^{0,\gamma}}

\def\del{\partial}
\def\ze{\zeta}
\def\Herm{{\operatorname{Herm}}}
\def\Cech{$\check{\h{C}}$ech }

\def\la{\lambda}
\def\beq{\begin{equation}}
\def\eeq{\end{equation}}
\def\beqno{\begin{equation*}}
\def\eeqno{\end{equation*}}
\def\eaeq{\end{aligned}}
\def\baeq{\begin{aligned}}

\def\bpf{\begin{proof}}
\def\epf{\end{proof}}
\def\Rico{\Ric\,\!\o}

\font\itnotsosml=cmti7

\def\thhit#1{\hbox{${\hbox{#1}}^{\,{\hbox{\itnotsosml th}}}$}}

\def\opcit{\underbar{\phantom{aaaaa}}}
\def\III{\h{\rm I}}

\def\fo{f_\omega}
\def\ga{\gamma}
\def\Ho{\calHo}
\def\HO{\calH_\O}
\def\G{\Gamma}
\def\etavp{\eta_{\varphi}}
\def\etavpn{\eta^n_{\varphi}}
\def\etan{\eta^n}

\def\BC{\h{\bf BC}}
\def\End{\h{\rm End}}
\def\eaeq{\end{aligned}}
\def\baeq{\begin{aligned}}
\def\phgs{polyhomogeneous}
\chardef\inodot="10

\def\HamMD{\h{\rm Ham}(M,D,\o)}
\def\hamMD{\h{\rm ham}(M,D,\o)}

\def\gC{g_{\hbox{\sml C}}}

\def\ovptn{\o_{\vp_t}^n}

\def\BC{\hbox{\bf BC}}

\def\Ent{\hbox{Ent}}

\def\osc{\operatorname{osc}}
\def\Sla{S_\lambda}

\def\la{\lambda}
\def\ze{\zeta}

\graphicspath{%
    {converted_graphics/}% inserted by PCTeX
    {C:/Users/me/Desktop/Dropbox/tex/Drafts/KECRM/figures-Ramos/}% inserted by PCTeX
}
\begin{document}
\title{Smooth and singular K\"ahler--Einstein metrics}

\author{ Yanir A. Rubinstein
}

\maketitle

\vglue-0.3cm
\centerline{\it Dedicated to Eugenio Calabi on the occassion of his
9\thhit{0}
 birthday }

\smallskip
\smallskip

\begin{abstract}
Smooth K\"ahler--Einstein metrics have been studied for the past 80 years.
More recently, singular K\"ahler--Einstein metrics have emerged as
objects of intrinsic interest, both in differential and algebraic geometry,
as well as a powerful tool in better understanding their smooth counterparts.
This article is mostly a survey of some of these developments.

\end{abstract}

\tableofcontents

\section{Introduction}

The \KE (KE) equation is among the oldest fully nonlinear equations
in modern geometry. A wide array of tools have been developed or applied
towards its understanding, ranging from Riemannian geometry,
PDE, pluripotential theory, several complex variables, microlocal analysis,
algebraic geometry, probability, convex analysis, and more.
The interested reader is referred to the numerous existing
surveys on related topics
\cite{Aubinbook,Biquard,Bl2007,Bl,Bourguignon1997,Bourguignon2003,CheltsovShramov,
Donaldson-toric,Futaki,Futaki2005,FutakiMabuchiSakane,
FutakiOno,FutakiSano,Gauduchon,Guedj,PhongSturmSurvey,
Siu,PhongSongSturm,Szekelyhidi,Thomas,Tianbook,Tian2002,Tian2012}.

In this article, which is largely a survey,
we do not attempt
a comprehensive overview, but instead focus on a rather subjective
bird's eye view of the subject.
Mainly, we aim to survey some recent developments on the study
of {\it singular} \KE metrics, more specifically, of {\it \KEE (KEE)} metrics,
and while doing so to provide a
unified introduction also to smooth \KE metrics. In particular, the new tools
that are needed to construct KEE metrics give a new perspective on construction
of smooth KE metrics, and so it seemed worthwhile to use this opportunity to
survey the existence and non-existence theory of both smooth and singular KE
metrics.
We also put an emphasis on understanding concrete new geometries that
can be constructed using the general theorems, concentrating on the
complex 2-dimensional picture, and touch on some relations to algebraic
geometry and non-compact Calabi--Yau spaces.

As noted above, most of the material in this article is a survey of existing results
and techniques, although, of course, sometimes the presentation may be different
than in the original sources.
For the reader's convenience, let us also point out the sections that contain
results that were not published elsewhere. First, some of 
the treatment of energy functionals in
\S\ref{StaircaseSubSec} and \S\ref{EquivalenceSubSec} extends some 
previous results of the author from the smooth setting to the edge setting, although
the generalization is quite straightforward and is presented here for the sake
of unity. The section on Bott--Chern forms \S\ref{BottChernSubSec} 
is essentially taken from the author's thesis, again with minor modifications to
the singular setting.  Second, some of the treatment of the a priori estimates
in \S\ref{APrioriSubSec} is new. In particular, the reverse Chern--Lu inequality
introduced in \S\ref{ReverseCLSubSec} is new, as is the proof in \S\ref{AYAsCorSubSec} 
that the classical Aubin--Yau Laplacian estimate follows from it.
Finally, we note the minor observation made in \S\ref{UniformityMetricTwoSubSubSec} 
that some of the classical Laplacian estimates
can be phrased also for \K metrics that need not satisfy a complex \MA equation.
Third, Proposition \ref{NegativeCaseProp} is a very slight extension
of the original result of Di Cerbo.

Many worthy vistas are not surveyed here, including:
toric geometry; quantization; log canonical thresholds,
Tian's $\alpha$-invariant, and Nadel multiplier ideal sheaves
\cite{CheltsovShramov,Tian1987,Nadel}; noncompact \KE metrics;
Berman's probabilistic approach to the KE problem
\cite{Berman-prob}; 
K\"ahler--Einstein metrics on singular varieties,
as, e.g., in \cite{EGZ,BermanG};
recent work on algebraic obstruction
to deforming singular KE metrics to smooth ones \cite{CDSun,Tian2013};
the GIT related aspects of the \KE problem, that we do not discuss in any detail
in this article, instead referring to the survey of Thomas \cite{Thomas}
and the articles of Paul \cite{Paul2012,Paul2013}.

\subsection*{Organization}
Historically, the \KE equation was first phrased locally as a complex
\MA equation.
This, and the corresponding global formulation, are discussed in Section \ref{KEEqSec}.
Section \ref{KahlerSection} is a condensed introduction to \K edge geometry,
describing aspects of the theory that are absent from the study of smooth
\K manifolds: new function spaces, a different notion of smoothness (polyhomogeneity),
a theory of (partial) elliptic regularity,
and new features of the reference geometry (e.g., unbounded curvature). We use these tools to describe the
structure of the Green kernel of \K edge metrics, and the
proof of higher regularity of KEE metrics \cite{JMR}.
Section \ref{Sec4} summarizes the main existence and non-existence results on KE(E)
metrics. First, it describes the classical obstructions due to Futaki and Mastushima coming
from the Lie algebra of holomorphic vector fields (and their edge counterparts),
and their more recent generalization to various notions of ``degenerations", that capture
more of the complexity of Mabuchi K-energy, the functional underlying the KE problem.
Second, it states the main existence and regularity theorems for KE(E) metrics.
The strongest form of these results appears in \cite{JMR,MR2012}, generalizing
the classical results of Aubin, Tian, and Yau from the smooth setting,
and those of Troyanov from the conical Riemann surface setting,  to the edge
setting. We also describe other approaches to this problem.
The main
objective of \S\ref{RCMSec}--\S\ref{APrioriSubSec} is to describe the analytic tools
needed to carry out the proof of these theorems, in a unified manner,
both in the smooth and the edge settings, and regardless of the sign of the curvature.
We describe this in more detail below; before that, however, Section \ref{EnergySection} reviews the variational theory underlying
the \MA equation in this setting. In particular, it reviews the basic properties
of the Mabuchi K-energy and the functionals introduced by Aubin, Mabuchi, and others.
The alternative definition of these via Bott--Chern forms is described in detail. A relation
to the Legendre transform due to Berman is also described. Finally, we describe
the properness and coercivity properties of these functionals. The former is needed for
the actual statement of the existence theorem in the positive case from \S\ref{Sec4}.

Section \ref{RCMSec} describes a new approach, the {\it Ricci continuity method},
 developed by Mazzeo and the author in \cite{JMR} to
prove existence of KE(E) metrics in a unified manner, and with only one-sided
curvature bounds. This approach works in a unified manner for all cases
(negative, zero, and positive Einstein constant) and for both
smooth and singular metrics.
Section \ref{APrioriSubSec} describes the a priori estimates needed to
carry out the Ricci continuity method. In doing so, this section
also provides a unified reference for a priori estimates for a wide
class of singular \MA equations.
Section \ref{AsympSection} describes work of Cheltsov and the author \cite{CR}
on classification problems in algebraic geometry related to the KEE
problem. Here, the notion of {\it asymptotically log Fano varieties}
is introduced. This is a class of varieties  much larger than Fano
varieties, but where we believe there is still hope of classification and
a complete picture of existence and non-existence of KEE metrics.
Section \ref{logCalabiSec} then builds on these classification
results to phrase a logarithmic version of Calabi's conjecture.
We then briefly describe a program to relate this conjecture concerning KEE metrics to
Calabi--Yau fibrations and global log canonical thresholds. Finally, some progress towards this conjecture
is described, in particular giving new KEE metrics on some explicit pairs,
and proving non-existence on others.

\section{The \KE equation}
\label{KEEqSec}

A \K manifold is a complex manifold $(M,J)$ equipped with a closed
positive 2-form $\o$ that is $J$-invariant, namely
$\o(J\,\cdot\,,J\,\cdot\,)=\o(\,\cdot\,,\,\cdot\,)$.
Here $M$ is a differentiable manifold, and
$J$ is a complex structure on $M$, namely
an endomorphism of the tangent bundle $TM$
satisfying $J^2=-I$ and $[T^{1,0}M,T^{1,0}M]\subseteq T^{1,0}M$,
where $TM\otimes_\RR\CC=T^{1,0}M\oplus T^{0,1}M$, with
$T^{1,0}M$ the $\i$-eigenspace of $J$, and
$T^{0,1}M$ the $-\i$-eigenspace of $J$.
Associated to $(M,J,\o)$ is a Riemannian metric
$g$, defined by $g(\,\cdot\,,\,\cdot\,)=\o(\,\cdot\,,J\,\cdot\,)$,
whose Levi-Civita connection we denote by $\nabla$.
Equivalently, one may also define $(M,J,g)$ to be a \K manifold when
$g$ is $J$-invariant, $J^2=-I$, and $\nabla J=0$
(i.e., $g$-parallel transport preserves $T^{1,0}M$).

Schouten called this new type of geometry ``unitary," 
\cite{Schouten} and this can be understood from the fact
that not only does one obtain a reduction of the structure group to $U(n)$
(this merely characterizes almost-Hermitian manifolds)
but also the holonomy is reduced from $O(2n)$ to $U(n)$.
In retrospect this name seems
quite fitting, on a somewhat similar footing
to ``symplectic geometry," a name first suggested by Ehresmann.
However, only a few authors in the 1940's and 1950's were familiar with Schouten-van Dantzig's work.
The first few accounts of Hermitian geometry, mainly by Bochner, Eckmann and Guckenheimer
referred only to K\"ahler's works, and the name ``\K geometry" became rooted
(for more on this topic, see \cite[\S2.1.4]{RThesis}).

One of Schouten--van Dantzig's and K\"ahler's discoveries was that on the eponymous manifold
the Einstein equation is implied, locally,
by the complex \MA equation
\begin{equation}
\label{LocalMAEq}
\det[u_{i\b j}]:=\det \bigg[\frac{\pa^2 u}{\pa z^i \pa{\bar z}^j}\bigg]= e^{-\mu u},\q \h{on\ } U,
\end{equation}
where $\mu\in\RR$ is the Einstein constant \cite{Kahler,SchoutenVanDantzig}.
Here $u$ is a plurisubharmonic
function defined on some coordinate chart $U\subset M$, with holomorphic coordinates
$z_1,\ldots,z_n:U\ra\CC^n$; Indeed, if we let
$g_{\h{\sml H}}|_U=\i u_{i\b j}dz^i\otimes \overline{dz^j}$
denote the Hermitian metric associated to $u$ (with summation over repeated indices), then the Riemannian metric
$g=2\Re\, g_{\h{\sml H}}$ is Einstein  if \eqref{LocalMAEq} holds because
$\Ric\, g = 2\Im((\log\det[u_{i\b j}])_{k\b l}dz^k\otimes\overline{dz^l})$.
Thus, in the ``unitary gauge", the Einstein system of equations reduces to a single,
albeit fully nonlinear, equation.

For some time though, it was not clear how to patch up these equations
globally. The crucial observation is that both sides of \eqref{LocalMAEq}
are the local expressions of global Hermitian metrics on two different
line bundles: the
anticanonical line bundle $\Lambda^n T^{1,0}M$ on the left hand side,
and the line bundle associated to $\mu[\o]$ on the right.
In particular this implies
that these line bundles must be isomorphic, and
\begin{equation}
\label{c1Eq}
c_1(M,J)-\mu[\o]=0.
\end{equation}
This of course puts a serious cohomological restriction on the problem.
But assuming this, the local \MA equation can be converted to a global one.
Let ${\bf dz}:=dz^1\w\cdots \w dz^n$. The expressions
$$
\i^{n^2}\det[u_{i\b j}]{\bf dz}\w\overline{\bf dz}
= e^{-\mu u} \i^{n^2}{\bf dz}\w\overline{\bf dz}
$$
represent globally defined volume forms on $M$, appropriately interpreted.
We choose a representative \K form $\o$ of $[\o]$. Since
$\i\ddbar u$ is the curvature of $e^{-u}$ it lies in $[\o]$.
On each $U$, $\i\ddbar u = \o+\i\ddbar\vp$ for a globally defined $\vp\in C^\infty(M)$:
this follows from the Hodge identities \cite[p. 111]{GH} by setting
$\vp=\tr_\o G_\o (\i\ddbar u - \o)$, where $G_\o$ denotes the Green
operator associated to the Laplacian (on the exterior algebra).

Thus,
$$
(\o+\i\ddbar\vp)^n
= \o^n e^{-\log\det[v_{i\b j}]-\mu v-\mu\vp}
$$
where locally $\o=\i v_{i\b j}dz^i\w \overline{dz^j}$.
But by \eqref{c1Eq} $e^{-\log\det[v_{i\b j}]-\mu v}$ is the quotient
of two Hermitian metrics on the same bundle, hence a globally defined positive
function that we denote by $e^{f_\o}$.
We thus obtain the \KE equation,
\begin{equation}
\label{KEEq}
\ovpn=\on e^{\fo-\mu\vp}, \quad \h{ on } M
\end{equation}
for a global smooth function $\vp$ (called the \K potential of $\ovp$ relative to $\o$).
The function $\fo$, in turn, is given in terms of the reference geometry and is thus known. It
is called the Ricci potential of $\o$,
and satisfies $\i\ddbar\fo=\Ric\,\o-\mu\o$, where it is
convenient to require the normalization $\int e^{\fo}\on=\int \on$.
Equivalently, \eqref{KEEq} says that $f_{\ovp}=0$. Thus, in the language
of the introduction to \S\ref{Sec4}, the KE problem
has a solution precisely when  the vector field $f:\phi\mapsto f_{\o_\phi}$ on the space
of all \K potentials has a zero.

\section{\K edge geometry}
\label{KahlerSection}

This section introduces the basics of \K edge geometry. First, we
describe some general motivation for introducing these more general
geometries in \S\ref{KahlerEdgeIntroSubsec}.
The one-dimensional geometry of a cone is the topic of \S\ref{ConSubSec}.
It is fundamental, since a \K edge manifold looks like a cone transverse
to its `boundary' divisor. Subsection \ref{KEEEqSubSec} phrases
the KEE problem as a singular complex \MA equation, generalizing
the discussion in \S\ref{KEEqSec}.
What is the appropriate smooth structure on a \K edge space?
This is discussed in \S\ref{phgSubSec} where the notion
of polyhomogeneity is introduced.
The edge and wedge scales of \Holder function space
are defined in \S\ref{EdgeWedgeSubSec}.
Next, the various
\Holder domains relevant to the study of the {\it complex} \MA
equation are introducted in \S\ref{HolderSubSec}.

\subsection{A generalization of \K geometry}
\label{KahlerEdgeIntroSubsec}

The cohomological obstruction
\eqref{c1Eq} is necessary for the local KE geometries to
patch up to a global one. However, the condition \eqref{c1Eq}
is very restrictive. Is there a way of constructing a KE
metric at least on a large subset of a general \K manifold? Of course,
such a question makes sense and is very interesting also
in the Riemannian context.
The interpretation of the local KE equation in terms of
Hermitian metrics on line bundles, though, distinguishes between
these two settings.

Thus, suppose that $c_1(M,J)-\mu[\o]$ is not
zero (nor torsion) but that this difference, or `excess
curvature' can be
decomposed as follows: there exist divisors $D_1,\ldots,D_N$
and numbers $\be_i>0$ such that
%\begin{equation}
%\label{cohomEq}
$$c_1(M,J)-\mu[\o]
=\sum_{i=1}^N (1-\be_i)c_1(L_{D_i}),
%\end{equation}
$$
with $L_F$ denoting the line bundle associated to a divisor $F$.
At least when $M$ is projective, by the Lefschetz theorem on (1,1)-classes \cite[p. 163]{GH}
this can always be done when the left hand side belongs
to $H^2(M,\ZZ)\hookrightarrow H^2(M,\RR)$, and therefore
also when it is merely a rational class, or a real class.
Of course this
means that we might need to take a non-reduced divisor on
the right hand side, and limits the usefulness of such a generalization
to the case when the divisor $D=\sum D_i$ is not too singular,
and the $\be_i$ are not too large.

Thus, necessarily any KE potential $u$ must satisfy locally,
\begin{equation}
\label{LocalMAEEq}
\det[u_{i\b j}]= e^{-\mu u}\prod_{i=1}^N|e_i|^{2\be_i-2},\q \h{on\ } U,
\end{equation}
where $e_i$ denotes a local holomorphic section of $D_i$.
This equation is then quite similar to the local KE equation
\eqref{LocalMAEq}
on neighborhoods $U$ contained in the complement of the $D_i$.
However, if $U$ intersects any of the $D_i$ nontrivially, the
equation becomes singular or degenerate.

To understand this equation better, it is helpful to consider the
model case with $N=1$, $M=\CC^n$ and $D=\{e_1=z_1=0\}=\{0\}\times\CC^{n-1}$.
Then,
$u=\frac12(\frac1{\be^2}|z_1|^{2\be}+|z_2|^2+\ldots+|z_n|^2)$ is a solution of \eqref{LocalMAEEq} with $\mu=0$.
Note that $u$ corresponds to a singular, but continuous, Hermitian metric $e^{-u}$.
We call the associated curvature form the {\it model edge form
on $(\CC^n,\CC^{n-1})$,}
\beq
\label{obetaEq}
\omega_\beta = 
-\i\ddbar\log e^{-u}
=
\frac12 \sqrt{-1}
\Big(
|z_1|^{2\beta-2} dz^1 \w \overline{dz^1} + \sum_{j=2}^n dz^j \w \overline{dz^j}
\Big),
\eeq
and also denote by
\begin{equation}
\label{modelmetric}
g_\beta  =   |z_1|^{2\beta-2} |dz^1|^2 + \sum_{j=2}^n |dz^j|^2,
\end{equation}
the {\it model edge metric on $(\CC^n,\CC^{n-1})$}.
More generally, if
$D=D_1+\ldots+D_N$ is the union of $N$ coordinate hyperplanes in $\CC^n$ which intersect
simply and normally at the origin,
we set
$\be=(\be_1,\ldots,\be_N)$ and denote
the {\it model edge form on $(\CC^n,D)$} by
$$
\omega_\beta =
\frac12 \sqrt{-1}
\Big(
\sum_{i=1}^N\la_i|z_{\eps(i)}|^{2\beta_i-2}
dz^{\eps(i)} \w \overline{dz^{\eps(i)}}
+
\sum_{j\in\{1,\ldots,n\}\sm\{\eps(i)\}_{i=1}^N} dz^j \w \overline{dz^j}
\Big),
$$
where
$\la_i=1$ if $0\in D_i$ and $\la_i=0$ otherwise,
and where $\eps(i)=0$ if $\la_i=0$, and otherwise
$\eps(i)\in\{1,\ldots,n\}$ is such that $\{z_{\eps(i)}=0\}=D_i$.

This model case motivates the following generalization of
a \K manifold. Some visual references are given in Figure
\ref{FigRiemannSurfaces} in \S\ref{RSSubSec} below.

\begin{definition}
A \K edge manifold is a quadruple $(M,D,\be=(\be_1,\ldots,\be_N),\o)$, with
$M$ a smooth K\"ahler manifold, $D=D_1+\ldots D_N\subset M$ a simple normal crossing (snc) divisor,
$\be_i:D_i\ra\RR_+$ a function for each $i=1,\ldots,N$, and, finally, $\o$ a \K current on $M$
that is smooth on
$M\sm D$ and asymptotically equivalent to $\o_{\be(p)}$
near each  point $p$ of $D$.

\end{definition}
For the notion of a K\"ahler current (also called a positive $(1,1)$ current)
we refer to \cite[Chapter 3]{GH}.
One could make the definition more general, e.g., by allowing $D$ to be more
singular, but we do not explore that here. Furthermore,
in our discussion below, we will always assume that $\be_i$ is constant
on each component $D_i$. (This assumption is present in essentially
all works on the subject so far.)
Lastly, for the moment, we are deliberately vague on the meaning of ``asymptotically equivalent."
Several working definitions are given in \S\ref{KedgeSubSec} (see in particular Lemma
\ref{KedgeLemma}).

The study of \K edge metrics was initiated by Tian \cite{Tian1994}
motivated in part by the possibility of endowing more \K manifolds
with a generalized KE metric, when the obstruction \eqref{c1Eq}
does not vanish. Of course, the possibility of uniformizing more
\K manifolds is exciting in itself, and there are many possible
applications of such metrics in algebraic geometry (see, e.g., \cite{Tsuji,Tian1994}).
However, as we will see later, this generalization sheds considerable
light also on the theory of smooth KE metrics. 
Finally, \K edge manifolds are also a natural generalization of conical Riemann
surfaces, who were first systematically studied by Troyanov \cite{Troyanov};
see \S\ref{RSSubSec} for more on this topic.

\subsection{The geometry of a cone}
\label{ConSubSec}

In the previous subsection we arrived at a generalization of
\K geometry by seeking a generalization of the \KE equation.
Before going back to the latter, let us first say a few words on
the geometry described by \eqref{obetaEq}.

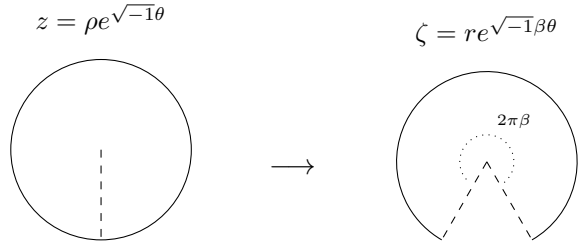
\begin{figure}
\centering
	\begin{tikzpicture}[scale=1.2]
		\draw (0,0) circle (1);
		\draw [dashed] (0,0) -- (0,-1);
		\node at (0,1.45) {\small $z=\rho e^{\i\theta}$}; %,\; \theta\in[0,2\pi)$};
	\end{tikzpicture}
	\ \ \ \ \ \ \ \raise25pt\hbox{$\longrightarrow$}\ \ \ \ \ \ \
	\begin{tikzpicture}[scale=1.2]
		\draw (0,0)+(-60:1) arc (-60:240:1);
		\draw [dashed] (0,0) -- (-60:1);
		\draw [dashed] (0,0) -- (240:1);
		\draw [dotted] (0,0)+(-50:.3) arc (-50:230:.3);
		\node at (.3,.45) {\tiny $2\pi\beta$};
		\node at (0,1.45) {\small $\zeta=r e^{\i\be\theta}$}; %,\; \theta\in[0,2\pi)$};
	\end{tikzpicture}
	\caption{The map $z\mapsto z^\beta=:\zeta$ maps the disc to a wedge of angle $2\pi\beta$
(the inverse map identifies the two `sides' of the wedge).
This map pulls-back the flat metric $|d\zeta|^2$ on the wedge to the singular metric
$\be^2|z|^{2\be-2}|dz|^2$ on the disc.
}
\label{FigureCone}
\end{figure}

The basic observation is that
\beq\label{CbetaEq}
C_\be:=(\CC,|z|^{2\be-2}\i dz\w\overline{dz})
\eeq
is a cone with tip at the origin. For instance,
when $\be=1/k$ with
$k\in\NN$, we obtain an orbifold.
We emphasize that here we equip the smooth manifold $\CC$ with
a singular metric, and we are claiming that thus equipped it is
isometric to a singular metric space, a cone.
Perhaps the simplest way to see
this is to recall the construction of a cone out of a wedge. Starting
with a wedge of angle $2\pi\be$ one identifies the two
sides of the boundary. This can of course be done in many ways,
in other words there are many maps of the form
$(r,\th)\mapsto(f(r,\th),\th/\be)$ (so the angle in the target
indeed varies between $0$ and $2\pi$, or in other words, the target
indeed `closes up' and is homeomorphic to $\RR_+\times S^1(2\pi)$).
However, there is a unique
such map that preserves angles, or in other words is holomorphic.
By the Cauchy--Riemann equations
(i.e., 
$$
r\del_r (\Re F,\Im F)=\del_\th(\Im F,-\Re F)
$$ 
for a holomorphic function $F$)
 it must be of the form $f(r,\th)=Cr^{1/\be}$, for some constant $C$.
The inverse of this map, from the cone to the wedge,
is given by $(\rho,\phi)\mapsto(\rho^\be,\be\phi)$, or
simply $z=\rho e^{\i\phi}\mapsto z^\be=:\ze$ (see Figure \ref{FigureCone}).
Declaring this map to be an isometry determines the geometry of
the cone; pulling back the Euclidean metric $|d\zeta|^2=\i\ddbar|\ze|^2$
endows the cone with the metric $\i\ddbar|z|^{2\be}
=\be^2|z|^{2\be-2}\i dz\w\overline{dz}$.

\subsection{The \KE edge equation}
\label{KEEEqSubSec}

We now return to our discussion of
the generalized {\it local} \KE equation \eqref{LocalMAEEq}.

To turn \eqref{LocalMAEEq} into a global equation we seek, at least
formally (in the case we are not dealing
with $\QQ$-line bundles),
an equality of two continuous Hermitian metrics: one on $\mu\O$
and the other on $-K_M+\sum_i(1-\be_i)L_{D_i}$.
In analogy with the discussion of \S\ref{KEEqSec}, we now choose
a reference metric that is locally asymptotic to the model edge
geometry, instead of the Euclidean geometry on $\CC^n$ that
was implicitly the model there. Suppose that $h_i$ is a smooth Hermitian
metric on $L_{D_i}$, and that $s_i$ is a global holomorphic section of
$L_{D_i}$, so that $D_i=s_i^{-1}(0)$.
We let
\beq
\label{oEq}
\o:=\o_0+\i\ddbar\phi_0,
\eeq
with
\beq
\label{phi0Eq}
\phi_0:=c\sum(|s_i|^2_{h_i})^{\be_i}.
\eeq
An easy computation shows that $\o$ is locally equivalent to $\o_\be$
and that for small enough $c>0$ it defines a \K metric away from $D$
\cite[Lemma 2.2]{JMR}. Moreover, in the coordinates above,
near a smooth point of $D_1$, for example,
$|z_1|^{2\be_1-2}\det \psi_{i\b j}$ is continuous, as desired.
Similar properties hold near crossing points of $D$.
Thus, we may proceed exactly as in \S\ref{KEEqSec} to obtain
an (seemingly!) identical equation,
\begin{equation}
\label{KEEEq}
\ovpn=\on e^{\fo-\mu\vp}, \quad \h{ on } M\sm D,
\end{equation}
only now away from $D$, and with respect to the reference form $\o$, where
$\vp$ is required to lie in the space of \K edge potentials
\begin{multline}
\label{KEdgePotentialsSpace}
\calH_{\o} := \{\vp\in C^\infty(\MsmD)\cap C^0(M)
\, :\,   \o_{\vp} := \o+\i\ddbar\vp>0 \h{\ on\ } M  \\
\mbox{and $\o_\vp$ is asymptotically equivalent to}\  \omega_\be\}.
\end{multline}
Equation \eqref{KEEEq} is the {\it \KE edge (KEE) equation}.

By construction, the twisted Ricci potential of $\o$,
$f_\o$, satisfies
\beq\label{foEq}
\i\ddbar\fo=\Ric\,\o-(1-\be)[D]-\mu\o,
\eeq
where it is
again convenient to require the normalization $\int e^{\fo}\on=\int \on$.
Differentiating \eqref{KEEEq} leads to the following equivalent formulation
of the KEE equation.

\begin{definition}
\label{KEEDef}
With all notation as above, a K\"ahler current $\oKE$ is called a
K\"ahler--Einstein edge current with angle $2\pi \beta$ along $D$ and Ricci
curvature $\mu$ if $\oKE\in\calH_\o$ (see (\ref{KEdgePotentialsSpace})), and if
\begin{equation}
\label{EdgeKEEq}
\Ric \oKE - (1-\beta) [D]= \mu \oKE,
\end{equation}
where $[D]$ is the current associated to integration along $D$,
and
where $\Ric\o$ denotes the Ricci current (on $M$) associated to $\o$, namely, in local coordinates
$\Ric\o=-\i\ddbar\log\det[g_{i\bar j}]$ if $\o=\i g_{i\bar j} dz^i\wedge \overline{dz^j}$.

\label{Kcurrent}
\end{definition}

The KEE equation may also be rewritten in terms of the background
smooth geometry. We carry this out for pedagogical purposes,
since it allows to unravel the difference between
equations \eqref{KEEEq} and \eqref{KEEq},
that are formally the same but involve different geometric objects.

 Let $e$ be a local holomorphic frame for $L_D$ valid in a neighborhood intersecting $D$, such that $s=z_1e$
on that neighborhood, and denote by
$$
a_i:=|e_i|^2_{h_i}
$$
a smooth positive function on that neighborhood.
Define $F_{\o_0}$ (up to a constant, for the moment) by
$\i\ddbar F_{\o_0} = \Ric\o_0-\mu\o_0+\sum_i(1-\beta_i)\i\ddbar\log a_i$.
Setting $\tilde\vp:=\phi_0+\vp$,
it easy to see that
\begin{equation}
\label
{fomegaSecondEq}
(\o_0+\i\ddbar\tilde\vp)^n
=
\prod_i|s_i|_{h_i}^{2\be_i-2}\o_0^n  e^{F_{\o_0}-\mu\tilde\vp},
\end{equation}
where we think of this equation as determining a normalization of $F_{\o_0}$
such that both sides have equal integrals   (cf. \cite[(53)]{JMR}).

For much of the rest of this article we will be concerned with
solving this equation with optimal estimates on the solution.
In this regard, we note that global solutions, smooth away from $D$, of this
(when $\mu\le0$) and
quite more general \MA equations were constructed already
in much earlier work of Yau \cite[\S7--9]{Yau1978}.
Our goal though, is to explain how to go beyond such
weak solutions and obtain estimates that show the metric
in fact has edge singularities near the `boundary' $D$.
To start, we define the appropriate function spaces
to prove such estimates.

\begin{remark}
{\rm
A point we glossed over in the discussion is the
fact that a necessary condition for the existence
of a KEE metric is that the first Chern class of the adjoint bundle,
\beq
\label{cohomEq}
-K_M+\sum_i(1-\be_i)L_{D_i},
\h{\ \ is\ $\mu$ times an ample class}.
\eeq
This is
obvious in the smooth case ($\be_i=1$)
from the KE equation. In general, the KEE equation
only guarantees that this class is represented
by ($\mu$ times) a \K current. But the edge singularities
are sufficiently mild that the divisor $D$ has zero
volume \cite[\S6]{JMR}, and the
Lelong numbers of this current along $D$ are zero.
As observed by Di Cerbo \cite{DiCerbo2012} this implies
that the class is actually ($\mu$ times) an ample class, i.e., represented
by a smooth \K metric.
}
\end{remark}

\subsection{Three smooth structures, one polyhomogeneous structure}
\label{phgSubSec}

Naturally, if one's goal is to obtain existence and regularity of
certain geometric objects (such as KEE metrics) on a \K edge manifold,
then understanding the differentiable structure underlying such a space
should play a key r\^ole.

First, let us consider the simplest compact closed example, that of $M=\PP^1$ with a
single cone point $D=p$. Then $M$ of course has its natural
conformal structure, and there is a corresponding smooth structure.
The former is represented locally near $p$ by the holomorphic coordinate
$z$, and the latter can be represented by the associated polar coordinates
$(\rho,\theta)$, where of course all the points $(0,\theta)$ are identified
with $p$.

Recall from \S\ref{ConSubSec} that
an alternative conformal structure is represented locally by the coordinate
$\zeta:=z^\be$. Of course, $\zeta$ is multi-valued, and one must choose a
branch of the Riemann surface associated to $z\mapsto z^\be$. Whenever
we work with $\zeta$, we assume such a choice was made. More specifically
we slit the disc, and work with the associated polar coordinates denoted by $(r,\theta)$,
where now $\theta\in[0,2\pi\be)$, and these represent the associated smooth
structure. That is, smooth functions are smooth functions of $r,\theta$.
Clearly, this smooth structure is incompatible with the one
of the previous paragraph.

Which smooth structure should we work with? The point of view
stressed in \cite{JMR}, following Melrose's general framework
\cite{Melrose}, is that it is rather the {\it polyhomogeneous
structure} that is central to the problem, and not either
of these smooth structures. The purpose of this section is to explain what
is meant by this. To carry this out, it will be preferable to work with
a bona fide {\it singular space} associated to $M$ (recall $M$ itself
is smooth).
That is, we consider the singular metric space $M_\sing$ obtained as the metric
completion from the underlying distance
function associated to $M$ equipped with its \K edge metric.
It is readily seen that this space is homeomorphic to $M$ (since
the divisor is at finite distance from any point by \eqref{modelmetric}).
What is important
to the development of the theory described below is that the smooth locus of this
singular space $M_\sing$ is precisely the complement of a simple normal crossing divisor in
a nonsingular space. Thus, as we describe below, it will admit an {\it edge structure}.

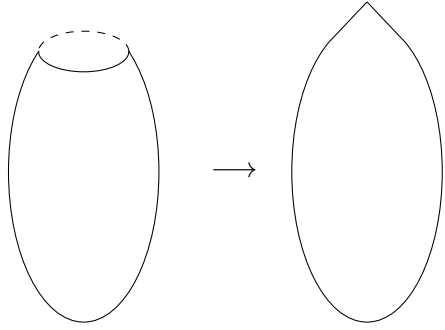
\begin{figure}
\centering
	\begin{tikzpicture}
		\begin{scope}
			\clip (-1.2,-2.2) rectangle (1.2,1.6);
			\draw (0,0) ellipse (1 and 2);
			\draw (0,1.6) ellipse (.6 and .27);
		\end{scope}
		\begin{scope}
			\clip (-1.2,1.6) rectangle (1.2,2.2);
			\draw [dashed] (0,1.6) ellipse (.6 and .27);
		\end{scope}
		\node at (1.5,0) {\hglue1cm$\longrightarrow$};
	\end{tikzpicture}
	\begin{tikzpicture}
		\begin{scope}
			\clip (-1.2,-2.2) rectangle (1.2,1.8);
			\draw (0,0) ellipse (1 and 2);
%			\draw (0,1.6) ellipse (.6 and .27);
		\end{scope}
		\begin{scope}
			\clip (-1.2,1.6) rectangle (1.2,2.26);
			\draw (-.43,1.805) -- (0,2.26);
			\draw (0,2.26) -- (.43,1.805);
		\end{scope}
	\end{tikzpicture}
\caption{The real blow-up of the tear-drop space}
\label{RealBlowUpFig}
\end{figure}

The first step we take is actually to {\it desingularize} $M_\sing$, i.e., resolve its singularities to obtain
a manifold with boundary. The edge structure we will define shortly, will, by definition, be on this desingularized space.
The desingularization is done, as usual, by a series of resolutions by real blow-ups;
when $D$ is smooth and connected a single such blow-up suffices.
In the simplest example of $(\PP^1,p)$ considered just above
this amounts to a single real blow-up at $p$, resulting in
the manifold with boundary $X$ consisting of the disjoint union
of $\PP^1\sm\{p\}$ and $S^1$ (see Figure \ref{RealBlowUpFig}).
In other words, the smallest manifold with boundary on which
the polar coordinates are well-defined, without identifying
the points $\{(0,\theta)\,:\,\theta\in[0,2\pi)\}$.
In general,
the manifold $X$ is the real blow-up of $M_\sing$ at $D$,
i.e., the disjoint union $M \setminus D$ and the circle normal bundle of $D$ in $M$,
endowed with the unique smallest topological and differential structure so that
the lifts of smooth functions on $M$ and polar coordinates around $D$ are smooth.

The advantage of working with functions on $X$ rather than on $M$ or $M_\sing$ is convenience.
For instance, it is much easier to keep track of singularities of distributions
(such
as the Green kernel) on the desingularization--this is explained in detail in
\S\ref{GreenSubSec}. When $D$ has crossings this desingularization also serves
to give a reasonable description for distributions near crossing points \cite{MR2012,MR}.

The relevant
coordinates on $X$ are $(r,\theta,y)$ where $y=(y_1,\ldots,y_{2n-2})$ denotes
the `conormal' coordinates, i.e., coordinates on $D$.
The relevant smooth structure to our problem turns out to be
the smooth structure of $X$ as a manifold with boundary.
Namely, smooth functions are smooth functions of $r,\theta,$ and $y_1,\ldots,y_{2n-2}$.

Next, we defined the associated `edge structure' of $X$ in the sense of Mazzeo \cite[\S3]{Mazzeo}.

\begin{definition}
The edge structure on the smooth compact manifold with boundary $X$ is
the Lie algebra of vector fields $\calV_e(X)$ generated by the smooth vector fields
on $X$ that are tangent to the $S^1$ fibers on $\del X$. That is,
$$
\calV_e(X)=\h{\rm span}_\RR\Big\{
r\frac{\del}{\del r}, \frac{\del}{\del \th}, r\frac{\del}{\del y_1},\ldots,
r\frac{\del}{\del y_{2n-2}}\Big\}.
$$
\end{definition}

The space $C^\infty_e$ is defined to be the space of $C^\infty$ functions on $X$
with respect to this edge structure, i.e., functions that are infinitely-differentiable
with respect to the vector fields in $\calV_e(X)$; see \S\ref{EdgeWedgeSubSec} for
an alternative, somewhat more geometric, definition of this space.
Note that a function on $X$ can be pushed-forward, using the blow-down
map, to a function on $M\sm D$, and, by abuse of notation, we often make this identification
without mention. Much care is needed here, however, since
a function in  $C^\infty_e$ need not correspond even to a continuous
function on $M$ (consider, e.g., the function $\sin(\log r)$)!

In \K edge geometry there is a clear distinction between
differentiation in the direction normal to the edge
and in the complementary directions. Thus, it is convenient
to define the Lie algebra $\calV_b(X)$ of all smooth vector
fields on $X$ that are tangent to $\del X$, i.e.,
$$
\calV_b(X)=\h{\rm span}_\RR\Big\{
r\frac{\del}{\del r}, \frac{\del}{\del \th}, \frac{\del}{\del y_1},\ldots,
\frac{\del}{\del y_{2n-2}}\Big\},
$$
and introduce the following terminology.

\begin{definition}
The space of {\it bounded conormal functions on X} is
$$
\calA^0=\calA^0(X):=\{ f\in L^\infty(X)\,:\,
\h{for all $k\in\NN$ and $V_1,\ldots,V_k\in\calV_b(X)$,\ }
V_1\cdots V_kf\in L^\infty(X)\}.
$$
\end{definition}

 In other words, these are bounded functions on $X$ that are infinitely
differentiable with respect to vector fields in $\calV_b(X)$.
Thus, they are infinitely-differentiable in directions tangent to $D$ (the `conormal
directions'), but still potentially
rather badly behaved with respect to the vector field $\frac{\del}{\del r}$.
Another geometric description of $\calA^0$ is given in \S\ref{EdgeWedgeSubSec}
below.

Now,
let us finally define the {\it polyhomogeneous structure} of $X$.
For that we first set
$$
\calA^{\gamma,p}:=r^\gamma(\log r)^p\calA^0.
$$

\begin{definition}
\label{phgDef}
The space of {\it bounded polyhomogeneous functions on X} is
$$
\calA_\phg^0:=\Big\{ f\in \calA^0(X)\,:\,
f \sim \sum_{j=0}^\infty \sum_{p=0}^{N_j} a_{jp}(\theta, y) r^{\sigma_j} (\log r)^p\Big\},
$$
with $\Re\,\sigma_j$ increasing to $\infty$, and  $\Re\,\sigma_j\ge0$ with $N_j=0$ if $\Re\,\sigma_j=0$.
\end{definition}
\noindent
Similarly, one defines $\calA_\phg^{\gamma,p}:=r^\gamma(\log r)^p\calA^0_\phg$,
and the {\it polyhomogeneous structure of $X$} is defined as the union
$$
\calA_\phg^{*}:=\bigcup_{\gamma,p}\calA_\phg^{\gamma,p}.
$$
For later use, we define the index set associated to a function $u\in \calA^0_\phg$ to be the set
\beq
\label{EuEq}
E_u:=\{(\sigma_j,p)\,:\, a_{jp}\not\equiv0
\}.
\eeq
In Definition \ref{phgDef}, the $\sim$ symbol means that
$f$ admits an asymptotic expansion in powers of $r$ and $\log r$.
By definition, this means that for every $k\in\NN$ and $m\in\{0,\ldots,N_k\}$,
\beq\label{remainderEq}
f-\sum_{j=0}^{k-1} \sum_{p=0}^{N_j} a_{jp}(\theta, y) r^{\sigma_j} (\log r)^p
-\sum_{p=0}^{N_k-m} a_{kp}(\theta, y) r^{\sigma_k} (\log r)^p
\in
\begin{cases}
\calA^{\sigma_k,N_k+1-m}_\phg & \h{if } m\in\NN,\cr
\calA^{\sigma_{k+1},N_{k+1}}_\phg & \h{if } m=0,\cr
\end{cases}
\eeq
and moreover that corresponding remainder estimates hold whenever any number of
vector fields in $\calV_b(X)$ are applied to the left hand side of \eqref{remainderEq}.
These are referred to as remainder estimates since if a function $u$ lies in
$\calA^{\gamma,p}$, it satisfies $|u|\le Cr^\gamma(\log r)^p$ on $X$.

Several remarks are in order.
First, on a smooth space, Taylor's theorem implies that any smooth function
admits a Taylor series expansion. In our setting,
however, there is a certain `gap' between the space $C^\infty_e$
(or even $\calA^0$) and
the space of functions admitting a bi-graded expansion as in Definition \ref{phgDef}.
Second, the push-forward of a function in $\calA^0_{\phg}$
(unlike a function in $C^\infty_e$)
can be considered as a continuous function on $M$,
provided its leading term is independent of $\theta$
 (and not just on $M\sm D$).
Third, expansions of polyhomogeneous functions are rarely convergent, but only give
`order of vanishing' type estimates, as described above.
Fourth, the remainder estimate \eqref{remainderEq} can be written as
an equality
\beq\label{remainder2Eq}
f-\sum_{j=0}^{k-1} \sum_{p=0}^{N_j} a_{jp}(\theta, y) r^{\sigma_j} (\log r)^p
-\sum_{p=0}^{N_k-m} a_{kp}(\theta, y) r^{\sigma_k} (\log r)^p=u,
\eeq
with $u$ belonging to one of the two spaces in the right hand side of \eqref{remainderEq},
and this equality is understood to hold
on some ball (in the coordinates $r,y,\theta$) centered at a point on $D$.
It is usually difficult to control the
size of this ball, i.e., to give lower bounds for its radius.
However, since such a positive radius
exists for each $p\in D$ and $D$ is compact, the equality \eqref{remainder2Eq}
actually holds on some tubular neighborhood (of positive distance---as
$r$ is uniformly equivalent to the distance function close
enough to $D$) of $D$.
Finally, note that the polyhomogeneous structure associated to either
the $(\rho,\th,y)$ or the $(r,\th,y)$ coordinates is the same, precisely because
we are allowing fractional powers of $r$ and $\be\log\rho=\log r$.

\subsection{Edge and wedge
function spaces}
\label{EdgeWedgeSubSec}

There are two distinct scales of function
spaces naturally associated to a \K edge metric, that
we will denote by $C^{k,\al}_s$
$$
\hbox{\sl with $s$ equal to either $w$ or $e$}.
$$
We consider both of these as subspaces of $L^\infty(M)$.

The first, the wedge spaces $C^{k,\al}_w$, are the usual $C^{k,\al}$ spaces
on $M\sm D$
with respect to a(ny) \K edge metric, intersected with $L^\infty(M)$ (and,
in fact, are contained in $C^0(M)$). That is, we consider
$M\sm D$ as an {\it incomplete} Riemannian manifold with the usual distance
to the edge being a distance function.

The second, the edge spaces $C^{k,\al}_e$ \cite{Mazzeo},  are the
intersection of $L^\infty(M)$ with the usual $C^{k,\al}$ spaces
on $M\sm D$
with respect to a conformal deformation of a \K edge metric for
which the logarithm of the usual distance
to the edge is now a distance function. Thus, we consider
$M\sm D$ as a {\it complete} Riemannian manifold.
Equivalently, $C^{k,\al}_e$ consists of functions on $M\sm D$ that
are push-forwards of functions on $X$ that are $C^{k,\alpha}$ with respect
to the edge differentials $\calV_e(X)$ of \S\ref{phgSubSec}.

Let us now explain how to define these function spaces
explicitly in the model edge $C_\be\times\RR^{2n-2}$.
In the notation of \S\ref{ConSubSec}, the flat metric
on a cone $C_\be$ takes the form
$dr^2 + r^2 d\theta^2$, with $\th\in[0,2\pi\be)$.
Thus the flat model edge metric $\o_\be$ is given by
$dr^2 + r^2 d\theta^2+ dy^2$, with $y=(y_1,\ldots,y_{2n-2})$
coordinates on $\RR^{2n-2}$, and it is this metric that
defines the spaces $C^{k,\al}_w$. The conformally rescaled
metric $(d\log r)^2 +  d\theta^2+ r^{-2}dy^2$ defines the spaces
$C^{k,\al}_e$.
It follows that the
defining vector fields for the spaces $C^{k,\al}_s$
are $r^{\delta_e(s)}\{\frac{\del}{\del r}, \frac1r\frac{\del}{\del \th},
\frac{\del}{\del y_j}\}$, where $\delta_e(e)=1,\delta_e(w)=0$.
We will say a bit more about these function spaces below, but
for a thorough discussion of these function spaces
as well as the different coordinate  choices involved we
refer to \cite[\S2]{JMR}.
For the moment we observe the obvious inclusion
$C^{k,\gamma}_w \subset C^{k,\gamma}_e$; the wedge
spaces are in fact much smaller than their edge
counterparts. E.g., as noted earlier, $\sin\log r\in C^\infty_e$ shows that
$C^\infty_e\not\subset C^0(M)$ (though it is contained in $L^\infty(M)$),
and $r^{k+\eps}\in C^{k,\eps}_w\cap C^\infty_e$ but is not contained
in any higher wedge space.

Note, finally, that the space $\calA^0$ defined
in \S\ref{phgSubSec} admits a similar geometric description.
Namely it is the space of bounded functions that are
$C^\infty$ space with respect to the metric
$d(\log r)^2+d\th^2+dy^2$, namely the product metric
obtained by conformally rescaling the flat metric
on the cone together with the standard (non-rescaled)
Euclidean metric on $\RR^{2n-2}$.

\subsection{Various \Holder domains}
\label{HolderSubSec}
Perhaps surprisingly, it turns out that the complex
\MA equation, and, as a special case, the Poisson equation,
cannot be solved in general in $C^{2,\al}_w$.
In fact, as already noted, 
$\Re\, z_1$ is pluriharmonic and belongs merely
to $C^{1,\frac1\be-1}_w$, which fails to lie in
$C^{2,\al}_w$ when $\be>1/2$. On the other hand,
$C^{2,\al}_e$ is certainly a large enough space
to find solutions; however, it is not even contained in $C^0$
and so seems to give too little control/regularity
to develop a reasonably strong existence theory
for the \MA equation.

First, we introduce the {\it maximal H\"older domains}
$$
\Ds(\Delta_\o)=\Ds:=\{u\in C^{0,\gamma}_s\,:\, \Delta_\o u\in C^{0,\gamma}_s\}.
$$
These are Banach spaces with associated norm
$$
||u||_{\Ds(\Delta_\o)}:=
||\Delta_\o u||_{C^{0,\gamma}_s}+||u||_{C^{0,\gamma}_s}.
$$
We also define the {\it little H\"older domains}
$$
\tildeDs(\Delta_\o),
$$
as the closure of the space $\calA^0_\phg$ of \phgs\ functions
in the $\Ds(\Delta_\o)$ norm.
The name `maximal' comes from the analogy with the usual definition
of the maximal domain of the Laplacian in $L^2$, namely
$\calD_{\max}(\Delta):=\{u\in L^2\,:\,\Delta u\in L^2\}$.
On the other hand, the `little' spaces are defined similarly to the
usual little H\"older spaces (see, e.g., \cite{Lunardi}).
The latter are separable, unlike the former, and of course
$\tildeDs(\Delta_\o)\subset \Ds(\Delta_\o)$.

The space $\Dw$ was introduced by Donaldson \cite{Donaldson-linear} (where it is
denoted $C^{2,\gamma,\be}$); it
gives wedge \Holder control of the wedge Laplacian, which
is a sum of {\it certain} second wedge derivatives of type (1,1).
Motivated by this, the space $\De$ was introduced in \cite{JMR}, and gives
edge \Holder control of the wedge Laplacian. Thus, unlike $\Dw$, it is
a sort of hybrid space. In fact, $\Ds$
can be characterized by requiring the (1,1)
wedge Hessian (and not only its trace!) to be H\"older
(in fact, a slightly stronger characterization holds, see
Theorem \ref{DsThm} below).
When $\be\le1/2$ an even sharper
characterization holds for $s=e$: the full (real) wedge Hessian
is H\"older.
Let $P_{1\bar1}:=\del_r^2 + \frac{1}{r}\del_r + \frac{1}{\beta^2 r^2} \del_\theta^2$, and define
$$
\calQ:=\{
\del_r, r^{-1} \del_\theta, \del_{y_a}, \del^2_{y_a y_b},
\del_r \del_{y_a}, \del_r \del_\theta,
r^{-1}\del_\theta \del_{y_a}, P_{1\bar{1}}, a,b=1,\ldots2n-2\}.
$$

\begin{thm}
\label{DsThm}

%\hfill\break
(i) There exists a constant $C>0$ independent of $u$
such that
$$
||Tu||_{s;0,\gamma}\le C(||\Delta_\o u||_{s;0,\gamma}+ ||u||_{s;0,\gamma}),
$$
for all $T\in \calQ$.
Thus, $u\in \Ds$ if and only if $Tu\in C^{0,\gamma}_s$ for all $T\in \calQ$.
\hfill\break
(ii)
If $\be\in(0,1/2]$, the previous statement holds with $\calQ$ replaced
by
$$\calQ\cup\{
\del^2_r, r^{-1}\del_r,
 r^{-2}\del_\th^2,\hfill\break r^{-1}\del_r\del_\th,
r^{-1}\del_r\del_{y_a},  a,b=1,\ldots2n-2\}.
$$
Also, $\Dw\subset C^{2,\min\{\frac1\be-2,\gamma\}}_w$.
\end{thm}

\begin{remark}{\rm
Part (ii) quantifies the difference between
the ``orbifold" regime $\be\in(0,1/2]$ and the harder
$\be\in(1/2,1)$ regime. In fact, one can write down stronger
and stronger characterizations of $\De$ under mild assumptions when $\be\in(0,\be_0)$
as $\be_0\le1/2$ approaches 0, however these will not involve
any new second order operators beyond those in (ii)
except $r^{1-\frac1{\be_0}}\del_r\del_\th$ and
$r^{-\frac1{\be_0}}\del_{y_a}\del_\th$.
For instance, when $\be\in(0,1/3)$,
the operators
$T=r^{-3}\del_\th,
r^{-2}\del_r\del_\th,r^{-3}\del_{y_a}\del_\th$ are also allowed.
}
\end{remark}
For the proof of (i) and (ii) with $s=e$, we
refer to \cite[Proposition 3.3]{JMR}, and for (i) and (ii) with $s=w$ to
\cite[Proposition 3.8]{JMR}. The key ingredient in the proofs
is a precise description of the singularities (or, in other words,
of the the polyhomogeneous structure) of the Green kernel
of the Laplacian of the (curved) reference metric $\o$ \eqref{oEq}, considered as
a distribution on certain blow-up of $X\times X$. This is
explained in \S\ref{GreenSubSec}.
For the proof of (i) with $s=w$ for most of
the operators in $\calQ$, we also refer to Donaldson \cite[Theorem 1]{Donaldson-linear};
the verification for the remaining operators is straightforward
from his arguments. Donaldson's approach is more elementary in
that he only obtains the polyhomogeneous structure of the Laplacian
of the flat model metric $\o_\be$ \eqref{obetaEq} (which can be done
by separation of variables arguments, and also goes back to the work
of Mooers \cite{Mooers}); using the Schauder estimate for this flat
model he then obtains (i) with $s=w$ by a partition of unity argument.
The advantage of the extra work done in \cite{JMR} is that
the polyhomogeneous structure of the Green kernel of the curved metric
itself allows much more refined information, e.g., (ii), and also leads
eventually to the higher regularity of solutions of the \MA equation---see
\S\ref{HigherRegSubsec}---that does not seem to be accessible using the arguments of \cite{Donaldson-linear}.

Next, similarly to the definition of the spaces $C^{k,\al}_s$,
also the spaces $\Ds$ can be defined with respect to {\it any}
\K edge metric of the form $\ovp$ with $\vp\in\tildeDs$.
This property is absolutely crucial in applications to `openness' along the
continuity method
(\S\ref{ConvergenceSubSec}), and to higher regularity (\S\ref{HigherRegSubsec}).
\begin{thm}
\label{DsThm2}
Suppose that $u\in\tildeDs$.
Then
$$
\Ds(\Delta_\o)=\Ds(\Delta_{\o_u}):=\{v\in C^{0,\gamma}_s\,:\,\Delta_{\o_u}v\in C^{0,\gamma}_s\}.
$$
\end{thm}

This is proved in \cite[Corollary 3.5]{JMR}; 
Donaldson \cite{Donaldson-linear} does not state any such result explicitly, 
but briefly sketches related ideas in the case $s=w$ 
in \cite[p. 64]{Donaldson-linear}.
This result is crucial in our approach and we describe the proof in some
detail below.

The goal is to show that
$\Ds(\Delta_\o)=\Ds(\Delta_{\o_u})$.
First, one observes that this is true when $u$ is \phgs, by
the explicit polyhomogeneous structure of the Green kernel
of \phgs\ elliptic edge operators (see
Theorem \ref{MazzeoThm} below).
 Next, note that
$$
\Delta_{\o_u}f=\frac{\o^n}{\o_u^n}\frac{n\i\ddbar f\w\o_u^{n-1}}{\o^n}.
$$
Thus, if $u\in\Ds$ then $|\D_{\o_u}f|\le C_1(|\D f|+(\i\ddbar f,\alpha)_\o)$,
where $|\alpha|_\o\le C_2$, and $C_1,C_2$ are both controlled by a polynomial
function of
$\sum_{i,j}[u_{i\b j}]_{s;0,\gamma}$. By Theorem \ref{DsThm} (i) these constants are then
also controlled by $||u||_{\Ds}$. Thus,
$\Ds(\Delta_\o)\subset\Ds(\Delta_{\o_u})$.

Since $\Delta_{\o_u}$ is injective
on the $L^2$-orthogonal complement to the constants
in $\Ds(\Delta_{\o_u})$, that we denote by
$\Ds(\Delta_{\o_u})'$, it is also injective
on $\Ds(\Delta_{\o})'$. The main point is that
it is also surjective when restricted to the latter space;
since $\Delta_{\o_u}:\Ds(\Delta_{\o_u})'\ra C^{0,\gamma}_s$ is by definition a bijection,
this immediately gives the reverse inclusion and concludes the proof. We now explain
the proof of the surjectivity statement.

The surjectivity follows in several steps. First, given
any nonzero $f\in C^{0,\gamma}_s$ there is a unique
solution $\phi\in \Ds(\Delta_{\o_u})'$ to $\Delta_{\o_u}\phi=f$.
Approximate $u$ in $\Ds(\Delta_{\o})$ by \phgs\ $u^{(k)}$,
and let $\phi^{(k)}$ be the associated solution of
$\Delta_{\o_{u^{(k)}}}\phi^{(k)}=f$.
Since $u^{(k)}$ are \phgs, there is a unique constant $c_k$
such that $\phi^{(k)}-c_k\in \Ds(\Delta_{\o})'$, and by
Theorem \ref{DsThm} (i) (which, as remarked in the beginning
of the proof, applies to $\o$ replaced by any \phgs\ metric)
this implies an estimate
\beq\label{APrioriEstLapkEq}
\sum_{Q\in\calQ}||Q\phi^{(k)}||_{s;0,\gamma}\le C^{(k)}(||f||_{e;0,\gamma}+||\phi^{(k)}||_{e;0,\gamma}).
\eeq
We claim that the constant $C^{(k)}$ in \eqref{APrioriEstLapkEq}
 is locally
uniform in the $C^{0,\gamma}_s$ norm of $u^{(k)}_{i\b j}$
(i.e., that if $u^{(k)}$ were supported in some small ball
then the constant on the right hand side of \eqref{APrioriEstLapkEq}
would be controlled uniformly in terms of the $C^{0,\gamma}_s$ norm of $u^{(k)}_{i\b j}$).
Of course,
since
$||u^{(k)}_{i\b j}-u_{i\b j}||_{s;0,\gamma}$ is uniform in $k$,
this would
imply also that the $C^{(k)}$ are locally uniformly controlled,
independently of $k$, by
the $C^{0,\gamma}_s$ norm of $u_{i\b j}$.
The claim is true for \phgs\ $u$ by the Green kernel construction
of Theorem \ref{MazzeoThm}. For a more general $u\in\tildeDs$, we freeze coefficients
of $\Delta_{\o_u}$ only in a small neighborhood of
any point in $M$ and approximate this operator by a global \phgs\ operator.
The estimate above follows for such a concatenated operator by
the standard method of proof that makes clear that the constant
in the estimate depends only on the local $C^{0,\gamma}_s$ norm of $u_{i\b j}$.
Using a partition of unity to paste these estimates then concludes the
proof of the claim (since both $D$ and $M$ are compact).

Thus, the $\phi^{(k)}$ are uniformly in $\Ds(\Delta_{\o})'$
(note that the constants $c_k$ are uniformly controlled).
When $s=w$ we can now take a subsequence converging
in $\calD_w^{0,\gamma'}$ (for any $\gamma'\in(0,\gamma)$)
to a function $\phi\in\Dw(\Delta_\o)'$ that solves $\Delta_{\o_u}\phi=f$.
When $s=e$, we cover $M\sm D$ by a countable collection of Whitney cubes so that
on each such cube $T\phi^{(k)}$ converges in $C^{0,\gamma'}_e$ to $T\phi$,
for all $\gamma'\in(0,\gamma)$. Taking a diagonal sequence then
produces a solution $\phi$ to $\Delta_{\o_u}\phi=f$ that belongs
to $\Ds(\Delta_{\o})'$, by Theorem \ref{DsThm} (i). In either
case ($s=e$ or $s=w$), (iii) follows.

Finally, we mention several further basic regularity properties of the \Holder domains.

\begin{thm}
\label{DsThm3}
%\hfill\break
(i) For any $T\in\calQ\sm\{P_{11}\}$,
$T$ maps $\De$ into $C^0(M)$. In particular,
if $u\in\De$, then $Tu$ is continuous up to $D$, and
has a well-defined restriction to $D$, independent of $\th$.
\hfill\break
(ii) Let $\gamma\ge0$. Then, $\De\subset C^{1,\al}_w$ for $\al\in(0,\frac1\be-1]\cap(0,1)$.
\hfill\break
(iii) If $u\in\tildeDe$ then $u$ has the partial
expansion near $D$
$$
u=a_0(y)+
(a_{01}(y)\sin\th+b_{01}(y)\cos\th)r^{\frac1\be}
+O(r^2).
$$

\end{thm}

For (i) refer to \cite[Corollary 3.6]{JMR},
for some of the operators mentioned, while the proof for the
remaining ones follows from (ii).
Part (ii) is the singular analogue for the usual $C^{1,\al}$ estimates
for a function in $W^{2,\infty}$, and can be proved using
\cite[Proposition 3.8]{JMR}---we sketch the proof in
Lemma \ref{HeBoundedLemma}. In fact, (ii) implies
the statement of (i) holds for a larger class of operators.
Finally, (iii) is a corollary of Theorem \ref{DsThm2} and the polyhomogeneous
structure of the Green kernel stated in Theorem \ref{MazzeoThm}.

\subsection{\K edge metrics}
\label{KedgeSubSec}

We start with a completely elementary, yet important, lemma.

\begin{lem}
\label{KedgeLemma}
Let $g$ be a continuous K\"ahler metric
on $M \setminus D$, and such that
in any
local holomorphic coordinate system near $D$ where $D = \{z_1 = 0\}$, and $z_1 =  \rho e^{\i\theta}$,
\begin{equation}
\label{gEdgeDef}
g_{1\bar{1}} = F \rho^{2\beta-2},\ g_{1 \b j } = g_{i \bar{1}} = O(\rho^{\beta - 1 + \eps}),\ \mbox{and all other}\ g_{i\b j } = O(1),
\end{equation}
for some $\eps \ge 0$, where $F$ is a bounded nonvanishing function which is
continuous at $D$. Then there exists some $C>0$
such that in any such coordinate
chart
\beq
\label{QIIneq}
\frac1C g_\be\le g \le C g_\be.
\eeq
Moreover, the converse implication holds.
\end{lem}

\bpf
This amounts to showing that $\frac1CI_\be\le[g_{i\b j}]\le CI_\be$,
where
$I_\be:=\diag(\rho^{2\be-2},1,\ldots,1)$.
Let $v=(v_1,\ldots,v_n)\in\CC^n$ be any vector. Then,
$$
\sum_{i,j=1}^n g_{i\b j}v_i\overline{v_j}
\le
g_{1\b1}|v_1|^2+\sum_{i,j=2}^ng_{i\b j}v_i\overline{v_j}
+
\sum_{j=2}^n|g_{1\b j}|^2|v_1|^2+\sum_{j=2}^n|v_j|^2
\le
C\rho^{2\be-2}|v_1|^2+C\sum_{j=2}^n|v_j|^2,
$$
for some $C>0$, proving the second inequality.
The first inequality follows similarly, since
\eqref{gEdgeDef} implies that
\begin{equation}
\label{gEdgeInvDef}
g^{1\bar{1}} = F \rho^{2-2\beta},\
g^{1 \b j }, g^{i \bar{1}} = O(\rho^{1-\beta + \eps}),\ \mbox{and all other}\ g^{i\b j } = O(1).
\end{equation}

Conversely, choosing $v=(0,v')$ shows that
$CI\le g':=[g_{i\b j}]_{i,j=2}^n\le C'I$, and thus
$g_{i\b j}=O(1)$, for all $i,j\ge2$.
\eqref{QIIneq} implies $C^{-n}\o_\be^n\le\o^n\le C^n\o_\be^n$,
but
$$
\det[g_{i\b j}]
=
g_{1\b 1}\det[g_{i\b j}]_{i,j=2}^n-|(g_{1\b2},\ldots,g_{1\b n})|^2_{g'},
$$
where $|v|^2_{g'}:=v^H[g_{i\b j}]_{i,j=2}^nv$.
It follows that $|g_{1\b j}|\le C\rho^{\be-1},$ for all $j\ge2$.
\epf

We define
\beq
\label{KEdgePotentialsEtaSpace}
\calH^\eps_{\o} := \{\vp\in C^\infty(\MsmD)\cap C^0(M)
\, :\,   \ovp>0 \h{\ on\ } M,
\mbox{\ and $\o_\vp$ satisfies \eqref{gEdgeDef}}\}.
\eeq
We call
$$
\calH^e_\o:=\cup_{\eps>0}\calH^\eps_\o
$$
the {\it space of  \K edge metrics}.

We note that another, simpler, way of deriving the preceeding lemma,
but which conceals some of what is going on,
is by working in the singular coordinate chart
$(\zeta,z_2,\ldots,z_n)$.
In a nutshell,
in terms of the singular coordinates,
the definition of a \K edge metric
simply means that the
cross terms in the matrix $[g_{i\b j}]$ (i.e., terms for which
precisely one of $i$ and $j$ equals 1) have a fixed
rate of decay $O(r^{\eps/\beta})$ near $D$, so that the metric
is asymptotically a product. In particular, this means the corresponding
Laplacian is also approximated in a certain precise sense by the
Laplacian of the product of a flat cone and $\CC^{n-1}$ (i.e.,
the Laplacian of $\o_\be$).
This is very important in proving existence of asymptotic expansions,
and structure results for the Green kernel, as we discuss later.
In fact, as will follow from the general theory
we explain later, in essentially all of
the discussion below one may take concretely
$\eps=\min\{1-\be,\beta\}$.
To be more precise, we will show that solutions
to essentially any reasonable complex \MA equation always
lie in $\calH^e_\o$, and KEE metrics in fact even
lie in $\calH^\be_\o$.

We end this subsection with a lemma that shows that
$\calH^e_\o$ is a natural choice of space
of \K metrics in this context. Indeed, once we have found
one reference \K edge metric, we may produce many other
such metrics by adding a \K potential with {\it merely bounded}
complex Hessian. This fact may seem counterintuitive at first. It
is absolutely crucial for all that follows.
Its proof relies crucially on the fact that $\be<1$,
or ultimately on the fact that no nonzero $2\pi\be$-periodic
function $b$ can satisfy $b_{\th\th}+b=0$, i.e., $-1$ is
not in the spectrum of the Laplacian on $S^1(2\pi\be)$,
which equals $-\NN_0/\beta^2$.

\begin{lem}
\label{HeBoundedLemma}
Suppose that $\eta\in \calH^\eps_\o$ and that
$[u_{i\b j}]$ is bounded. Then
$\eta+\i\ddbar u\in \calH^{\tilde\eps}_\o$,
where $\tilde\eps=\min\{\eps,\be,1-\be\}$.

\end{lem}

This is a corollary of Theorem \ref{DsThm3} (ii).
The proof of this result uses the basic Schauder type
estimates in the edge spaces.
Indeed, the assumption
implies of course that $\Delta_\o u=f\in L^\infty$.
Thus $u=Gf+k$ with $k$ an element of the $L^2$ kernel
of $\D_\o$. One may show that any element of the kernel
can be written in the form $a_0(y)+a_1(y,\th)r^{1/\be}
+O(r^{2+\delta})$, and thus
$\del_y\del_r k=O(r^{\min\{\frac1\be-1,1\}})$.
On the other hand, it is not difficult to show that
both $\del_r\circ G$ and $\frac1r\del_\th\circ G$ map
$L^\infty$ to itself (see Theorem \ref{DsThm} (i)).
Thus, $\frac1r\del_\th u\in L^\infty$
and $\del_r u\in L^\infty$. We will use both of these facts
in a moment. Now, by our assumption,
$\del_y\del_r u=O(1)$.
Integrating twice we find that $u=v(y,\th)+a(r,\th)+O(r)$.
Applying the two boundedness results from above it follows
that $\frac1r\del_\th(v+a)=0$. But this implies that
both $\del_\th v(y,\th)$ and $\del_\th a(r,\th)$ are functions of $\th$ alone,
but plugging that back in we see that necessarily then
$v+a$ must be a equal to $v_1(y)+a_1(r)$, with $a_1(r)=O(r)$,
so $u=O(r)+v_1(y)$. We fix a value of $y$ and regard $u$
as a function on the flat cone $\{y\}\times C_\be$. We know
that $\D_\be u=P_{11}u=F(y)\in L^\infty$, and so $u=G_\be F+k(y)$.
But now, by the properties of $G_\be$ we conclude that
$u$ must vanish to order $O(r^{\min\{1/\be,2\}})$, and so
it follows that $\del_r\del_y u=O(r^{\min\{1/\be-1,1\}})$, from
which the lemma follows.

\begin{remark}
\label{tildeRemark}
{\rm
Essentially, if one assumes that $F$ in \eqref{gEdgeDef}
is also in the little \Holder space associated to $C^{0,\gamma}_e$,
and one defines $\widetilde\calH_\o^e$ accordingly,
then it would follow that $\widetilde\calH_\o^e\subset\tildeDe$. It would be interesting
to also prove upper bounds on $\tildeDe$.
}
\end{remark}

\subsection{The reference geometry}

{\sl ``[A] Hermitian metric has the peculiarity of favoring negative curvature
over positive curvature."}

\bigskip
\hfill
%glue3.5in 
---  Solomon Bochner \cite[p. 179]{Bochner1947}.

\bigskip

One of the novelties of \cite{JMR} was to  prove
a priori second order estimates for a fully nonlinear PDE
without curvature bounds on the reference geometry, but only with
a one-sided curvature bound.
We are not aware of other situations where this is possible.
Indeed, the next lemma shows that there is no uniform bound,
in general, for the curvature of the reference geometry.
But, it provides a one-sided bound, which turns out to be
very useful for the Laplacian estimates (see \S\ref{CLSubsubsec}).

\begin{lem}
\label{LiRLemma}
The bisectional curvature of the reference metric $\o$
(recall \ref{oEq}) is bounded from above. In general, it is not
bounded from below.
\end{lem}

The proof, due to Li and the author, appears in \cite[Appendix]{JMR} and its
generalization to the case of normal crossings in \cite{MR2012,MR}.
It relies on a careful computation, using an adapted normal coordinate
system appearing in the work of Tian--Yau \cite{TY1990}.
The lack of the lower bound can be seen directly
from the computations in \cite[Appendix]{JMR} by considering
the $R_{1\b 1 1\b 1}$ component of the bisectional curvature,
and observing that the upper bound can be made as negative
as one wishes.
If one is only interested in proving the upper bound
(but without proving that no lower bound exists),
the proof of Lemma \ref{LiRLemma} can be simplified considerably, as pointed out to
the author by J. Sturm in November 2013. Indeed, working locally on a ball $B\subset M$
centered at a point in $D$, one considers the holomorphic map
$F:(z_1,\ldots,z_n)\mapsto (z_1,\ldots,z_n,z_1^{\frac1\be})$. Let
$\pi:\CC^{n+1}\ra\CC^n$ denote the projection to the first
$n$ components. The \K form $\pi^\star\o_0+\i\ddbar(H(z_1,\ldots,z_n)|z_{n+1}|^{2})$
is a smooth metric on a ball in $\CC^{n+1}$ whose
bisectional curvature is uniformly bounded.
Pulling this form back from the complement of $\{z_1=0\}\cup\{z_{n+1}=0\}$
under
the holomorphic map $F$ yields the \K form
$\o_0+\i\ddbar(H|z_1|^{2\be})$  on $\CC^n$
whose bisectional curvature can only decrease
by the classical expression for the second fundamental form
of a complex submanifold \cite[p. 79]{GH}. Even though
the map $F$ is only multivalued, the pull-back of the \K
form above is continuous single-valued since $H$ is
a smooth function of $z_1$ and the expression $|z_1|^{2\be}$
is single-valued under this pullback.

A corollary of both statements in the Lemma is that
the Ricci curvature of $\o$ is unbounded from below: if it were not, the upper
bound on the bisectional curvature would force a lower one as well.
As alluded to above, this fact motivates using the Ricci continuity
method that starts out, morally, with a metric whose Ricci curvature bound
is $-\infty$, and gradually produces a better lower bound until,
eventually, arriving
at the \KE metric.

In the case of the model edge $C_\be\times\CC^{n-1}$ the curvature
is identically zero outside the divisor. In general, it is natural to expect
that the holomorphic submanifold geometry of the divisor $D$
in $M$ should be related to whether a reference metric exists with
bounded bisectional curvature. A first step in this direction
was taken by Arezzo--Della Vedova--La Nave \cite{ArezzoVedovaLaNave2014}.

\subsection{The structure of the Green kernel}
\label{GreenSubSec}

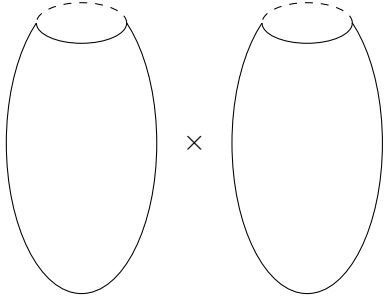
\begin{figure}
\centering
	\begin{tikzpicture}
		\begin{scope}
			\clip (-1.2,-2.2) rectangle (1.2,1.6);
			\draw (0,0) ellipse (1 and 2);
			\draw (0,1.6) ellipse (.6 and .27);
		\end{scope}
		\begin{scope}
			\clip (-1.2,1.6) rectangle (1.2,2.2);
			\draw [dashed] (0,1.6) ellipse (.6 and .27);
		\end{scope}
%%%%%%%%%%%%%%%%%%%%%%%
		\begin{scope}[shift={(3,0)}]
			\begin{scope}
				\clip (-1.2,-2.2) rectangle (1.2,1.6);
				\draw (0,0) ellipse (1 and 2);
				\draw (0,1.6) ellipse (.6 and .27);
			\end{scope}
			\begin{scope}
				\clip (-1.2,1.6) rectangle (1.2,2.2);
				\draw [dashed] (0,1.6) ellipse (.6 and .27);
			\end{scope}
		\end{scope}
		\node at (1.5,0) {$\times$};
	\end{tikzpicture}
\caption{The space $X\times X$.}
\label{FigXtimesX}
\end{figure}

In this subsection we describe in detail the structure of the Green
kernel in the case $M$ is $S^2$ with $D$ a single point. The general
case is not much more complicated, but we prefer concreteness over
generality in this discussion.

Thus we start with the singular manifold $(S^2,p)$ and blow-up the
point $p$ as in Figure \ref{RealBlowUpFig}. What this means is that we `separate' the different
directions in which one may approach the point $p$. Each of these
directions is now a separate point in the blow-up.
As in \S\ref{phgSubSec}, what this real
blow-up amounts to is to introduce polar coordinates $(r,\th)$ around $p$,
with $r(p)=0$. Then $X=\h{Bl}_p^{\RR_+}S^2$ is the disjoint union
of $S^2\sm\{p\}$ and a circle $S^1_{2\pi\be}$ of radius $2\pi\be$.
It comes equipped with a blow-down map $X\ra S^2$ and the inverse
image of $p$ is $S^1_{2\pi\be}$.
(In the real setting a real blow-up sometimes also refers to an
$\RR$ blow-up where the resulting fiber over $p$ is a half circle
or $\RR\PP^1$; therefore we used the superscript $\RR_+$ for
our `oriented' blow-up.)

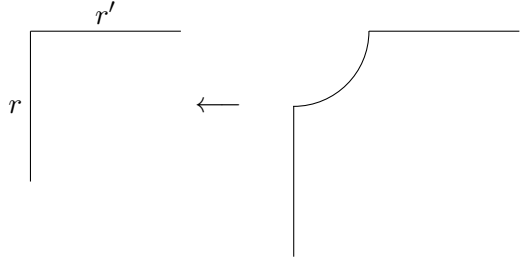
\begin{figure}
\centering
	\begin{tikzpicture}
		\begin{scope}[shift = {(0,0)}]
			\draw (0,-2) -- (0,0);
			\draw (2,0) -- (0,0);
		\end{scope}
		\begin{scope}[shift = {(3.5,0)}]
			\draw (0,-3) -- (0,-1);
			\draw (3,0) -- (1,0);
			\clip (0,-1.1) rectangle (1.1,0);
			\draw (0,0) circle (1);
		\end{scope}
	\node at (2.5,-1) {$\longleftarrow$};
	\node at (-.2,-1) {$r$};
	\node at (1,.25) {$r'$};
	\end{tikzpicture}
\caption{The real blow-up of the origin in the quadrant; the origin represents
any point in the diagonal of $\del X\times\del X$.}
\label{Figrrprimeblowup}
\end{figure}

We do not describe the Green kernel $G_\o$ associated to $\o$ (or
any phg edge metric) on $X\times X$ (see Figure \ref{FigXtimesX}),
but instead pull it back under one more blow up map.
The purpose of this is to `separate' certain directions near the boundary diagonal.
In other words, if $\{(0,\th)\}\times\{(0,\th')\}$ is
a point on $\del X\times \del X$ there are many different ways
to approach it from the interior: if $r,r'$ are the radial variables
on each of the two copies of $X$, we may let $r$ approach zero
faster than $r'$ or vice versa. Thus we consider $r,r'$
as coordinates on the positive orthant $\RR^2_+:=\RR_+\times\RR_+$
and blow-up the origin (see Figure \ref{Figrrprimeblowup}). The result
$\h{Bl}_{(0,0)}^{\RR_+}\RR_+^2$ is the disjoint union
of $\RR_+^2\sm\{(0,0)\}$
and a quarter circle $S^1_{++}$ of radius $\pi/2$.
The quarter circle parametrizes various values of the ratio $r/r'$,
via the map $\tan^{-1}:[0,\infty)\ra S^1_{++}$.
This blow-up is a purely local construction; we may thus perform
it on $X\times X$, which amounts to blowing-up each point
on the boundary diagonal $\del X\times\del X$, or in other words
blowing-up that whole submanifold. We denote the resulting
space by the edge double space
$$
X^2_e:=\h{Bl}_{\del X\times\del X}^{\RR}X\times X
$$
(here the $\RR$ and $\RR_+$ blow-ups coincide). In higher dimensions,
there are an additional $2n-2$ conormal directions $y$. Then we only
blow-up the fiber diagonal $\{r=r'=0, y=y'\}$ which is
a strict submanifold of $\del X\times\del X$.
We denote by
$$
\pi:X^2_e\ra X^2
$$
the blow-down map.
Observe that $\pi^{-1}(\del X\times\del X)$ is the union
of three hypersurfaces:
$$
\h{rf}:=\pi^{-1}\big(\{r=0 \}\big), \quad
\h{lf}:=\pi^{-1}\big(\{r'=0 \}\big),
$$
called the right and left faces,
and the new hypersurface
$$
\h{ff}:=\pi^{-1}\big(\{r=r'=0 \}\big),
$$
diffeomorphic to $S^1\times S^1\times S^1_{++}$,
called the front face.
These hypersurfaces all have coordinate descriptions
in terms of polar coordinates about the corner $\{r=r'=0\}$.
Namely, denote $R:=\sqrt{r^2+r'^2}$ and set
$(\psi,\psi'):=(r/R,r'/R)\in S^1_{++}$. The double edge space
then is parametrized near its boundary by $(R,\psi,\psi',\th,\th')$
(while, by comparison, $X^2$ was parametrized by
$(r,r',\th,\th')$), and
$\h{rf}=\{\psi=0\},
\h{lf}=\{\psi'=0\},
\h{ff}=\{R=0 \}.
$
Away from its boundary $X^2_e$ is locally diffeomorphic
to $M\times M$, of course, and is parametrized there
by coordinates on the latter.

Now, we return to the general setting of a K\"ahler edge
manifold, and describe $G_\o$ as the push-down
of a distribution on $X^2_e$ under the blow-down map $\pi$.
For concreteness, the reader may focus on the example given above,
even though the result below is stated in the general setting.
The first main reason for doing so
is that $G_\o$ is of course singular along the diagonal of $X^2$,
however the diagonal intersects the boundary nontransversally,
while the lifted diagonal
$$
\diag_e:=\{(R,1/\sqrt2,1/\sqrt2,\th,\th)\,:\, R\in\RR_+,\th\in S^1\}
$$
intersects ff transversally at $\{(0,1/\sqrt2,1/\sqrt2,\th,\th)\}\cong S^1$,
and the intersection points lie in the interior of ff.
Another reason for working on $X^2_e$
is the dilation invariance structure.
Finally, here is the description of $G_\o$, a corollary
of a general result of Mazzeo \cite[Theorem 6.1]{Mazzeo},
see \cite[Proposition 3.8]{JMR}.

\begin{thm}
\label{MazzeoThm}
Let $g$ be a polyhomogeneous edge metric with index
set $E_g$ with angle
$\beta$ along $D$ and denote by $G$
the generalized inverse to
the Friedrichs extension of $-\Delta_g$.
Then the Schwartz kernel $K_G$ of $G$ is a
distribution on $X^2_e$ that can be decomposed
as $K_G=K_1+K_2$ with the following properties:
\hfill\break
$K_1$ is supported on a neighborhood of $\diag_e$
disjoint from lf and rf;
$R^{2n-2}K_1$ has a classical
pseudodifferential singularity of order $-2$
on $\diag_e$
that remains conormal when extended to ff.
\hfill\break
$K_2$ is polyhomogeneous on $X^2_e$ with index
sets $2-2n$ at ff, and $E$ at both rf and lf,
where
$$
E\subset \{ (j/\beta + k, \ell): j, k, \ell \in \ZZ_+\ \mbox{and}\ \ell = 0\ \mbox{for}\ j+k \leq 1,\ (j,k,\ell) \neq (0,1,0) \},
$$
when $E_g=\{0\}$, and otherwise $E$ is contained in a larger set that is
determined by $E_g$ and the above set. In particular,
$$
T\circ G:C^{0,\gamma}_e\ra C^{0,\gamma}_e,
$$
is a bounded operator where $T\in\calQ$ is as in Theorem \ref{DsThm}.
\end{thm}

The proof of this occupies a good part of \cite{Mazzeo} and we
will leave a detailed exposition of it in our complex setting to a separate
exposition.

\subsection{Higher regularity for solutions of the \MA equation}
\label{HigherRegSubsec}

The standard theory of elliptic regularity applies directly to the
\MA equation, despite it being fully nonlinear
\cite[\S17]{GT}. Indeed, the linearization of the \MA
equation is simply the Poisson equation (with
a potential)---we review this shortly. It implies that
any $C^{2,\al}$ solution of \eqref{KEEq} is automatically smooth
(for an alternative approach see \cite{DH}).
This assumes, of course, that the reference form $\o$ is smooth.
A penetrating feature of \K edge geometry is that although $\o$ (see
\eqref{oEq}) is no longer smooth, it is possible to develop
an analytical theory that isolates, so to speak, the only
direction in which the geometry is singular, or, equivalently,
in which the associated Laplacian is degenerate (fails to be elliptic).
Thus, KEE metrics turn out to be smooth in all
conormal directions. One could stop here and argue
that this is already the right analogue of higher regularity in the singular
setting. However, as already indicated in \S\S\ref{phgSubSec},
Taylor's theorem does not apply to a merely conormal function.
The analogue of smoothness is thus the space of polyhomogeneous
conormal distributions. This is more than an academic
difference. One of the thrusts of this 
approach is that polyhomogeneity yields 
basic geometric information.
Also, the proofs
of these facts are rather standard by now---essentially a matter
of high-level bookkeeping, somewhat analogously to
tools in algebraic geometry. These proofs have their origin in the
fundamental work of Melrose, developed in-depth in the case
of real edges by Mazzeo, and further developed in the case of
complex codimension one edges in \cite{JMR} and in any codimension
for crossing complex edges of codimension one \cite{MR2012,MR}.
 For the rest of this subsection
we describe, in broad strokes,
how higher regularity is proved in the edge setting
and mention some of its basic analytic and geometric consequences.

A basic fact is that solutions of a general class of complex
\MA equations are automatically polyhomogeneous as soon as they
lie in the maximal \Holder domain $\De$. For concreteness,
we concentrate on the special class of equations that we
will consider in \S\ref{RCMSec}.
We denote by $\PSH(M,\o)$ the space of $\o$-plurisubharmonic
functions on $M$, namely upper semi-continuous functions $u$
from $M$ to $[-\infty,\infty)$ such that $\o_u=\o+\i\ddbar u$
is a nonnegative current on $M$.

\begin{thm}
\label{phgThm}
Let $\o$ be a polyhomogeneous K\"ahler edge metric.
Suppose that $\vp \in \tildeDs\cap \PSH(M,\o)$, satisfies
\begin{equation}
\label{genma}
{\o_\vp^n} = {\o^n}e^{f+c\vp}, \quad \h{\  on $\MsmD$},
\end{equation}
where $c\in\RR$ and $f\in\calA_{\phg}^0$.
Then
$\vp \in \calA_{\phg}^0$.

\end{thm}

The proof, that we now outline, relies mostly on the linear
theory for $\De$. The detailed proof appears in \cite[\S4]{JMR}.

\bpf

Differentiating (the logarithm of \eqref{genma})
in a conormal direction
$y_a$ (the same discussion applies to differentiation in $\th$), the \MA equation becomes a seemingly linear equation,
\beq\label{LinearizedEq}
\Delta_{\ovp}\del_{y_a}\vp=c\del_{y_a}\vp+\del_{y_a}f\in C^{0,\gamma}_e,
\eeq
where we used Theorem \ref{DsThm} (i) for the inclusion.
The nonlinearity comes, of course, from the fact
that the operator $\Delta_{\ovp}$ itself depends
on $\vp$. However, by Theorem \ref{DsThm2}
this is irrelevant, as we assume that $\vp\in\tildeDs\subseteq\tildeDe$! Thus, we conclude that $\del_{y_a}\vp,\del_\th\vp\in\De$.
Again because $\vp\in\De$ this implies that $\del_y(\Delta_{\ovp}\vp)\in
\De$. By induction, $\del_\th^{l}\del_y^k(\Delta_{\ovp}\vp),
\Delta_{\ovp}(\del_\th^{l}\del_y^k\vp) \in C^{0,\gamma}_e$
for all $l,k\in\NN\cup\{0\}$ with $l+k>0$, where $\del_y^k$ denotes any operator
of the form $\del_{y_{a(1)}}\circ\cdots\circ\del_{y_{a(k)}}$, with
$a(i)\in\{1,\ldots,2n-2\}$. Thus, by
composing with $G_{\o_\vp}$, we see that
$\del_y\vp, \del_\th\vp$ are infinitely differentiable with
respect to $\del_\th$ and $\del_y$.
Moreover, by \eqref{LinearizedEq} and standard elliptic
regularity in the edge spaces \cite{Mazzeo} it follows that
$\del_y\vp, \del_\th\vp$ are infinitely differentiable with
respect to $r\del_r$. And the inductive argument
above then gives the same also for $\del^l_\th\del_y^k\vp$.
Namely, $\del_y\vp, \del_\th\vp\in\calA^0$.

Next, we write
$\del_{y_a}\vp=G_{\o_\vp}(\del_{y_a}f+c\del_{y_a}\vp)+\kappa$, where $\kappa\in\ker
\Delta_{\o_\vp}$
is simply a constant.
We show now that $\del_{y_a}\vp\in\calA^0_\phg$.
We claim that the
Green operator
$G_{\o_\vp}$ maps a conormal function to a function in
$r^2\calA^0$, at least modulo something polyhomogoneous (phg).
This is an improvement on the general fact
(cf. \cite[Proposition 3.28]{Mazzeo}) that
such an operator (we do not go into the details of what
``such" means here, but refer to \cite[Lemma 4.2]{JMR} for a
precise statement) maps $\calA^0$ to itself.
By induction, this is all we need to conclude the proof, as we have already
showed that $\Delta_{\ovp}\del_{y_a}\vp=\del_{y_a}f+\del_{y_a}\vp\in\calA^0$; applying $G_{\o_\vp}$ to both
sides would then show $\del_{y_a}\vp\in r^2\calA^0+\calA^0_\phg$.
But, since $f$ is phg it then follows that
$\del_{y_a}f+\del_{y_a}\vp\in r^2\calA^0+\calA^0_\phg$.
And so, applying the previous
claim again we conclude $\del_{y_a}\vp\in r^4\calA^0+\calA^0_\phg$,
and thus by induction $\del_{y_a}\vp\in \calA^0_\phg$.

Here is the idea behind the proof of the claim: by
a theorem of Mazzeo \cite[Theorem 6.1]{Mazzeo} we know that $G_\be$,
the Green kernel of the reference metric $\o$,
has a polyhomogeneous kernel.
If all terms in its
expansion, with one variable frozen,
are $O(r^2)$ or better we would be done. In general though
there will be finitely many (positive, of course)
exponents $\gamma_1,\ldots,\gamma_N$ in the polyhomogeneous
expansion of $G$ that are smaller than $2$.
(In our setting, there is only one such, $1/\be$.)
But then post-composing $G_\be$ with $\Pi_{i=1}^N(r\del_r-\gamma_i)$
eliminates these terms. Integrating the equation
$\Pi_{i=1}^N(r\del_r-\gamma_i)\circ G_\be v=O(r^2)$ in $r$ thus
yields that $G_\be v=O(r^2)+\sum u_i(y,\th)r^{\gamma_i}$.
When $v\in\calA^0$ this yields therefore $G_\be v\in r^2\calA^0+\calA^0_\phg$.
However, we are dealing with the Green kernel $G_{\ovp}$ which
does not necessarily have a phg expansion (since $\vp$ is not
phg as of yet). However, by Theorem \ref{DsThm2}
the domain $\tildeDe(\Delta_{\ovp})$ is independent of $\vp\in\widetilde\calH^e_\o$,
and simply equals $\tildeDe$ (see Remark \ref{tildeRemark}).
In particular, $G_{\ovp}$ has the same
asymptotic behavior near the lf and rf of $X$. In other
words, while it is not phg, it has a partial expansion
(the regularity implied by belonging to $\De$), of
the form $a_0(y)+(a_1(y)\cos\th+a_2(y)\sin\th)r^{1/\be})+O(r^2)$.
Thus, the same argument as above implies that for $v\in\calA^0$ one has
$G_{\ovp} v\in r^2\calA^0+\calA^0_\phg$.

Thus, $\nabla_y\vp\in\calA^0_\phg$, and by the same reasoning
also $\del_\th\vp\in\calA^0_\phg$.
Thus, integrating, $\vp=\vp_0+\vp_1$ with $\vp_0$ a function of $r$ alone,
and $\vp_1(r,\th,y)\in\calA^0_\phg$.
But we already know that
$\vp\in C^\infty_e$.
 Thus $\vp_0$ lies in $C^\infty_e$, but then necessarily also
in $\calA^0$ since it is independent of $y,\th$. Thus,
$\vp_0(r)$ and hence $\vp$ lie in $\calA^0$ and, moreover,
$\vp-\vp_0(r)\in\calA^0_\phg$. Thus, to obtain
Theorem \ref{phgThm} it remains to prove
that $\del_r\vp\in \calA^0_\phg$.

To that end, we apply $\del_r$
to the logarithm of \eqref{genma}. This yields
$\del_r\vp = G_{\ovp}(\del_r f+c\del_r\vp)+\kappa_2$, where now
the term in paranthesis belongs to $r^{-1}\calA^0$, and $\kappa_2\in\RR$.
Since the volume form with respect to which $G_{\ovp}$
is defined is asymptotically equivalent to
$rdrd\th dy$ we may still apply $G_{\ovp}$ to a distribution
in $r^{-1}\calA^0$ and obtain a distribution that
lies in $r^{-1}\calA^0$. But now we know that
$G_{\ovp}\del_r f$ is actually bounded, since
by Theorem \ref{DsThm} (i) $\del_r\vp$ is as $\vp\in\De$. Thus,
again, we may treat the equation as an ODE in $r$,
and show that in fact $\del_r\vp$ lies
in $r\calA^0+\calA^0_\phg$. As before, an induction
results in $\del_r\vp\in\calA^0_\phg$.
In conclusion then $\vp\in\calA^0_\phg$,
completing the proof.
\epf

\begin{remark}{\rm
A precursor to this sort of argument is the work of Lee--Melrose
\cite{LeeMelrose} on the complex \MA equation on
a domain, and related work of Mazzeo on the singular
Yamabe problem \cite{MazzeoYamabe}. We also refer
to Rochon--Zhang \cite{RochonZhang} for recent work. All
of these articles deal with several quite different situations,
but are all on complete spaces, as opposed to our incomplete setting.
}
\end{remark}

We end this subsection by discussing the significance of the terms
that appear in the asymptotic expansion.
When the function on the right hand side $f$ has a certain index set $E_f$
(see \eqref{EuEq})
the solution then has an index set that
will depend on $E, E_\o$, and $E_f$ where $E$ is the index
set of $G_\be$, and $E_\o$ that of $\o$ (see \eqref{oEq}),
by which we mean the index set of $\psi_0+\phi_0$ where
$\psi_0$ is any \K potential for $\o_0$ and $\phi_0$ is
defined in \eqref{phi0Eq}.
In the setting of the \KE equation,
one may determine the index set of the solution from this
observation. This is a somewhat tedious, albeit completely
inductive routine---one simply treats the \MA
equation as an ODE in $r$ and computes the terms
that can occur in the expansion, relying on
our alternative characterization of $\De$ (Theorems \ref{DsThm}
and \ref{DsThm2}).
It follows \cite[Proposition 4.3]{JMR} that
solutions to
essentially all the complex \MA equations considered in the setting
of the \KE equation  have the following
expansion:
when $0 < \be < 1/2$,
\begin{equation}
\label{ExpansionBetaLessHalfEq}
\vp(r,\th,y) =  a_{00}(y) + a_{20}(y)r^2 + (a_{01}(y)\sin\th+b_{01}(y)\cos\th)r^{\frac1\be}
+ a_{40}(y)r^4 + O(r^{4+\eps})
\end{equation}
for some $\eps=\eps(\be) > 0$; when $\be = 1/2$, the asymptotic sum on the
right includes an extra term $(a_{02}(y)\sin2\th+b_{02}(y)\cos2\th)r^4$; finally,
if $1/2 < \be < 1$, then
\begin{equation}
\label
{ExpansionAllBetaEq}
\vp(r,\th,y) = a_{00}(y) + (a_{01}(y)\sin\th+b_{01}(y)\cos\th)r^{\frac1\be} + a_{20}(y)r^2  + O(r^{2+\eps})
\end{equation}
for some $\eps=\eps(\be) > 0$.

\begin{remark}{\rm
When $\be\in(1/3,1/2]$, one may readily show that the term in the expansion after $O(r^{1/\be})$
is $O(r^{1/\be+1})$. This is roughly equivalent 
to the divisor $D$ being totally geodesic for all $\be\in(0,1/2)$
(for $\be\in(0,1/3)$ this is related to Atiyah--LeBrun \cite{AtiyahLeBrun}).
This is another justification for calling the regime $\be\in(0,1/2]$
the orbifold regime.
}
\end{remark}

Note that the term of order $r^{1/\be}$
in \eqref{ExpansionBetaLessHalfEq}--\eqref{ExpansionAllBetaEq}
is annihiliated
by $\del^2/\del\zeta\overline{\del\zeta}$.
From this, and the explicit formulas for the reference metric $\o$
and its derivatives \cite[\S2.2]{JMR} one immediately deduces the following geometric information \cite[Theorem 2]{JMR}.

\begin{cor}
\label{ExpansionCurvCor}
Let $\vp\in\widetilde\calH^e_\o\cap \De$, and suppose that $\ovp$
is KEE or else is a solution of a complex \MA equation
of the type (\ref{genma}) with $E_F$ (see (\ref{EuEq})) suitably small.
Then $\ovp\in C^{0,\min\{1,\frac1\be-1\}}_w$, but in general
its Christoffel symbols and curvature tensor are unbounded.
More precisely, in the coordinates
$(r,\th,y)$, $g_{i\b j}=O(1+r+r^{\frac1\be-1}+r^{\frac2\be-2})$, $R_{i\b j k\b l}=
O(1+r^{\frac1\be-2}+r^{\frac2\be-4})$, while
$\Gamma^k_{ij}=g^{k\b l}g_{i\b l,j}=O(1+r^{\frac1\be-2}+r^{\frac2\be-3})$.
Moreover, $\ovp|_{D}$ is a smooth \K metric.
\end{cor}

As another corollary, Song--Wang observed that the expansion \eqref{ExpansionAllBetaEq}
directly implies that the curvature tensor of a KEE metric
is in $L^2$ with respect to the metric \cite[\S4.1]{SongWang}.

%%%%%%%%%%%%%%%%%%%%%%%%%%%%%%%%%%%%%%%%%%%%%%%

\section{Existence and non-existence}
\label{Sec4}

In this section we survey necessary and sufficient conditions for existence
of KE(E) metrics. We start in \S\ref{NonposSubSec} with the easier nonpositive curvature regime,
where the cohomological criterion \eqref{cohomEq} is necessary and sufficient.
We then move on to obstructions (\S\ref{ReductivitySubSec}--\ref{MabuchiFutakiSubSec}),
and finally, an essentially optimal sufficient condition
for existence in the positive case (\S\ref{PositiveSubSec}). In \S\ref{RSSubSec} we pause to
discuss the Riemann surface case. The existence theorems described in
\S\ref{NonposSubSec} and \S\ref{PositiveSubSec} are the main existence results
for KE(E) metrics, and their proof uses results from Sections \ref{KahlerSection}, \ref{RCMSec}, and
\ref{APrioriSubSec}.

We now review some background that is useful in describing some of
the obstructions below. Denote by
$$
\calHo=\{\vp\in C^\infty(M)\,:\, \ovp:=\o+\i\ddbar\vp>0\}
$$
the (moduli) space of \K potentials representing \K forms
(equivalently, metrics) in a fixed cohomology class $\O=[\o]$.
When we discuss \K edge metrics, the natural replacement is
$\calH_\o^\eps$ defined in \eqref{KEdgePotentialsEtaSpace}.
The corresponding space of \K forms is denoted
by
$$
\calH_\O=\{\al \h{\ is a \K form with } [\a]=\O\}.
$$
These are infinite-dimensional Fr\'echet manifolds,
whose transformations, functions, forms, and vector fields
control different aspects of the KE problem. We first introduce
some basic objects on this space.

The tangent bundle of $\calHo$ is isomorphic to $\Ho\times C^\infty(M)$,
and similarly $T^\star\Ho\isom \Ho\times\G(M,\Lambda^{2n}T^\star M)$, the latter factor
denoting the space of top degree forms on $M$, with
the fiberwise pairing given by integration over $M$.
The Mabuchi metric on $\Ho$ is defined by \cite{Mabuchi1987,Semmes92,Donaldson99}
\begin{equation}
\label{MabuchiMetricEq}
\gM(\nu,\eta)|_\vp:=\frac1V\int_M\nu\eta\,\o_\vp^n,
\quad \nu,\eta\in T_\vp\calH_\o\isom C^\infty(M),
\end{equation}
where $V=[\o]^n/n!$.
Note that the constant vector field $\bf 1$$(\vp)=1$ is of unit norm.
This induces a Riemannian splitting $\Ho=\iota_\o(\HO)\times\RR$, with
$
d\iota_\o(T\HO)\perp {\bf 1},
$
$\HO$ thus identified with a totally geodesic submanifold of $\Ho$
passing through $0\in\Ho$.
Thus, with some abuse of notation, we may speak of geodesics in $\Ho$ or in $\HO$ interchangeably, with the latter meaning geodesics in $\iota_\o(\HO)$.
Remarkably, $\bf 1$ is a gradient vector field for $\gM$
\cite[Theorem 2.3]{Mabuchi1987},
and the corresponding potential
\beq
\label{LEq}
L:\Ho\ra\RR
\eeq
is thus a distance function
for $\gM$ (in the sense that the norm of its gradient is one; as an aside
we remark that it would be interesting to find an interpretation
of this fact in terms of the Mabuchi distance), as first observed by Mabuchi
\cite[Theorem 2.3]{Mabuchi1986},\cite[Remark 3.3]{Mabuchi1987}. It is known as the Aubin--Mabuchi
or \MA energy
since $dL|_\vp=\ovp^n$
(see also \S\ref{NonlineaerEnergiesSubSec}),
and sometimes also referred to as the  Aubin--Yau functional,
and has been studied by many authors, see, e.g.,
\cite{BermanBoucksom,BBGZ,Berman2010,PhongSturmSurvey,BBEGZ}.
Consequently, since $\bf 1$ is constant,
for any $C^2$ curve $\gamma(s)$ in $\Ho$,
\begin{equation}
\label{GeodCharEq}
\ddot L(\gamma(s))=\gM(1,\nabla^{\gM}_{\dot\ga}\dot\ga),
\end{equation}
with the notation $\dot f:=df/{ds}, \ddot f:=d^2f/{ds^2}$.
That is, if $\nabla^{\gM}_{\dot\ga}\dot\ga\ge0$ then $\gamma(s)$ is a geodesic iff $L(\gamma(s))$ is linear in $s$.
Note also that $\iota_\o(\HO)=L^{-1}(0)$.

\subsection{Nonpositive curvature: the Calabi--Tian conjectures}
\label{NonposSubSec}

As reviewed in \S\ref{KEEqSec} the early history of \KE metrics
revolved around the local version of this equation.
Two decades later, Calabi first formulated an ambitious program for
constructing KE metrics on {\it closed} manifolds. He explicitly
formulated in writing the case of zero Ricci curvature, but in analogy
with the uniformization theorem, he expected that the negative case
should then follow as well. The following statement was first formulated
as a theorem in \cite{Calabi54},
before the existence part quickly became a conjecture
thanks to discussions between Calabi and Nirenberg on the need
of a priori estimates.

\begin{conj}
{\rm (Calabi's conjecture 1953 \cite{Calabi54,Calabi54ICM,Calabi1957})}
Let $(M,J)$ denote a closed \K manifold with $c_1(M,J)<0$ or
$c_1(M,J)=0$. Then $(M,J)$ admits a \KE metrics that are unique
up to homothety in the former case, and unique in each \K class
in the latter.
\end{conj}

When $\mu<0$, the uniqueness is immediate from the maximum
principle: if $u,v$ satisfy $\o_u^n/\o_v^n=e^{|\mu|(u-v)}$
then if $u-v$ is maximized in $p\in M$, the form $\i\ddbar(u-v)$
has non-positive eigenvalue with respect to $\o_v$, and
so $u-v\le0$; by symmetry $v-u\le0$, so $u=v$.
When $\mu=0$ the uniqueness part was proved by Calabi \cite{Calabi1957}
by exploiting the algebraic structure of the \MA equation
(more specifically, properties of determinants and integration
by parts). Indeed, if $\o_u^n=\o_v^n$ then
$0=\o_u^n-\o_v^n=\i\ddbar(u-v)\w T$ with $T$ a positive $(n-1,n-1)$-form.
Multiplying by $u-v$ and integrating by parts shows $u-v$ is constant.

The existence part of Calabi's conjecture was established by Aubin in the negative case
and by Yau in both cases \cite{Aubin1978,Yau1978} (the zero case
under the restrictive assumption that the manifold admits a reference \K metric with nonnegative bisectional
curvature was established earlier by Aubin \cite{Aubin1970}).
The main innovation was to establish
the a priori $C^0$ and Laplacian estimates conjectured by Calabi and Nirenberg.
The higher derivative estimates then followed by work of Calabi
described in \S\ref{Calabi3rdSubSec}, and by elliptic bootstrapping.

Four decades later, motivated by application to algebraic geometry,
Tian formulated a generalization for pairs of Calabi's conjecture.

\begin{conj}\label{TianConj}
{\rm (Tian's conjecture 1994 \cite{Tian1994})}
Let $(M,J)$ denote a closed \K manifold and $D=D_1+\ldots+D_r\subset M$ a divisor
with simple normal crossing support. Suppose that
$c_1(M,J)-\sum_{i=1}^r(1-\be_j)[D_j]$ is negative or zero.
Then $(M,J)$ admits a KEE metric with angle $2\pi\be_j$
along $D_j$ that is unique
in the former case, and unique in each \K class
in the latter.
\end{conj}

This conjecture, similarly to Calabi's original conjecture, was also first stated as
a theorem in Jeffres' Ph.D. thesis \cite{Jeffresthesis,Tian1994}, but subsequently only the uniqueness
was published \cite{Jeffres2} due to the need of a priori estimates, as well as
a good linear theory. A program outlining a collaboration between Jeffres and Mazzeo
toward such a linear theory  was announced by Mazzeo \cite{Mazzeo1999}. As described in
Theorem \ref{DsThm}, two independent sets of linear estimates, in the case of
a smooth divisor,
were obtained, finally, first by Donaldson \cite{Donaldson-linear}, and subsequently
in \cite{JMR}.
What was lacking from the approach suggested in \cite{Mazzeo1999} was the simple, but crucial,
observation of \cite{Donaldson-linear} that if one simply considers \Holder control only 
of the complex $(1,1)$-type
derivatives, then a Schauder theory can be established; if one considers the
full Hessian then this is not the case, due to the harmonic function $z^{1/\beta}$.
In the case of a smooth divisor, this conjecture was resolved  by
Mazzeo and the author \cite{JMR}, that also established higher regularity of the metric.
The key here was, first, to establish the a priori estimates in the absence of curvature bounds
and, second, to develop a good linear theory in the edge spaces that, among other things, allows
to obtain higher regularity (polyhomogeneity).
A different proof of existence was later given by Guenancia--P\v aun \cite{GP} relying on
work of Guenancia--P\v aun--Campana \cite{CGP} that appeared around the same
time as \cite{JMR} (also around the same time, Brendle \cite{Brendle} obtained the special
case of a smooth divisor and angle $\be\in(0,1/2]$ using Donaldson's linear
estimate, Theorem \ref{DsThm} (i) with $s=w$).
As observed by Datar--Song \cite{DatarSong}, the existence in the general snc case actually
follows easily by combining the Chern--Lu approach of \cite{JMR} and a standard regularization
argument \`a la Demailly; a more
computationally involved proof by Guenancia--P\v aun of the general case
appeared slightly earlier than \cite{DatarSong}, and around the same time also another approach based on regularization
was developed by Yao \cite{Yao} . Finally, the general snc
case has also been settled independently by
\cite{MR2012,MR},
where, in addition, the higher regularity of solutions is proved
by establishing a linear theory also in the snc case; such a theory
is considerably more complicated than in the smooth divisor case.

Much of Sections \ref{KahlerSection}, \ref{RCMSec}, and  \ref{APrioriSubSec}
are devoted to describing the proof of Conjecture \ref{TianConj}.
In fact,
the proof of Conjecture \ref{TianConj} is, in fact, a corollary
of the proof of the parallel result in the harder positive case---Theorem
\ref{KEEExistenceThm}. The strategy of proof, following \cite{JMR}, is described
in \S\ref{PositiveSubSec} below, where the latter result is stated.
We also survey other approaches in \S\ref{OtherApprSubSec}.

Fulfilling one of the original motivations
for introducing edges---see \S\ref{KahlerEdgeIntroSubsec},
the following result was first obtained in 2011 in \cite{JMR} (an
alternative proof was given in 2013 by Guenancia--P\v aun \cite{GP}):

\begin{cor}
\label{KEEallProjCor}{\rm (KEE metrics on all projective manifolds)}
Every projective manifold admits a KEE metric of negative curvature
for any angle $2\pi\be\in(0,2\pi)$.
\end{cor}

In fact, on any projective manifold $M$ there exists an ample class $H$.
Therefore, given $\be\in(0,1)$, for large enough $m=m(\be)\in\NN$,
both $|mH|$ admits a smooth representative by Bertini's theorem, and
the class $K_M+(1-\be)m_\be H$ is positive. Note that if $M$ is {\it minimal},
i.e., $K_M$ is nef, then one may take $m_\be=1$ independently of $\beta$
by Kleiman's criterion, so that the fixed pair $(M,H)$ admits
KEE metrics for all $\be\in(0,1)$.

\begin{problem}
Let $M$ be minimal and $H$ a smooth ample divisor.
Describe the limiting behavior of the KEE metrics on $(M,H)$ as $\be$ tends to 1.
\end{problem}

\begin{example}
{\rm
To illustrate some of the new metrics arising from Corollary \ref{KEEallProjCor}
consider the projective plane $\PP^2$ which is typically thought of as `positively curved'.
By choosing a smooth plane curve $D$ of degree $d\ge 4$ and any
$\be\in(0,(d-3)/d)$, one can thus construct a KEE metric of negative
Ricci curvature. As $d\ra\infty$ the angle $2\pi\be$ can be made arbitrarily close
to $2\pi$. Of course, this example can be generalized to any $\PP^n$,
with $D$ a smooth hypersurface of degree $d$, and $\be\in(0,(d-1-n)/d)$.
}\end{example}

\begin{problem}
As $d\ra\infty$, does the divisor $D$ becomes dense in $\PP^n$ with respect
to the KEE metric? How do the diameters of $D$ and $M$ behave?
What is a Gromov--Hausdorff limit of such a metric space 
(for different limiting values of $\be$)?
\end{problem}

The existence problem in the positive Ricci curvature realm is  more complicated,
as the cohomological assumption is not sufficient.
In the rest of this section we first describe obstructions to existence,
and then, finally, formulate a general existence criterion.

\subsection{Reductivity of the automorphism group}
\label{ReductivitySubSec}

The (connected) isometry group of the round 2-sphere complexifies to
give the (identity component of the) M\"obius group of its conformal transformations.
A similar fact also holds for the isometry group of a football:
it complexifies to the conformal transformations that fix the poles,
or cone points.
Matsushima's theorem states that the same is true
for any KEE manifold.

\begin{thm}
\label{reductiveThm}
Let $(M,J,D,g)$ be a KEE manifold.
Let $\h{\rm Isom}_0(M,g)$ denote the identity component of the isometry group and denote by
$\Aut_0(M,D)$ the identity component of the Lie group of automorphisms preserving $D$. Then
$\Aut_0(M,D)=\h{\rm Isom}_0(M,g)^\CC$. In particular,
$\Aut_0(M,D)$ is reductive.
\end{thm}

Matsushima's original result was proved for KE manifolds
and extended to constant scalar curvature manifolds
by Lichnerowicz \cite{Matsushima,Lichnerowicz}. The
singular version appearing in \cite{CR} follows the same lines, by
using regularity results from \cite{JMR}. A more general result (assuming
the manifold only has a weak KE metric in the sense of \cite{EGZ})
but in the special case $M$ is Fano and $D\in|-K_M|$ can be found in
\cite{CDSun,Tian2013}.

Instead of reviewing the proofs above, we describe a formal proof
due to Donaldson in the smooth setting, in the spirit of
infinite-dimensional geometry \cite{Donaldson-CIMAT}.

First, let $\h{\rm diff}(M,D)$ denote the Lie algebra of
vector fields on $M$ that vanish on $D$, and denote
by $\hamMD$ the Lie subalgebra of fields $X$ satisfying
$[\iota_X\o]=0$. The corresponding Lie subgroup, $\HamMD$, consists of
 Hamiltonian diffeomorphisms of $(M\sm D,\o)$ that preserve $D$.
 The group $\HamMD$ acts in a natural
way on the space $\calJ$ of $\o$-compatible complex structures.
Furthermore, $\calJ$ comes equipped with a natural symplectic
structure. A result of Donaldson in the smooth setting shows
that the $\HamMD$-action on $\calJ$ is Hamiltonian, with a
moment map given by the scalar curvature $s(\o,J)$ minus
its average.
As long as we restrict to diffeomorphisms $F$ that preserve polyhomogeneity
and for which further $F^\star \o-\o=\i\ddbar u_F$, with $u_F\in\Dw$,
it is not hard to show that this result extends to the edge setting.

Now, suppose a holomorphic Hamiltonian group action $\calG$
can be complexified to $\calG^\CC$.
A basic fact in finite-dimensional moment map geometry says that
for an element belonging to the zero of the moment map, the isotropy group of
 the $\calG^\CC$-action is the complexification
of the isotropy group of the $\calG$-action.
Now for any $J\in \calJ$ the isotropy group of $\HamMD$ are those
diffeomorphisms of $M\sm D$ that preserve both $J$ and $\o$, i.e., isometries.
Similarly, if one is willing to think of $\HamMD^\CC$ as all diffeomorphisms
preserving $D$ and the (1,1)-type of $\o$, then the isotropy group
of $J$ can be identified with $\Aut(M,D)$ (some motivation for
such a formal identification is discussed in \cite{Semmes92,Donaldson99}).
Thus, if $J\in\calJ$ is a constant scalar curvature structure, the
Matsushima--Lichnerowicz criterion follows.

\begin{example}
\label{CR}
{\rm
The automorphism group of all strongly asymptotically log del
Pezzo pairs (these pairs are defined in Definition \ref{definition:log-Fano})
are computed in \cite{CR}. Here are some explicit examples.
The
pair $(\PP^n, H)$ with $H$ a hyperplane in $\mathbb{P}^n$,
satisfies
$$
\Aut(\PP^n, H)\cong\Aut(\PP^n,
p)\cong\Aut(\h{Bl}_p\PP^n)\cong\mathbb{G}_a^n\rtimes\mathrm{GL}_{n}(\mathbb{C}),
$$
for a point $p\in\mathbb{P}^n$, where $\h{Bl}_p\PP^n$ denotes the
blow-up of $\PP^n$ at $p$, and where $\mathbb{G}_a$ the additive
group $(\mathbb{C},\,+\,)$.
 The latter group is not reductive. This generalizes the
well-known obstruction to the existence of a
constant curvature metric on the teardrop ($S^2$ with one cone
point), see \S\ref{RSSubSec} below (another way to see this pair
is obstructed is to use the Bogomolov--Miyaoka--Yau inequality \cite{Tian1994,SongWang}).
For the pair $(\PP^2,D)$, with $D\in|2H|$ a smooth quadric,
any automorphism of $D$ lifts to automorphism of $\PP^2$,
since $\PP^2(H^0(D,-K_D))=\PP^2$ (again, note that $\Aut(\PP^2)\ra\Aut(D)$
is injective, since an automorphism of $\PP^2$ can only fix a linear
subspace).
Thus $\Aut(M,D)$ equals $\mathrm{PGL}_{2}(\mathbb{C})$ (and is reductive).
This pair, however, turns out to be obstructed only when  $\be\in(0,1/4]$
\cite{LiSun}.
}
\end{example}

\subsection{Mabuchi energy, Futaki character, and their relatives}
\label{MabuchiFutakiSubSec}

When \eqref{c1Eq} holds, the KE problem reduces to finding
a constant scalar curvature metric in $\Ho$. Indeed,
$\dbar s_\o=\del^\star\Ric\,\o$, where $s_\o=\tr_\o\Rico$,
denotes the scalar curvature (more precisely, half of the standard convention in
Riemannian geometry).
Or simpler, $s_\o\o^n=n\Ric\o\w \o^{n-1}$ and thus $s_\o$ is constant
iff $\Ric\o$ is harmonic, i.e., a multiple of $\o$.
This can be considered as an easy K\"ahlerian analogue
of the Obata theorem in conformal geometry.

Thus, it is natural to consider the following vector field on $\Ho$:
\beq\label{svecEq}
{\bf s}:\o\mapsto s_\o-\gM(s_\o,1)=s_\o-nc_1.[\o]^{n-1}/[\o]^n.
\eeq

\subsubsection{Mabuchi K-energy}
The zeros of the vector field {\bf s} are
the constant scalar curvature (csc) metrics in $\calH$, and its integral curves are
the trajectories of the Calabi flow \cite{Calabi1982}.
A remarkable fact is that ${\bf s}$ is in fact a
gradient vector field $\nabla^M E_0={\bf s}$ \cite{Mabuchi1986}.
Its potential $E_0:\Ho\times\Ho\ra\RR$
$$
E_0(\o_0,\o_1)=\int_\gamma{\bf s}^{\flat},
$$
is known as the Mabuchi K-energy,
where ${\bf s}^\flat$ denotes the 1-form associated to $\bf s$ via
$\gM$, and $\gamma$ is any simple path between $\o_0$ and $\o_1$.

We now mention several fundamental properties of Mabuchi's K-energy
due to Bando--Mabuchi \cite{BandoMabuchi}.
Any critical point of $E_0$ is KE as already remarked; and in fact
all such critical points lie in a connected finite-dimensional
totally geodesic submanifold of $\HO$ parametrizing all KE metrics
in $\HO$ and isometric to the symmetric space
$\Aut(M,J)/\h{Isom}(M,\oKE)$---the orbit of a single KE metric $\oKE\in\HO$
under the action of the identity component of $\Aut(M,J)$
on $\HO$ by pull-back.
Moreover, a computation shows that the second variation of $E_0$
at a critical point is nonnegative, so KE metrics
are local minima. A more difficult result is that KE metrics
are in fact global minima of $E_0$.
Thus, the KE problem is equivalent to determining when $E_0$
attains its absolute minimum.
In particular, for a solution to exist,
$E_0$ must be bounded from below.
In practice, it is easier to use this criterion to prove non-existence,
as we discuss in the next paragraph.
A stronger criterion introduced by Tian \cite{Tian1997}, ``properness of $E_0$", does turn out to
be equivalent to existence; see \S\ref{EnergySection}.

\subsubsection{Futaki character}

A na\"\inodot ve way to show that $E_0$ is unbounded from below is to find
a path along which its derivative is uniformly negative.
The simplest kind of path arises from a holomorphic vector
field, and is the pull-back of a
fixed metric by a one-paramter group of automorphisms.
To spell this out, denote by ${\psi^X}$ the vector field on $\Ho$ associated to the
holomorphic vector field $X$,
$$
{\psi^X}:\o\mapsto \psi^X_\o\in C^\infty(M),
$$
with $\calL_X\o=d\iota_X\o=\i\ddbar\psi^X_\o$, and $\int\psi^X_\o\o^n=0$.
Thus, ${\psi^X}$ is tangent to $\iota_\o(\HO)$. Moreover,
if $\ga(s)$ is an integral curve of ${\psi^X}$ in $\iota_\o(\HO)$,
then $\o(s):=\iota_\o^{-1}\ga(t)=(\exp sX)^\star\o(0)$ and
$\dot L(\ga(s))=\gM({\bf 1},\psi^X)|_{\o(s)}
=\int\psi^X_{\o(0)}{\o(0)}^n=\gM({\bf 1},\psi^X)|_{\o(0)}$
since $\psi^X_{\o(s)}=(\exp sX)^\star\psi^X_\o(0)$ from the definition of
$\psi^X$ and the fact that the pull-back by an automorphism commutes with
$\ddbar$. Thus, $\ga(s)$ is a geodesic by \eqref{GeodCharEq}.

Along such geodesics 
$\dot E_0(\o(0),\o(s))=\gM({\bf s},\psi^X)(\ga(s))$ is manifestly constant in $s$.
What is less obvious is that this latter constant
does not depend on the initial metric $\o(0)$.
To see this, note
that $\bf s^\flat$$=dE_0$ is evidently a closed form
(here $\flat$ and
$\sharp$ denote the ``musical operators" associating to a vector field
a one-form via the metric $\gM$ and vice versa).
Denote by $\exp t\psi^X$ the flow on $\Ho$ associated to $\psi^X$,
given by pull-back by $\exp tX$. Since the scalar curvature is natural
then $(\exp t\psi^X)^\star {\bf s^\flat}={\bf s^\flat}$.
Combining these facts, $\calL_{\psi^X}{\bf s^\flat}=d\iota_{\psi^X}{\bf s^\flat}
=0$. Thus, $\psi^X(E_0)$ is constant on $\Ho$, as claimed.
It is thus denoted by $F(X)$, and is called the Futaki invariant of $X$
\cite{Futaki1983,Calabi1985,Bourguignon1986}.
Since
$F(-X)=-F(X)$, it follows that $F:\aut(M,J)\ra\RR$ must vanish
identically or else $E_0$ is unbounded from below.

Like most objects in \K geometry, also the above discussion
generalizes in a straightforward manner to include edges when \eqref{cohomEq}
holds. Considering the 1-form
$$
{\bf s_\be^\flat}:\o\mapsto
(s_\o-nc_1.[\o]^{n-1}/[\o]^n)\on-(1-\be)\on|_D,
$$
proofs conceptually similar to the ones in the smooth setting show that this 1-form is
closed, admits a potential $E_0^\be$, and an associated Futaki
character on the Lie algebra $\aut(M,D)$ of holomorphic vector
fields vanishing on $D$, and that $E_0^\be$ must be bounded from
below when a KEE metric exists. To obtain these results one
works on the space 
$\calH^\eps_\o$ \eqref{KEdgePotentialsEtaSpace},
where integration by parts arguments, as in the smooth setting,
are justified.

We refer to \cite{Futaki,DingTian1992,Lu1999,Tianbook,Nakagawa,ChiLithesis}
for explicit computations of Futaki invariants.

\begin{remark}
\label{FutakiMatsushimaRemark}
{\rm
In the smooth 2-dimensional Fano setting Futaki's and Matshusima's obstructions
coincide. But in higher dimensions, there exist examples where
only one of the obstructions appears. For instance, $\PP(E)$,
where $E=\calO_{\PP^1}(-1)\oplus\calO_{\PP^2}(-1)$ is the rank 2 bundle
over $\PP^1\times\PP^2$, is a Fano 4-fold with reductive automorphism
group and nonvanishing Futaki character \cite[pp. 24--26]{FutakiMabuchiSakane}
(see also \cite{Futaki1983,Wang1991}).
On the other hand, there exists a Fano 3-fold whose connected automorphism
group is $\mathbb{G}_a=(\CC,\,+\,)$, see case 2) in the main result of Prokhorov \cite{Prokhorov},
and so its Futaki character must vanish by a theorem of Mabuchi \cite[Theorem 0.1]{Mabuchi-poisson}.
We also remark that for some time it was believed that these two obstructions
should be also sufficient \cite[\S4]{Futaki1983},
\cite[p. 575]{Mabuchi1986}, \cite[Conjecture E]{HwangMabuchi1993}. This was verified by Tian
for del Pezzo surfaces \cite{Tian1990} (see \S\ref{logCalabiSec}),
but disproved for 3-folds, as we discuss next.
}
\end{remark}

\subsubsection{Degenerations and geodesic rays}

One-parameter subgroups of automorphisms
are particularly amenable to computations, as we saw in the previous
subsection. Yet, generic complex manifolds do not admit
such automorphisms, and the complexity of $E_0$ goes
beyond automorphisms: Tian constructed Fano 3-folds
with no nontrivial one-parameter subgroups of automorphisms that
admit no KE metrics \cite{Tian1997}.
As a replacement, Tian suggested to consider $M$ as embedded in
a one-parameter family of complex manifolds and consider
a $\CC^\star$ action on this whole family. He observed that
one may still associate a Futaki type invariant to this family,
often referred to as a {\it special degeneration} or {\it
test configuration}, and this reduces to the ordinary invariant when the family is a product.
(Of course, we are glossing over many technical details here, among them
that the action should lift to an action on a polarizing line bundle,
and that the invariant is computed by considering high powers of
this bundle. As noted in the Introduction, we refer to Thomas \cite{Thomas}
for GIT aspects of KE theory.)
Furthermore, he conjectured, in what became later known as
the Yau--Tian--Donaldson conjecture, that if the sign of this
invariant is identical for all special degenerations,
and is zero only for product configurations, then the manifold
should admit a KE metric \cite{Tian1997}. This criterion is
called {\it K-polystability}. In Tian's
original definition only certain singularities were allowed
on the ``central fiber" of such a degeneration. Donaldson
extended this to much more general ones \cite{Donaldson2002Toric}. Li--Xu finally showed
that the original definition suffices \cite{LiXu}.
A different definition of stability has been introduced
by Paul \cite{Paul2012,Paul2013,Tian-Paulsurvey2013}, who also conjectured
its equivalence to the existence of KE metrics.
Very recently, a solution to these
conjectures has been announced by Chen--Donaldson--Sun
and Tian \cite{CDSun,Tian2013}, crucially building upon the theory
of KEE metrics described in this article, and in particular on \cite{Donaldson-linear,JMR}.
Shortly after the appearance of  \cite{CDSun,Tian2013},
Sz\'ekelyhidi observed that those may be adapted
to give similar conclusions without using KEE metrics 
\cite{Szekelyhidi2014}.

Recently, Ross--Witt-Nystr\"om introduced the notion of an
analytic test configuration \cite{RWN}. Very roughly, in their approach
the singularities of the central fiber are replaced
by a ``singularity type" of a curve of $\o$-psh function, and
a (usually singular) generalized geodesic ray is constructed out of this data.
This generalized previous constructions of Arezzo--Tian, Phong--Sturm, Song--Zelditch,
that constructed such rays out a degeneration in the sense of the previous paragraph
\cite{ArezzoTian,PhongSturm,SongZelditch}.
A generalized, or weak, geodesic existing for time $s\in[0,T]$,
is a solution to the homogeneous complex
Monge--Amp\`ere equation on the product $S_T\times M$, where
$M$ is a strip $[0,T]\times\RR$ of width $T$,
but it need not be a path in $\calH$,
but only in $\PSH(M,\o)$. These can be regarded as
bona fide geodesics if one enlarges the space of \K metrics
appropriately, as considered by Darvas \cite{Darvas}.
A different approach to construction of geodesic rays is
suggested by Zelditch and the author via the Cauchy problem
for the \MA equation \cite{RZ123}.

Coming full circle, a conjecture of Donaldson states that the
existence problem should actually be completely characterized by
generalizing Futaki's criterion to all geodesic rays.
Namely, the non-existence of a csc \K metric is conjectured to
be equivalent to the non-existence of a geodesic
ray (whose regularity is not specified in the conjecture) along
which the derivative of the K-energy is negative \cite{Donaldson99}.

\subsubsection{Flow paths and metric completions}

Another conceivable way to ``destabilize" a Fano manifold is
to use Hamilton's (volume normalized) Ricci flow \cite{Hamilton},
\begin{equation}
\label{RFEq}
\frac{\del \o(t)}{\del t}=-\Ric \o(t)+\o(t),\quad \o(0)=\o\in\calH_{c_1},
\end{equation}
which preserves
the space $\calH_{c_1}$ of K\"ahler forms cohomologous to $c_1$
and exists for all time \cite{Cao}.
A theorem of Perleman asserts that
the flow will converge to a \KE metric if and only if such
a metric exists \cite{TianZhu2007}. Thus, it is tempting
to hope that Ricci flow trajectories could actually
give another way of constructing geodesic rays. The following
conjecture was suggested by La Nave--Tian  \cite[\S5.4]{LaNaveTian},
based on a description of the Ricci flow as a Monge--Amp\`ere equation
in one extra dimension, that should bear a relationship to
the homogeneous complex \MA equation governing geodesics in the Mabuchi metric
\eqref{MabuchiMetricEq}.

\begin{conj}
\label{KRFConj}
The \KRF is asymptotic to a
$\gM$-geodesic in a suitable sense.
\end{conj}

The first attempt to understand the metric completion of $\calH_\o$
is due to Clarke and the author \cite[Theorem 5.6]{ClarkeR} where the metric completion
with respect to the Calabi metric is computed, motivated
by an old research announcement of Calabi \cite{Calabi54} that speculated
a different answer. Calabi's metric is the $L^2$ metric on the level of
\K forms, that also takes the form
\begin{equation}
\label{CalabiMetricEq}
\gC(\nu,\eta)|_\vp:=\int_M\D_\vp\nu\D_\vp\eta\,\frac{\o_\vp^n}{n!}.
\end{equation}
Interestingly, Calabi's metric was his original motivation for introducing the Calabi
conjecture \cite{Calabi54}, see \cite[Remark 4.1]{ClarkeR}.
The following result relates degenerations arising from the Ricci flow,
the existence problem of KE metrics, and the metric geometry of $\calH_\o$.

\begin{thm}
\label{RFMetricThm}
Let $(M,J)$ denote a Fano manifold.
The following are equivalent:\hfill\break
(i) $(M,J)$ admits a KE metric;\hfill\break
(ii) Any Ricci flow trajectory in $\calH_{c_1}$ converges in the Calabi metric
$\gC$ \eqref{CalabiMetricEq}.
\end{thm}

Note that the derivative of the K-energy is negative
along the Ricci flow. Also, any divergent path must have
infinite length in the corresponding metric.
Thus, this result stands in precise analogy to Donaldson's
conjecture on ``geodesic stability," with Ricci flow paths taking the place of
$\gM$-geodesic rays. This gives some intuitive motivation both
for Donaldson's conjecture and for Conjecture \ref{KRFConj}.
We conjecture that these statements are also equivalent
to the statement that any Ricci flow trajectory in $\calH_{c_1}$ converges in the Mabuchi metric
$\gM$ \eqref{MabuchiMetricEq}. A weaker statement was
obtained by McFeron \cite{McFeron} where ``converges" is replaced
by ``has finite length".

Convergence in the metric sense is
rather weak from a PDE point of view. Thus,
the main point in proving this theorem is to show
that such convergence implies strong convergence.
A sequence of \K metrics $\o_i$ converges in the Calabi metric
precisely when the volume forms $\o_i^n$ converge in $L^1(M,\o_0^n)$ \cite[Corollary 5.5]{ClarkeR},
and in particular do not charge Lesbegue null sets.
Now, the limit along the Ricci flow converges to a metric that charges
an analytic subvariety \cite[Lemma 6.5]{ClarkeR} (this builds
on results of Nadel \cite{Nadel,Nadel1996} and the author \cite{R09}), which proves that
(ii) implies (i). The converse is an immediate corollary of the
exponential convergence of the flow when a KE metric exists
\cite{TianZhu2007,PSSW2}.

Finally, we remark that both Donaldson's conjecture and Conjecture \ref{KRFConj}
should be related to the Hamilton--Tian conjecture,
stipulating that the K\"ahler--Ricci
flow on Fano manifolds should converge in a suitable sense to a K\"ahler--Ricci soliton
(with respect to a possibly different complex structure) with mild singularities
\cite[Conjecture 9.1]{Tian1997}. This relation could be related to the following problem.

\begin{problem}
Determine which $\gM$-geodesic rays (possibly non-smooth) come from
automorphisms of a manifold with a nearby complex structure.
\end{problem}

For smooth Riemann surfaces the Ricci flow converges
to a constant scalar curvature by results of Hamilton and Chow \cite{Hamilton1988,Chow1991},
and to a soliton in the case of orbifold Riemann surfaces (i.e., angles $\be_i$ of the form
$1/m_i$ with $m_i\in\NN$)  \cite{ChowWu}.
The Hamilton--Tian conjecture was also recently established in the setting of conical Riemann surfaces \cite{MRSesum} and
in the smooth 3-dimensional setting by Tian--Zhang \cite{TianZ2013}.
The smooth 2-dimensional case was previously known by Tian's uniformization
of del Pezzo surfaces relying on their classification (see Remark \ref{FutakiMatsushimaRemark}
and the introduction of \S\ref{logCalabiSec}) and
Tian--Zhu's generalization of Perleman's convergence result for the
K\"ahler--Ricci flow \cite{TianZhu2007}.

Finally, we mention that much of the discussion above
has an analogue for the Calabi flow \cite{Calabi1982}. Calabi conjectured
that the eponymous flow exists for all time and converges
to a csc or, in an appropriate sense, to an extremal metric.
The most general long time existence result known is due
to Streets \cite{Streets} though the question of
regularity of such weak flows is open (assuming curvature
bounds along the flow a smooth solution exists for
all time by Chen--He \cite{ChenHe}).
An analogue of Conjecture \ref{KRFConj} in this setting
is shown by Chen--Sun \cite{ChenSun}. It would be interesting
to understand an analogue of Theorem \ref{RFMetricThm}.
We mention in this context that
according to \cite[Theorem 4.1]{ChenHe}
(see \cite[\S3]{CarlottoChodoshRubinstein} for a recent alternative and conceptual proof) 
the Calabi flow converges exponentially fast as soon as it converges smoothly, say.
In addition, a result of Berman states
that, when the K\"ahler class is a multiple of the canonical class,
the Calabi flow, when it exists, 
converges to a K\"ahler--Einstein metric when one exists
\cite[Theorem 1.4]{Berman2010}. Finally, results of 
Chen--He, Feng--Huang, He, Huang, and Tosatti--Weinkove, among others,
give various conditions for convergence of the Calabi flow
\cite{ChenHeCalabiSurf,FengHuang,He-localCalabi,Huang,TosattiWeinkove}.

\subsection{Conical Riemann surfaces}
\label{RSSubSec}
It is interesting to contrast the general results surveyed so far
with the 1-dimensional picture: when can one
uniformize a conical Riemann surface? By uniformization
we refer to the construction, in a given conformal class,
 of a constant scalar curvature metric
with prescribed conical singularities. This question was studied,
in the negative case, already in Picard's work more than a century ago \cite{Picard1893,Picard1905}
and was first treated definitively in Troyanov's thesis \cite{Troyanov-thesis,Troyanov} (see
also McOwen \cite{McOwen,Mc-correction}) where a sufficient condition for existence
was established; uniqueness and necessity of this condition
were addressed by Luo--Tian \cite{LuoTian}.
We refer to \cite{MRSesum} for background and an alternative approach
via the Ricci flow (see also \cite{Yin} for the nonpositive cases).

\begin{figure}
\centerline{\includegraphics[width=1.79in,height=2.49in,keepaspectratio]{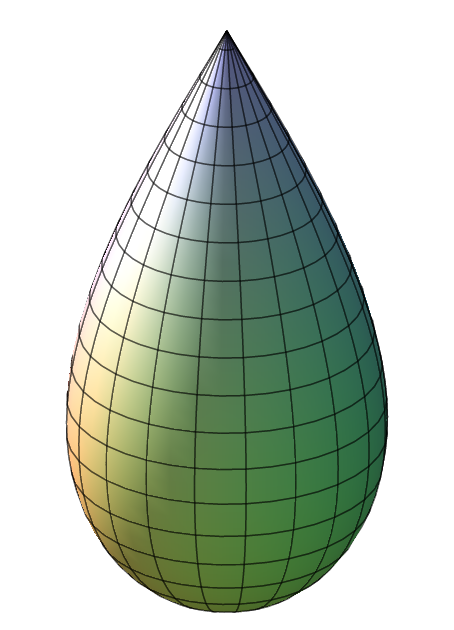}
\includegraphics[width=1.79in,height=3.5in,keepaspectratio]{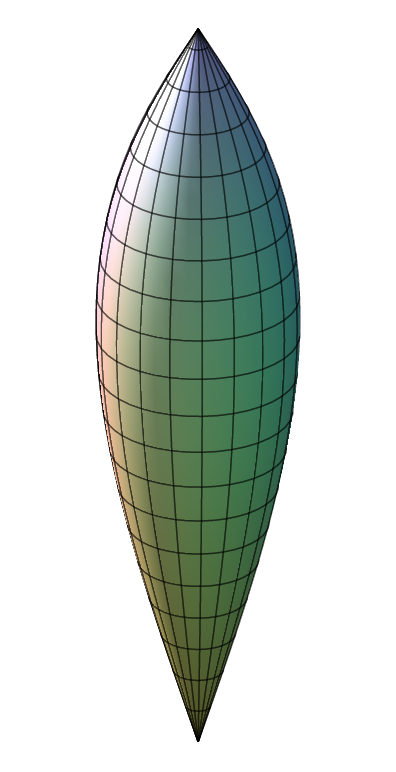}
\includegraphics[width=1.79in,height=3.5in,keepaspectratio]{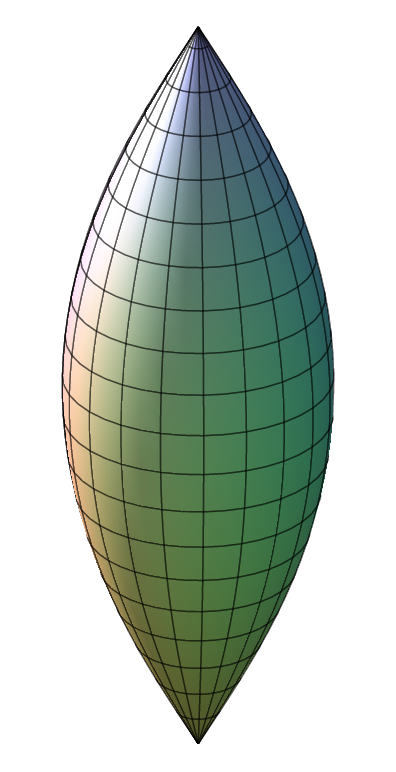}
}
\caption{The teardrop and football (with different angles) solitons with non-constant curvature, and
a constant curvature football with equal angles (courtesy of D. Ramos \cite{Ramos}).}
  \label{FigRiemannSurfaces}
\end{figure}

Fix a smooth compact surface $M$, along with a conformal, or equivalently, complex structure $J$.
A divisor $D$ is now a collection of distinct points $\{p_1, \ldots, p_N\} \subset M$
and the associated class $c_1(M,J)-\sum(1-\be_j)p_j$,
that can be thought of as a modified Euler characteristic, is
\begin{equation}
\label{cec}
2-2g(M)-\sum_{j=1}^N(1-\be_j).
\end{equation}
In particular, this is now a number, and so the cohomological condition \eqref{cohomEq}
is always satisfied. In a local conformal coordinate chart near each $p_j$,
$g=\i\gamma_j dz\otimes\overline{dz}=\i\gamma_j|dz|^2$,
with $\gamma=|z|^{2\be_j-2}F_j$ and $F_j$ bounded.
Suppose that $g$ has constant curvature away from the cone points.
The Poincar\'e--Lelong formula asserts
that $-\Delta_g\log|z|$ is a multiple of the delta function at $\{z=0\}$
(this can be seen by excising a small neighborhood near the cone point
and using Stokes' formula).
Together with the standard formula for the scalar curvature,
$K_g = -\Delta_g\log \gamma$ (up to a constant factor),
it follows that
\beq\label{CurvEq}
K_g-2\pi\sum(1-\be_i)\delta_{p_i}=\h{const},
\eeq
with the constant given by \eqref{cec}.
The existence in the nonpositive regime uses standard methods
as in the smooth setting,
e.g., the variational method of Berger \cite{Berger1971}, or the
method of sub- and super-solutions used by Kazdan--Warner \cite{KazdanWarner}.
A somewhat surprising discovery of Troyanov was the sufficient condition
\begin{equation}
\be_j-1 > \sum_{i \neq j}\be_i-1, \quad \h{for each $j=1,\ldots,N$,}
\label{TC}
\end{equation}
for existence in the positive case, that followed by generalizing
Moser's inequality \cite{Moser,Moser1973} to the singular setting.
For instance, if $N=2$, this is violated when $\be_1\not=\be_2$.
As observed by Ross--Thomas, this condition can be rephrased by
saying that the Futaki invariant of the pair $(M,\sum(1-\be_i)p_i)$
has the right sign, or that the pair is slope (poly)stable
\cite[Theorem 8.1]{RossThomasJDG} (see \cite[Remark 8.4]{RossThomasJDG}
for a careful treatment of the equality case in \eqref{TC}).
The only additional unobstructed pair is $N=2$ with $\be_1=\be_2$,
and such a csc metric can be constructed explicitly using ODE methods;
when $\be_1\not=\be_2$ or $N=1$ one can still construct a shrinking Ricci soliton
\cite{Hamilton1988} (see also \cite{BernsteinMettler,Ramos}). These are
depicted in Figure \ref{FigRiemannSurfaces}.
We also remark that \cite{JMR} and Berman's work gave a new proof of
Troyanov's original results \cite{Berman2010},
and Berndtsson's work gave a new approach to uniqueness \cite{Bern}.
Finally, the higher regularity of such a metric was only obtained much later \cite{JMR},
as a corollary of Theorem \ref{phgThm}.

The variational approach has recently been extended considerably through the work of Malchiodi et al.
to allow angles $\be>1$, even when coercivity fails, see, e.g., \cite{BDM,CM}.
For some higher regularity results also in this regime we refer to \cite{MRSesum}.

\subsection{Existence theorem for positive curvature}
\label{PositiveSubSec}

The following essentially optimal existence result in the positive case
is due to \cite{JMR}.
It parallels and generalizes Tian's theorem
from the smooth setting \cite[\S6]{Tianbook}.
We postpone the definition of properness to \S\ref{NonlinearSubSec}.

\begin{thm}
\label{KEEExistenceThm}
{\rm (\KE edge metrics with positive Ricci curvature)}
Let $(M,\o_0)$ be a compact \K manifold, $D\subset M$ a smooth divisor, and suppose that
$\be\in(0,1]$ and $\mu>0$ are such that
$$
c_1(M)-(1-\be)[D]=\mu[\o_0],
$$
and that the twisted K-energy $E_0^\beta$ is proper.
Then, there exists a \KE edge metric $\o_{\vpKE}$ with Ricci curvature $\mu$ and with angle $2\pi\be$ along $D$, that is unique up to automorphisms that preserve $D$.
This metric is polyhomogeneous, namely, $\vpKE$ admits a complete
asymptotic expansion  with smooth coefficients as $r \to 0$ of the form
\begin{equation}
\label{KEExpEq}
\vpKE(r,\th,z_2,\ldots,z_n) \sim \sum_{j, k \geq 0} \sum_{\ell=0}^{N_{j,k}} a_{ j k \ell }
(\th,z_2,\ldots,z_n)r^{j+k/\be}(\log r)^{\ell},
\end{equation}
where locally $D$ is cut-out by $z_1$, $r = |z_1|^\beta/\beta$ and $\th = \arg z_1$, and with each $a_{jk\ell} \in C^\infty$.
There are no terms of the form $r^\zeta (\log r)^\ell$ with $\ell > 0$ if $\zeta \leq 2$.
In particular, $\vpKE\in\calA^0 \cap \calD^{0,\frac1\be-1}_w$, i.e., $\o_{\vpKE}$  has infinite conormal regularity, and
is $(\frac1\be-1)$-H\"older continuous with respect to the reference edge metric $\o$.
\end{thm}

A parallel result in the snc case can be stated \cite{MR2012,MR}. In fact,
the a priori estimates of \cite{JMR} apply to the snc case without any change.
The essential new difficulty, however, as compared to the smooth divisor case,
is to extend the linear theory to one with crossing edges. This is new even
in the real setting, and goes beyond the original methods of Mazzeo \cite{Mazzeo}.
A different approach to the existence,
avoiding such linear theory,
but producing only $\calD^{0,0}_w$ solutions
 without
\Holder estimates on the metric (or even continuity of the metric up to $D$) nor
higher regularity,
and under the assumption that a $C^0$ (or even a weaker) solution exists,
has been developed by Guenancia--P\v aun \cite{GP} (a similar result is obtained by Yao \cite{Yao} who
derives the $C^0$ estimate based on the approximation scheme of \cite{CDSun,Tian2013}),
as described in \S\ref{NonposSubSec}.
Their approach starts from the $C^0$ solution constructed
by Berman \cite{Berman2010} and \cite{JMR}, and produces a $\calD^{0,0}_w$ estimate by
a careful approximation by smooth metrics\footnote{Just before posting this 
article, Guenancia--P\v aun have considerably revised their article \cite{GP}.
In this version of their article they also develop, among other things, 
a new alternative approach
for the $\Dw$ estimate.}.
Finally, the uniqueness is due to Berndtsson \cite{Bern}.

\begin{proof}[Strategy of proof of Conjecture \ref{TianConj} and Theorem \ref{KEEExistenceThm}]
The proof uses results from Sections \ref{KahlerSection}, \ref{RCMSec}, and
\ref{APrioriSubSec}. We now describe how these results piece together.

Solving the KE equation is equivalent to solving the \MA equation \eqref{KEEEq}.
We embed this equation in a one-parameter family of equations \eqref{RCMEq},
called the Ricci continuity path (Ricci CP), and define as usual the set
$$
I:=\{s\in(-\infty,\mu]\,:\,  \h{\eqref{RCMEq} with parameter value $s$
admits a solution in $\Dw\cap C^4(M\sm D)$}\}.
$$
Equation  \eqref{RCMEq} with parameter value $\mu$ is precisely \eqref{KEEEq}.
By Proposition \ref{NewtonProp}, there exists some $S>-\infty$ such
that $(-\infty,S)\subset I$, so that $I$ is not empty. Furthermore, $I$
is open. As a matter of fact, the linearization of  \eqref{RCMEq} at the parameter value $s$
is given by $\Delta_{\o_{\vp(s)}}+s$. This is clearly invertible when $s<0$, and is
also invertible when $s>0$ since, as shown in \S\ref{UniformSubSec}, the first
eigenvalue of $\Delta_{\o_{\vp(s)}}+s$ is positive whenever $s\in(0,\mu)$. Invertibility
at $s=0$ follows by working on the orthogonal complement of the constants \cite{Aubin1984}.
Thus, Theorem \ref{DsThm} (i) yields the openness. The fact that $I$ is closed (and hence
equal to $(-\infty,\mu]$) follows from Theorem \ref{APrioriEstTotalThm}.
Thus,   \eqref{KEEEq} admits a solution in $\Dw\cap C^4(M\sm D)$.
Theorem \ref{phgThm} implies that the solution is polyhomogeneous and belongs
to $\calA^0_\phg$; Equation \eqref{ExpansionAllBetaEq} further implies the precise
regularity statements about the solution.
\end{proof}

Note that the proof gives a new and unified proof for the classical results of Aubin, Tian and Yau,
on the existence of smooth KE metrics, in that it uses a single continuity method
for all signs of $\mu$. This
point is discussed in detail in \S\ref{OtherCMSubSec}.

Finally, consider the special case that $M$ is Fano and $D$ is a smooth anticanonical
divisor (the existence of such a divisor is highly nontrivial,
we refer to Problem \ref{FanoProb} and the discussion there).
Then, as noted by Berman \cite{Berman2010}, the twisted K-energy is proper for small $\mu=\be>0$.
Theorem \ref{KEEExistenceThm} thus gives the following corollary conjectured by Donaldson \cite{Don2009}.
\begin{cor}
\label{FanoCor}
Let $M$ be a Fano manifold, and suppose that there exists a smooth anticanonical divisor
$D\subset M$. Then there exists some $\be_0\in(0,1]$ such that for all $\be\in(0,\be_0)$
there exists a KEE metric with angle $2\pi\be$ along $D$ and with
positive Ricci curvature equal to $\be$.
\end{cor}

Li--Sun \cite{LiSun} observed that the exact same arguments actually prove existence
also in the plurianticanonical setting: supposing there exists a smooth divisor $D$ in $|-mK_M|$
(this always holds if $M$ is Fano and $m\in\NN$ is sufficiently large by Kodaira's theorem), there exists
a small $\be_0>0$ such that there exists a KEE metric
with $\be\in(1-\frac1m,1-\frac1m+\be_0)$ along $D$ and with $\mu=1-(1-\be)m$.

%%%%%%%%%%%%%%%%%%%%%%%%%%%%%%%%%%%%%%%%%%%%%%%

\section{Energy functionals}
\label{EnergySection}

The Mabuchi energy was crucial in stating the
main existence theorem (Theorem \ref{KEEExistenceThm}).
The purpose of this section is to describe several energy functionals, including
the Mabuchi energy, that
cast the \KE (edge) problem as a variational one.
(The variational approach can be applied to much more general settings,
and we do not seek full generality, for which we refer to
Berman et al. \cite{Berman2010,BBEGZ,BBGZ}.)
First, in \S\ref{NonlineaerEnergiesSubSec}, we describe the Aubin energy functionals, that
are nonlinear generalizations of the $W^{1,2}$-seminorm.
Using these functionals it is easy to construct the Ding functional
whose Euler--Lagrange equation is the inhomogeneous \MA equation.
Using the Aubin energies, one can then define a notion of relative
compactness of level sets of an energy functional, and this
leads to a calculus of variations formulation of the \KE problem
in \S\ref{NonlinearSubSec}.
Both the Mabuchi and Ding energies can be defined in a more conceptual
way, using Bott--Chern forms. We describe this in detail in
\S\ref{BottChernSubSec} since we are not aware of a single easy-to-read
reference for this (we refer to \cite{PRSturm,SongWeinkove} for the
related Deligne pairing formalism).
In \S\ref{StaircaseSubSec} we also
describe other natural functionals, the K\"ahler--Ricci energies, that lend themselves to a similar
description, mainly to illustrate the richness
of the theory, but also to show that there are many more-or-less equivalent
functionals whose variational theory underlies the KE problem.
Subsection \ref{LegendreSubsubsec} describes a relation between
the Ding energy and the Mabuchi energy, involving the Legendre
transform. Building on the preceeding
subsections, in \S\ref{EquivalenceSubSec} we prove the equivalence, in suitable senses, of
the Ding, Mabuchi, and K\"ahler--Ricci energies.

\subsection{Nonlinear Dirichlet energies and the Berger--Moser--Ding
functional}
\label{NonlineaerEnergiesSubSec}

The most basic functionals, going back to the work of Aubin
\cite{Aubin1984}, are defined by
\beq
\label{AubinEnergyEq}
\begin{aligned}
I(\eta,\etavp) & =\frac1V\int_M\i\del\vp\w\dbar\vp\w\sum_{l=0}^{n-1}\eta^{n-1-l}\w\etavp^{l}
=\frac1V\int_M\vp(\etan-\etavpn),
\cr J(\eta,\etavp) & =\frac{V^{-1}}{n+1}\int_M\i\del\vp\w\dbar\vp\w\sum_{l=0}^{n-1}(n-l)\eta^{n-l-1}\w\etavp^{l}.
\end{aligned}
\eeq
The functionals $I,J$, and $I-J$ are all equivalent
(and hence the latter is nonnegative), in the
sense that
$\frac1n J\le I-J\le \frac{n}{n+1} I\le nJ.$
Granted, these might not look so `basic' at a first glance. Let us
give some motivation to these definitions.

The first motivation comes from the calculus of variations:
is there a sort of nonlinear Dirichlet energy
whose Euler--Lagrange equation is \eqref{KEEq}? In
the case $n=1$ such an energy was studied by Berger and Moser
\cite{Berger1971,Moser1973}
\beq
\label{FOneDEq}
F(\eta,\etavp)
=
\begin{cases}
\dis\frac1V\int\frac12\i\del\vp\w\dbar\vp-\frac1V\int\vp\eta-\frac1\mu\log\frac1V\int e^{f_\eta-\mu\vp}\eta,
& \hbox{for $\mu\ne0$,}\cr
\dis\frac1V\int\frac12\i\del\vp\w\dbar\vp-\frac1V\int\vp\eta
+\frac1V\int \vp e^{f_\eta}\eta,
& \hbox{for $\mu=0$.}\cr
\end{cases}
\eeq
(Recall the definition of $f_\eta$ \eqref{foEq}.)
Its critical points are precisely constant scalar curvature metrics in the conformal class.
It is straightforward to generalize in higher dimensions the last term
to $\frac1\mu\log\frac1V\int e^{f_\eta-\mu\vp}\eta^n$
(respectively, $\frac1V\int \vp e^{f_\eta}\etan$)
so that its
variational derivative comes out to be the right hand side of \eqref{KEEq}
(up to the implicit normalization that the integral of the right
hand side is 1).
How to replace the first two terms so that the resulting Euler--Lagrange
equation comes out right?
This amounts to finding a functional whose differential is
exactly minus the \MA operator, i.e., the left hand side of \eqref{KEEq}.
This is precisely the functional
$-L$ (recall \eqref{LEq} and the line following it). Decomposing $-L$ into two terms,
one of which is $-\frac1V\int\vp\eta^n$, yields precisely
$J$ as the second term! (Incidentally, this also explains why
$L$, defined in \eqref{LEq}, is called the Aubin--Mabuchi functional.)

Thus, the Euler--Lagrange equation for the Berger--Moser--Ding energy
(or Ding energy for short) \cite{Ding1988}
\beq
\label{FbetaEq}
F^\be(\eta,\etavp)
=
\begin{cases}
\dis
J(\eta,\etavp)-\frac1V\int\vp\etan-\frac1\mu\log\frac1V\int e^{f_\eta-\mu\vp}\eta^n,
& \hbox{for $\mu\ne0$,}\cr
\dis
J(\eta,\etavp)-\frac1V\int\vp\etan+\frac1V\int e^{f_\eta}\eta^n,
& \hbox{for $\mu=0$,}\cr
\end{cases}
\eeq
is precisely \eqref{KEEEq}.

\subsection{A nonlinear variational problem}
\label{NonlinearSubSec}

This brings us to a second motivation for $J$. As we just noticed,
$-L(\vp)=J(\eta,\etavp)-\int\vp\etan$. 
Recall that $\pm L$ \eqref{LEq} is a distance function
for the Mabuchi metric.
Thus, it is tempting to think of
$J$ as a sort of approximate distance function.
This is of course not quite true, since the
submanifold
$\{\vp\in\calH^e_\o\,:\, \int\vp\etan=0\}$
of $\calH^e_\o$ is not totally geodesic.
But ignoring this subtlety, one is then tempted
to think of $J$
as a good way to define coercivity for our nonlinear
problem on $\calH_\o^e$; of course this temptation also
arises from the first motivation discussed earlier.
Following Tian \cite{Tian1997}, one says a functional $E$ on $\Ho\times\Ho$ is
proper provided that it dominates $J$
(or, by their equivalence, $I$ or $I-J$) on each $\Ho$ slice.
This is an analogue of the standard assumption
in the direct method in the calculus of variations, namely that sublevel
sets (of $E$) are compact.
In particular, if $E$ is proper, 
a sublevel set of $E$ is contained in some sublevel
set of $J$.
Theorem \ref{KEEExistenceThm} can be recast
as a nonlinear analogue of the fundamental
theorem of the direct method (cf. \cite[Theorem 1.1]{Struwe}) in this setting,
justifying the definition of properness.

\begin{thm}
\label{PropernessThm}
Suppose that $F^\be(\o,\,\cdot\,)$ is proper on $\calH^e_\o$.
Then it is bounded from below, and its infimum is attained
at a solution of (\ref{KEEq}).
\end{thm}

For the proof we refer to \cite[Theorem 2]{JMR} (that proof assumes that
the Mabuchi functional $E^\be_0$ is proper, however the arguments are
identical; a conceptual proof of the equivalence coercivity of $F^\be$
and that of $E_0^\be$ is given in Theorem \ref{EquivLowerBdThm} (ii),
although a direct proof (i.e., one that does not use the existence
of a minimizer) of equivalence of the properness of these functionals
seems to be unknown).
The special case when $\be=1$ goes back
to Ding--Tian \cite{DingTian1993} and Tian \cite{Tian1997} (see
also \cite[\S2.2]{Tian2012} for a generalization that allows for automorphisms).
On the other hand, Berman showed that under the
assumption in the theorem the infimum is attained
within the larger set $\PSH(M,\o)\cap C^0(M)$,
and while it is easy to see that the minimizer
is contained in $C^\infty(M\sm D)$ by local ellipticity
and the usual arguments as in the smooth case,
the proof that the minimizer lies in $\calH^\eps_\o$ \cite{JMR}
relies on the Ricci continuity method,
as described in \S\ref{RCMSec}--\ref{APrioriSubSec}.

\subsection{The staircase energies and K\"ahler--Ricci energies}
\label{StaircaseSubSec}

A generalization of the Aubin energy $J$ \eqref{AubinEnergyEq}
was introduced in \cite[(6)]{RJFA},
\beq
\label{IkEq}
I_k(\o,\ovp)=
\frac{V^{-1}}{k+1}\int_M\vp(k\on-\sum_{l=1}^k\o^{n-l}\w\ovp^l).
\eeq
Note that $I_n=J,\; I_{n-1}=({(n+1)J-I})/n$. These functionals
are nonnegative (see \eqref{IkSecondEq} below) and
are related to, but different from, functionals defined by
Chen--Tian \cite{ChenTian2002}---see \cite[p. 133]{RThesis}.
The functionals $I_k$ can be thought of as gradual nonlinear generalizations
of the Dirichlet energy that interpolate between the latter and $J$.
In fact, $I_1$ is simply a multiple
of the Dirichlet energy $\int \i\del\vp\w\bar\del\vp\w\o^{n-1}$,
and increasing $k$ is analogous to climbing a staircase
(see \cite[p. 132]{RThesis} for a pictorial description)
where $I_k$ incorporates one additional `mixed'
term proportional to $\int\i\del\vp\w\bar\del\vp\w\o^{n-k}\w \ovp^{k-1}$;
indeed, by integrating by parts,
\beq
\label{IkSecondEq}
I_k(\o,\ovp)=
\frac1V\int_M
\i\del \vp\w\dbar \vp\w\sum_{l=0}^{k-1} \frac{k-l}{k+1}\o^{n-1-l}\w\ovp^l.
\eeq

The definition \eqref{IkEq}
certainly makes sense for pairs of smooth \K forms, and by
Theorem \ref{DsThm3} (i) and
the continuity of the mixed
\MA operators on $\PSH(M,\o_0)\cap C^0(M)$ \cite[Proposition 2.3]{BT},
these functionals can be uniquely extended to pairs $(\o_{\vp_1},\o_{\vp_2})$,
with $\vp_1,\vp_2\in\calH^e_{\o}$ with $\o$ the
reference edge metric as in \eqref{oEq}.
Thus, these functionals are well-defined on
$\calH^e_{\o}\times\calH^e_\o$. Moreover, \eqref{IkSecondEq}
is still justified for edge metrics: working on a tubular
neighborhood of $D$ the boundary term one obtains from
integrating by part in \eqref{IkEq} tends to zero with the radius of the neighborhood about
$D$.

As we just saw, the functionals $I_k$ are a natural generalization of the Dirichlet
energy on the one hand, and Aubin's $J$ functional on the other.
It turns out that $I_k$ play a role in describing
a scale of energy functionals $E_k$ that similarly
depend gradually on additional `mixed terms' as $k$
increases, loosely speaking.
This can be made precise using the language of Bott--Chern
forms that we review
in \S\ref{BottChernSubSec}. Another, more down to earth way, of defining
this functional is as follows.

Denote the normalized elementary symmetric polynomials
of the eigenvalues of the twisted Ricci form
$\Ric\,\eta-(1-\be)[D]$ (with respect to $\eta$)
by
$$
\sigma^\be_k(\eta)=
\frac{(\Ric\,\eta-(1-\be)[D])^k\w\eta^{n-k}}{\eta^n},\quad k=1,\ldots,n,
$$
and their average (which is independent of the representative of $[\eta]$)
over $(M,\on)$ by
\beq
\label{mukEq}
\mu^\be_k:=
\frac{(c_1-(1-\be)c_1(L_D))^{k}\cup[\eta]^{n-k}([M])}{[\eta]^n([M])}.
\eeq
When $k=1$, $\sigma_1^\be(\eta)=s_\eta/n=\tr_\eta(\Ric\,\eta-(1-\be)[D])/n$,
where $s_\be(\eta)$ denotes the twisted scalar curvature of $\eta$.
When $k=n$ this is $\det_\eta(\Ric\,\eta-(1-\be)[D])$,
the determinant of the twisted Ricci curvature with respect to
the metric.

A straightforward generalization
of a theorem of the author from the case $\be=1$
\cite[Proposition 2.6]{RJFA} yields the following.

\begin{prop}
\label{EkProp}
Let $k\in\{0,\ldots,n\}$.
Suppose $\mu[\eta]=c_1-(1-\be)c_1(L_D)$.
The 1-form
$$
\eta\mapsto
\Big[\D_\eta\sigma^\be_k(\eta)-\frac{n-k}{k+1}\Big(\sigma^\be_{k+1}(\eta)-\mu^\be_{k+1}\Big)\Big]\etan
$$
is exact. Its potential, considered
as a function on $\calH_\eta^e\times\calH_\eta^e$,
is given by
\beq
\label{SecondBMRPropEq}
\mu^{-k}E^\be_k(\eta,\etavp)
 =
E^\be_0(\eta,\etavp)+\mu I_k(\etavp,
\mu^{-1}\Ric\,\!\etavp-\mu^{-1}(1-\be)[D])
-\mu I_k(\eta,
\mu^{-1}\Ric\,\!\eta-\mu^{-1}(1-\be)[D]),
\eeq
uniquely determined by the normalization
$E_k^\be(\eta,\eta)=0$.
\end{prop}

\bpf
As remarked earlier, $I_k$ is continuous in the topology
of uniform convergence (on the level of \K potentials)
by a result of Bedford--Taylor.
By definition (see
the second part of Definition \ref{KEEDef},
\eqref{KEdgePotentialsEtaSpace}, and
\eqref{gEdgeDef}), $\Ric\ovp-(1-\be)[D]$ admits
a continuous \K potential with respect to $\o$, and moreover
this potential can be suitably uniformly approximated.
Smoothly approximate
$(\o,\ovp)$ uniformly while approximating
$\Ric\ovp-(1-\be)[D]$ uniformly on
the level of potentials. Thus, it suffices to verify
the proposition when the \K forms are smooth, when integration
by parts, precisely those in the proof in the smooth case \cite{RJFA}, are justified. But then the proof in the smooth case can be applied
verbatim, with the definition of the twisted Ricci potential as in
\eqref{foEq} (recall $f_\o$ is continuous for all $\be$)
making the proof consistent with the introduction
of the terms $(1-\be)[D]$ in the right hand side
of \eqref{SecondBMRPropEq}.
\epf

Similar arguments allow to generalize the following
formula of Tian \cite{T1994} from the smooth case.

\begin{lem}
Let $\eta,\eta_\vp\in \calH^e_\o$. One has,
\begin{equation}
\label{EbetaFormula}
E^\be_0(\eta,\eta_{\vp})=
\frac1V\int_M\log\frac{\eta_{\vp}^n}{\etan}\eta_{\vp}^n
-\mu(I-J)(\eta,\eta_{\vp})
+\frac1V\int_M f_\eta(\etan-\eta_{\vp}^n).
\end{equation}
\end{lem}

%%%%%%%%%%%%%%%%%%%%%%%%%%%%%%%

The functionals $E_0$ and $E_n$ were introduced by Mabuchi \cite{Mabuchi1986}
and Bando--Mabuchi \cite{BandoMabuchi}, while the remaining ones by Chen--Tian \cite{ChenTian2002}.
Formula \eqref{SecondBMRPropEq} shows that the $E_k$ interpolate
between $E_0$ and $E_n$, to wit \cite[(17)]{RJFA}
$$
\begin{aligned}
E^\be_0(\eta,\eta_{\vp}) &=
((1-\hbox{$\frac{l}{k+1}$})E_0+\hbox{$\frac{l}{k+1}$}E_n)(\eta,\etavp)+(I_k-\hbox{$\frac{l}{k+1}$}J)(\etavp, \mu^{-1}\Ric\,\!\etavp-(1-\be)[D]/\mu)\qq&
\cr\cr
& \qq-(I_k-\hbox{$\frac{l}{k+1}$}J)(\eta,
\Ric\,\eta),\q\forall\,l\in\{0,\ldots,k+1\}.
\end{aligned}
$$
Since $E_0$ is known as the K-energy or \K energy, and $E_n$
as the Ricci energy, it is natural, following \cite{RJFA}, to refer to $E_k$
as the {\it K\"ahler--Ricci energies.}

The Ricci energy is special. It lends a geometric interpretation
to the Ding functional, via the inverse Ricci operator,
in the sense of \cite[\S9]{R08}. That is, suppose that
\eqref{cohomEq} holds with $\mu=1$.
Define by $\Ric^{-1}_\be:\calH_\o^\eps\ra\calH_\o^\eps$
the twisted inverse Ricci operator, by letting
$\Ric^{-1}_\be\eta:=\chi$ be the unique \K form
cohomologous to $\eta$ satisfying
$\Ric\chi-(1-\be)[D]=\eta$. Such a unique form exists
by Conjecture \ref{TianConj}, since this equation
can be written as a \MA equation of identical form to that
of the equation for a Ricci flat \K edge metric.

\begin{prop}
\label{EnFProp}
One has $(\Ric_\be^{-1})^\star E^\be_n=F^\be.$
\end{prop}

Again, the proof is a simple adaptation of the proof
in the smooth case \cite[Proposition 10.4]{R08}, following
the arguments in the proof of Proposition \ref{EkProp}
above.

\subsection{Bott--Chern forms}
\label{BottChernSubSec}

In this subsection we will describe the theory of  Bott--Chern forms and
energy functionals, inspired by work of Bott and Chern \cite{BottChern1}
and developed by Donaldson and Bismut--Gillet--Soul\'e
\cite{Donaldson1985,BismutGS}.
This will be applied to showing that the functionals $E_k^\be$
have an expression in terms of Bott--Chern forms, slightly
generalizing but very closely following the discussion in
\cite[\S4.4.5]{RThesis}.

The main idea of Bott--Chern forms is that given a moduli space of Hermitian
metrics on a bundle one may construct canonically defined ``universal"
functions on it associated to curvature. These functions arise via
a ``potential" for the curvature form of a ``universal" bundle over the whole
moduli space.

Let $E\ra M$ be a holomorphic vector bundle of rank $r$.%
\footnote{For the general discussion of Bott--Chern forms we will only make use
of the fact that
$
 M$ is complex (rather than K\"ahler).
In our applications we will always work with line bundles, i.e.,
$
 r=1$.
If one were interested only in that case the discussion below could be
slightly simplified. However we have chosen to maintain this level of generality.}
A vector bundle represents a \Cech cohomology class in $H^1(M,\calO_M(GL(r,\CC)))$.
Here by $\calO_M(GL(r,\CC))$ we mean the sheaf of germs of holomorphic functions
to $GL(r,\CC)$. When $r=1$ it is denoted by $\calO_M^\star$.
Let us identify $E$ with its \Cech class representative, i.e.,
by a collection of transition functions $g=\{g_{\a\be}\}$ that are
holomorphic maps from the intersection of any two coordinate
neighborhoods $U_\a, U_\be\subset M$ to $GL(r,\CC)$,
$g_{\a\be}:U_\a\cap U_\be\ra GL(r,\CC)$,  satisfying the \Cech
cocycle conditions \cite[p. 66]{GH}
$$
(\delta g)_{\a\be\gamma}:=g_{\a\be}\cdot g_{\be\gamma}\cdot g_{\gamma\a}=\III.
$$
Note that here the groups comprising the sheaf have a multiplicative structure
(and not an additive one) and hence $\delta g=\III$ expresses closedness, i.e.,
$[g]$ represents a \Cech cohomology class.

Denote by $\H_E$ the space of all Hermitian metrics on $E$. Let
$\Herm(r)$ denote the space of positive Hermitian $r\times r$ matrices.
Any Hermitian metric $H\in\H_E$ can be represented by smooth maps
$H_{\a}:U_\a\ra$ Herm$(r)$ such that with respect to
local bases of sections one has
$(H_\a)_{i\b j}=(g_{\a\be})_{ik}(H_\be)_{k\b l} \overline{(g_{\a\be})_{jl}}$,
or simply $H_\a=g_{\a\be}^\star H_\be g_{\a\be}$.
This is summarized in the notation $H=\{H_\a\}$.

To every $H\in \H_E$ there is associated a unique complex connection $D_H$.
For a general holomorphic vector bundle the connection $D_H$ is a 1-form on $M$ with
values in the bundle $\End(E)$ of endomorphisms of $E$. One may write globally
$$
D_H=H^{-1}\circ \del\circ H.
$$
The derivation of this formula follows  the same argument as in the line bundle case.
The exact meaning of how to understand this and other similar expressions involving
compositions of endomorphisms and differential operators
will be explained in detail below in several computations.
With respect to a local holomorphic frame $e_1,\ldots,e_r$ over $U_\a\subset M$
the endomorphism $H$
may be represented by a matrix and one has the special expression
$(D_H|_{U_\a})_j^i=\del h_{j\b l}\cdot h^{\b li}$. This expression is not valid
with respect to an arbitrary frame.

The global expression for the curvature
is then
\beq\label{CurvatureHermitianMetricBundleEq}
F_H:=D_H\circ D_H=\dbar\circ H^{-1}\circ\del\circ H.
\eeq
Implicit in this notation as well as in the sequel is the convention of working
with endomorphisms with values in the exterior algebra of differential forms on $M$.
Locally with respect to a holomorphic frame one has the special expression
$F_H|_{U_\a}=\dbar(\del h_{j\b l}\cdot h^{\b li})$. This expression is not
valid with respect to an arbitrary frame. However it demonstrates that $F_H$ is
$(1,1)$-form on $M$ (with values in $\End(E)$), since the type is independent
of the choice of frame.

Let $\phi$ denote an elementary symmetric polynomial on
$gl(r,\CC)\times\cdots\times gl(r,\CC)$ ($p$ times)
that is invariant under the adjoint action of $GL(r,\CC)$ (conjugation).
The idea behind Chern classes is that one may plug into such polynomials
matrices that have differential forms as their entries, for example $F_H$.
Since the polynomials
are $GL(r,\CC)$-invariant one obtains a differential form that is
invariant under a change of local trivializations for the bundle, hence is
intrinsically defined. Moreover it turns out that such forms are closed,
hence define intrinsic cohomology classes that depend only on the complex
structure of $E\ra M$. 

Now we come back to our original task of constructing functions
on $\calH_E$.
One may show that $\phi(F_H):=\phi(F_H,\ldots,F_H)$ is a closed
$2p$-form. It certainly depends on the metric $H$, however its cohomology
class in $H^{2p}(M,\ZZ)$
does not. This means that the difference
$$
\phi(F_{H_0})-\phi(F_{H_1})
$$
is exact. Moreover, and here we arrive at the main point of Bott--Chern theory,
one may find a $(p-1,p-1)$-form $\BC(\phi;H_0,H_1)$, well-defined up to $\del$- and $\dbar$-exact
forms, such that
$$
\dbar\del \BC(\phi;H_0,H_1)=\phi(F_{H_1})-\phi(F_{H_0}).
$$
The form $\BC(\phi;H_0,H_1)$ may then be integrated against a $(n-p+1,n-p+1)$-form on $M$
to give a number. Fixing $H_0$ and letting $H_1$ vary we therefore obtain
a function on $\H_E$ as desired. And, if $p-1=n$, we even do not need to make
a choice of such a form in order to integrate (this will be the case in our applications).
We now show how to construct the Bott--Chern form $\BC(\phi;H_0,H_1)$.
The proof below is a slow pitch version of the orignal one.

\begin{prop}
\label{BottChernProp}
{\rm (See \cite[Proposition 6]{Donaldson1985}.)}
Let $\phi$ be a $GL(r,\CC)$-invariant elementary symmetric polynomial.
Given $H_0,H_1\in \H_E$ and any path $\{H_t\}_{t\in[0,1]}$ in $\H_E$ connecting them,
the $(p-1,p-1)$-form
\beq\label{DifferenceCurvatureFormsEq}
\BC(\phi;H_0,H_1)
:=
p(\i)^{p-1}\int_{[0,1]}
\phi(H_t^{-1}\dot H_t,F_{H_t},\ldots,F_{H_t}) dt \q\h{\rm mod}\q \Im\del +\Im\dbar,
\eeq
is well-defined, namely does not depend on the choice of path.
In addition,
\beq\label{BottChernCocycleConditionEq}	
\BC(\phi;H_0,H_1)+\BC(\phi;H_1,H_2)+\BC(\phi;H_2,H_0)=0,
\eeq
and
\beq\label{ddbarUniversalCurvatureEq}
\dbar\del \BC(\phi;H_0,H_1)=\phi(F_{H_1})-\phi(F_{H_0}).
\eeq
\end{prop}

Notice that the first argument of $\phi$ is an endomorphism while the rest of its
arguments are endomorphism-valued $2$-forms.

\bigskip
\bpf
Note that $\BC(\phi;H_0,H_1)$ is given by integration over a path connecting
$H_0$ and $H_1$ of a globally defined 1-form on $\H_E$ with values in $(p-1,p-1)$-forms on $M$.
We call this form $\theta$. To show independence
of path we show that this form is closed modulo $\del$- and $\dbar$-exact terms.
Let $H\in \H_E$. Let $h,k\in T_H\H_E$ and extend them to constant vector fields near $H$.
Then,
\beq
\label{ClosedOneFormEndomorphismEq}
d\theta(h,k)
=
k\,\theta|_H(h)-h\,\theta|_H(k).
\eeq

First we obtain an expression for the infinitesimal change of the curvature under
a variation of a Hermitian metric.
Write $H+tk=H\circ(\III+t H^{-1}k)=:H\circ f$. We write $\circ$ to emphasize that
composition of endomorphisms is taking place (in coordinates multiplication of matrices).
According to \eqref{CurvatureHermitianMetricBundleEq} we have
\beq
\begin{aligned}
F_{H\circ f}
&
=
\dbar\big( (H\circ f)^{-1}\circ\del (H\circ f) \big)
\cr
& =
\dbar
\big( f^{-1}\circ(H^{-1}\circ\del H)\circ f+f^{-1}\circ\del f \big)
\cr
& = \dbar\circ f^{-1}\circ D_{H}^{1,0}\circ f.&
\label{CircCurvatureEq}
\end{aligned}
\end{equation}
Here it should be emphasized that $D_H$ decomposes according to type (of its 1-form part) into $D_H^{1,0}$ and
$D_H^{0,1}$ and that while originally the connection $D_H$ was defined on $E$ it may be extended
naturally to $\End(E)$ and it is this extension that we use in the equation above ($f$ is a section
of $\End(E)$ and not of $E$). The same applies to the operator $\dbar$ that we also extend
to act on $\End(E)$.
To understand the last equation better we let the endomorphism it defines act on a
holomorphic section $s$
of $E$, and compute:
$$
\begin{aligned}
F_{H\circ f}s
& =
\dbar\circ f^{-1}\circ D_{H}^{1,0} (fs)
\cr
& =
\dbar\circ f^{-1}\circ \big( (D_{H}^{1,0}f)s + f D_{H}^{1,0}s \big)
\cr
& =
\Big(\dbar\big(f^{-1}(D_{H}^{1,0}f)\big) + \dbar\circ D_{H}^{1,0}\Big) s.
\end{aligned}
$$
Therefore we write
$$
F_{H\circ f}=F_H+\dbar\big(f^{-1}(D_{H}^{1,0}f)\big),
$$
and the second term should be understood to be distinct from \eqref{CircCurvatureEq}.
This subtle notational issue can be a cause for great confusion when consulting
the literature on vector bundles and Yang-Mills theory.
Putting now $f=\III+tH^{-1}k$ we obtain
$$
F_{H+tk}=F_H+t\dbar\big(D_{H}^{1,0}(H^{-1}k)\big)+O(t^2).
$$

Hence the first term in \eqref{ClosedOneFormEndomorphismEq} is given by
\beq
\begin{aligned}
\frac1{p(\i)^{p-1}}
k\theta|_H(h)
& =
\frac d{dt}\Big|_0
\phi((H+tk)^{-1}h,F_{H+tk},\ldots,F_{H+tk})
\cr
& =
\phi(-\i Hk\i Hh,F_{H},\ldots,F_{H})
\cr
&\q +
{\text\sum}\phi(\i Hh,F_{H},\ldots,\dbar D_{H}^{1,0}(\i Hk),\ldots,F_H).
\end{aligned}
\eeq
Therefore one has, putting $\sigma=\i Hh,\; \tau=\i Hk$,
\beq
\baeq
\frac1{p(\i)^{p-1}}d\theta(h,k)
& =
\phi([\s,\tau],F_{H},\ldots,F_{H})
+{\text\sum}\phi(\s,F_{H},\ldots,\dbar D_{H}^{1,0}\tau,\ldots,F_H)
\cr
&\q -{\text\sum}\phi(\tau,F_{H},\ldots,\dbar D_{H}^{1,0}\s,\ldots,F_H).
&\label{dthetaeq}
%\dda
%
\eaeq
\eeq
$\s$ and $\tau$ are sections of the endomorphism bundle
of $E$. Note that as operators on this bundle one has
$$
\dbar \circ D_{H}^{1,0}+D_{H}^{1,0}\circ\dbar
=
(\dbar+D_{H}^{1,0})\circ(\dbar+D_{H}^{1,0})=D_H^2=F_H,
$$
since we already saw that the curvature is of type $(1,1)$.
Hence for example
\beq
\dbar\circ D_{H}^{1,0}\sigma=-D_{H}^{1,0}\circ\dbar\sigma+[F_H,\sigma].
\label{CommutatorCurvatureBundleEq}
\eeq
Note $[F_H,\sigma]\equiv F_H\s$, and the bracket notation simply emphasizes
that we have extended $F_H$ to act on the endomorphism bundle and so the endomorphism
part of $F_H$ will actually act by bracket on the endomorphism $\sigma$ (the 2-form part
will simply be multiplied along).
By the Bianchi identity $D_HF_H=0$ and
so $D_{H}^{1,0}F_H=0, \; \dbar F_H=0$.
Now,
\beq
\baeq
\phi(\s,\dbar D_{H}^{1,0}\tau,F_H,\ldots,F_H)
& =
\dbar\phi(\s,D_{H}^{1,0}\tau,F_H,\ldots,F_H)
\cr
& \q-\phi(\dbar\s,D_{H}^{1,0}\tau,F_H,\ldots,F_H)
\cr
&\q -{\text\sum}\phi(\s,D_{H}^{1,0}\tau,F_H,\ldots,\dbar F_H,\ldots,F_H)
\cr
& =\dbar\phi(\s,D_{H}^{1,0}\tau,F_H,\ldots,F_H)
\cr
&\q-\phi(\dbar\s,D_{H}^{1,0}\tau,F_H,\ldots,F_H).
&
\label{Firsteq}
\eaeq
\eeq
Using \eqref{CommutatorCurvatureBundleEq} the corresponding term in the second sum of \eqref{dthetaeq} is
\beq
\baeq
-\phi(\tau,\dbar D_{H}^{1,0}\s,F_H,\ldots,F_H)
& =
-\phi(\tau,-D_{H}^{1,0}\dbar\s+[F_H,\s],F_H,\ldots,F_H)
\cr
& =
-\phi(\tau,[F_H,\s],F_H,\ldots,F_H)
+\del\phi(\tau,\dbar\s,F_H,\ldots,F_H)
\cr
&\q
-\phi(D_{H}^{1,0}\tau,\dbar\s,F_H,\ldots,F_H)
\cr
&\q
-{\text\sum}\phi(\tau,\dbar\s,F_H,\ldots,D_{H}^{1,0}F_H,\ldots,F_H).
\cr
& =
\phi(\tau,[\s,F_H],F_H,\ldots,F_H)
+\del\phi(\tau,\dbar\s,F_H,\ldots,F_H)
\cr
&\q
-\phi(D_{H}^{1,0}\tau,\dbar\s,F_H,\ldots,F_H).
&\label{Secondeq}
\eaeq
\eeq
Note that it is not necessarily true that
\beq
-\phi(D_{H}^{1,0}\tau,\dbar\s,F_H,\ldots,F_H)\label{termoneeq}
\eeq
cancels with
$$
-\phi(\dbar\s,D_{H}^{1,0}\tau,F_H,\ldots,F_H).
$$
But, eventually taking the equations
\eqref{Firsteq} and
\eqref{Secondeq} for all pairs appearing in the sums \eqref{dthetaeq} then, e.g.,
the term \eqref{termoneeq} will cancel with the term
$$
-\phi(\dbar\s,F_H,\ldots,F_H,D_{H}^{1,0}\tau).
$$
Indeed we are only allowed to permute the arguments of $\phi$ cyclically (e.g., for three matrices
$A,B,C$ one has $\tr(ABC)=\tr(CAB)$, but in general $\tr(ABC)$ is different from $\tr(BAC)$).
Note also that while $\phi$ does not change when permuting matrices cyclically, when
we permute cyclically matrix valued 1-forms a sign appears, as usual. This explains the cancellation above.

Hence, modulo $\del$- and $\dbar$-exact terms, we are left with
$$
\sum\phi(\tau,F_H,\ldots,[\s,F_H],F_H,\ldots,F_H)
$$
which cancels with the first term in \eqref{dthetaeq}; this can be seen by using the
invariance of $\phi$ under the action of $GL(r,\CC)$ by conjugation
$$
\frac d{dt}\Big|_0\phi(e^{-tB}A_1e^{tB},\ldots,e^{-tB}A_pe^{tB})
=\sum\phi(A_1,\ldots,[A_j,B],\ldots,A_p),
$$
concluding the proof.
\epf

\begin{example}
{\rm
\label{FzeroBottChernExam}
Let $E$ denote an ample line bundle polarizing a \K class $\O=c_1(E)$.
We identify $\H_E$ with $\H_\o$ where $\o=-\i\ddbar\log h=\i\dbar(h^{-1}\del h)=\i F_h$
is a \K form with $h\in\H_E$
(and hence $[\o]=\O$).
Now $r=1$ and so no traces are needed (the matrices are all one-dimensional).
Put
\beq\label{PhiInvariantGLFunctionEq}
\phi(A_1,\ldots,A_{n+1}):=A_1\cdots A_{n+1}.
\eeq
Take a path of Hermitian metrics $h_t=e^{-\vp_t}h$.
Then $F_{h_t}=F_h+\ddbar\vp_t=\o_{\vp_t}/\i$. Then the
Bott--Chern form is
\beq\label{FzeroBottChernEq}
\BC(\phi;h_0,h_1)
=(n+1)\int_0^1(\i)^n\phi(-\dot\vp_t,F_{h_t},\ldots,F_{h_t}) dt
=-(n+1)\int_0^1\dot\vp_t\,\ovptn\w dt.
\eeq
This expresses a $2n$-form on $M$. Integrating this over $M$ gives
the function $-(n+1)L(\vp)$,
on $\H_\o$.
By \ref{BottChernProp} this is independent of the choice of path since any
two choices differ by $\del$- and $\dbar$-exact terms and hence by $d$-exact
terms since our expressions are real.
}
\end{example}

%%%%%%%%%%%%%%%%%%%%%%%%%

Now we turn to the setting of a \K manifold with an integral \K class $[\o]$. Let $L$ be
a line bundle polarizing the \K class, namely $c_1(L)=[\o]$.
Let $K_M^{-1}$ denote the anticanonical bundle polarizing the class $c_1$,
and $L_D$ the line bundle associated to $D$.
Given a Hermitian metric $h$ on $L$ of positive curvature,
and a global holomorphic section $s$ of $L_D$
one obtains a metric $\det F_h |s|^{2-2\be}=\det(\o/\i)$
on $K_M^{-1}\otimes L_D^{\be-1}$
(we mean that locally $\det F_h=\det (g_{i\b j})$ if
$\o=\i g_{i\b j}dz^i\w \overline{dz^j}$)
where $\o=-\i\ddbar\log h$. We write $h\in\Ho=:\H_L^+$. Note that
the Bott--Chern forms defined below are defined on $\H_L^+$ rather
than on all of $\H_L$ (or to be more precise on a set isomorphic to $\H_L^+$, see
\cite[p. 214]{Tian2000}).

The main result of this subsection is the following expression for the K\"ahler--Ricci
functionals in terms of Bott--Chern forms. While this generalizes a theorem of Tian
for the K-energy (the case $k=0$) \cite[\S2]{Tian2000}
and its extension by the author
to all $k$ \cite[\S4.4.5]{RThesis},
the computations are almost identical.
In essence, the (twisted) K\"ahler--Ricci functionals are realized as a linear combination of Bott--Chern forms, one
for each of the $\RR$-line bundles $E_j=K_M^{-1}\otimes L_D^{\be-1}\otimes L^{n-2j}$ for $j=0,\ldots,n$. One
of these terms (the simplest contribution) is
a multiple of the form appearing in Example \ref{FzeroBottChernExam}.

\begin{thm}
Let $k\in\{0,\ldots,n\}$ and let $\phi$ be defined by \eqref{PhiInvariantGLFunctionEq}.
Let $(M,J,\o)$ be a projective \K manifold, and let $L$ be a line bundle
with $c_1(L)=[\o]$. Let $\mu^\be_k$ be given by \eqref{mukEq}.
For each $k\in\{0,\ldots,n\}$,

\beq
\begin{aligned}
\label{dEkEq}
\Big[&\D_\o\sigma^\be_k-\frac{n-k}{k+1}\Big(\sigma^\be_{k+1}-\mu^\be_{k+1}\Big)\Big]\o^n
=\cr
&
\frac{2^{-n}}{(n+1)!}{n\choose k}^{-1}
\frac1V
\sum_{j=0}^n
(-1)^j{n\choose j} (n-2j)^k
\,
\BC(\phi;h_0^{n-2j} \det F_{h_0}|s|^{2-2\be},h_1^{n-2j}\det F_{h_1}|s|^{2-2\be})
\cr
&
\q\;\; -\frac1V\frac{\mu^\be_{k+1}}{n+1}\frac{n-k}{k+1}\BC(\phi;h_0,h_1).&
\end{aligned}
\eeq
\end{thm}

\bpf
The proof in the case $\be=1$ is given in \cite[Proposition 4.22]{RThesis}.
The general case follows from the same computations by subtracting
the {\it fixed} current $(1-\be)[D]$ from each Ricci current that appears in the computation.
\epf

Remark that from the proof it follows that similar ``twisted'' functionals can
thus be defined by replacing $(1-\be)[D]$ with some other fixed
curvature current of a line bundle.

\begin{remark}{\rm
Tian \cite[p. 255--257]{T1994}
gave an interpretation of the
``complex Hessian" of the K-energy $E_0$ in terms of a
certain ``universal" Hermitian metric $\bf h$ (see the references above for the notation and definitions):
\beq
-\i\ddbar E_0=\frac1V\int_M (\i)^{n+1}
\Big( F_{\det F_{\bf h}}-\frac{n\mu_0}{n+1}F_{\bf h} \Big)\w (F_{\bf h})^n.\label{TianFormulaEq}
\eeq
We note in passing the following generalization of this formula to the K\"ahler--Ricci functionals
\beq
-\i\ddbar E_k=\frac1V\int_M (\i)^{n+1}
\Big( (F_{\det F_{\bf h}})^{k+1}
-
\frac{(n-k)\mu_{k}}{(k+1)(n+1)}(F_{\bf h})^{k+1} \Big)\w (F_{\bf h})^{n-k}.\label{TianFormulaEq}
\eeq
We also remark that using the techniques of \cite{RJFA} one may generalize
appropriately Theorem \ref{PropernessThm} in terms of properness of
the functional $E_k^\be$ on the space of \K forms $\eta$ for which
$\Ric\eta-(1-\be)[D]$ is positive and cohomologous to $[\eta]$.
Finally, also $L$ has such a formula due to Tian \cite[p. 214]{Tian2000},
$$
\i\ddbar L=\frac1V\int_M(\i F_{\bf h})^{n+1},
$$
that has been extended by Berman--Boucksom \cite[(4.1)]{BermanBoucksom}
to \K potentials with low regularity, and can be used to characterize
weak geodesics in the Mabuchi metric.
}
\end{remark}

\subsection{Legendre transform} 
\label{LegendreSubsubsec}

In this section we restrict to the case $\mu>0$ to simplify
the notation, and describe work of Berman \cite{Berman2010}
and Berman--Boucksom \cite{BermanBoucksom}
that ties $F^\be$ with $E^\be_0$ via
the Legendre transform.

Defining the probability measure (recall the normalization
\eqref{foEq})
\beq
\nu_\eta:=(\Ric^{-1}\eta)^n=
{e^{f_\eta}\etan},
\eeq
we rewrite \eqref{FbetaEq} and \eqref{EbetaFormula} as
\begin{equation}
\label{FbetaEbetaFormula}
\begin{aligned}
F^\be(\eta,\etavp)
&=
-L_\eta(\vp)-\frac1\mu\log\frac1V\int e^{-\mu\vp}\nu_\eta,
\cr
E^\be_0(\eta,\eta_{\vp})
&=
\frac1V\int_M\log\frac{\eta_{\vp}^n}{\nu_\eta}\eta_{\vp}^n
-\mu(I-J)(\eta,\eta_{\vp})+\frac1V\int_M f_\eta\etan.
\end{aligned}
\end{equation}
The last term in $E^\be_0$ is a constant. It can be
eliminated if we had normalized $\int_M f_\eta\etan=0$,
and then set $\nu_\eta:=\frac{e^{f_\eta}\etan}
{\frac1V\int e^{f_\eta}\eta^n}$.
The first term in $E^\be_0$ is the entropy of $\etavp^n$ relative to the measure $\nu_\eta$
considered as a functional on the space of measures,
\beq
\Ent(\nu,\chi)=
\frac1V\int_M\log\frac{\chi}{\nu}\chi,
\eeq
which in terms of the density $d=\chi/\nu$
takes the familiar form for the entropy $\int d\log d\, \nu$.
On the other hand, it is classical that the last term in $F^\be$ is precisely the Legendre transform of the entropy, in the sense
that
\cite[p. 264]{DemboZ},
\beq
\label{EntropyDualityEq}
\begin{aligned}
\Lambda_\nu(-\mu\vp)=
\log\frac1V\int e^{-\mu\vp}\nu
=\Ent(\nu,\,\cdot\,)^\star(-\mu\vp)=\sup_{\chi\in\calV_V}\{\langle -\mu\vp,\chi\rangle - \Ent(\nu,\chi)\},
\end{aligned}
\eeq
where $\calV_V=\{\nu : 0\le\nu / \o^n \in C^0(M, \o^n), \int\nu=V\}$.
Conversely, by convexity,
$$
\Ent(\nu,\mu)=\sup_{\psi\in C^0(M)}\{\langle\mu,\psi\rangle
- \log\frac1V\int e^{\psi}\nu\}=\Lambda_\nu^\star(\mu).
$$
One of Berman's insights was that the remaining terms in $E^\be_0$ and
$F^\be$ are similarly related by the Legendre transform
\cite[\S2]{Berman2010} (cf. \cite[Theorem 5.3]{BBGZ}).
\begin{lem}
\label{BermanLemma}
Let $\mu>0$. Then,
$\sup_{\psi\in\PSH(M,\o)\cap C^0}\{ \langle -\mu\psi,\ovpn\rangle +\mu L_\o(\psi) \}
=\mu(I-J)(\o,\o_{\vp}). $

\end{lem}

We could have introduced a minus sign into the usual
inner product between functions and measures in order to
obtain that $(-L)^\star = I-J$.
Instead of doing that, we kept the usual
inner product but then the left hand side in Lemma \ref{BermanLemma}
is not precisely the Legendre transform.

\begin{proof}
The proof would be easier if we knew that the supremum
over all functions coincides with that over $\o$-psh ones.
Indeed in that case,
$F(\psi):=\langle -\mu\psi,\ovpn\rangle +\mu L_\o(\psi)$
is concave in $\psi$ in the sense that
$$
-\frac{d^2}{dt^2}\Big|_{t=0}F(\psi+t\phi)=
\mu\int n\i\del\phi\w\dbar\phi\o_\psi^n\ge0
$$
(recall that $dL_\o|_\psi=\o_\psi^n$).
Differentiating,
one sees a critical point $\psi$, necessarily a maximum, must satisfy
$\o_\psi^n=\o_\vp^n$
(here we are being a bit loose since $C^0(M)$ is infinite-dimensional;
however the same reasoning as for the Legendre transform in
finite-dimensions applies), hence by uniqueness $\vp=\psi+C$
\cite{Bl2003}.
Plugging back in, and using the formula
$L(\vp)=(I-J)(\o,\ovp)+\frac1V\int\vp\ovpn$ then yields the statement.

To make the reduction to the subset
$\PSH(M,\o)\cap C^0(M)\subset C^0(M)$, recall
the definition of the $\o$-psh envelope operator
$P_\o:\vp\mapsto \sup\{\phi\in \PSH(M,\o)\cap C^0(M)\,:\, \phi\le\vp\}$.
By a result of Berman--Boucksom \cite{BermanBoucksom},
$L_\o\circ P_\o$ is concave on $C^0(M)$, Gateaux differentiable,
 and
$dL_\o\circ P_\o|_\vp=\o_{P_\o\vp}^n$.
Thus,
$\sup_{\psi\in C^0}\{ \langle -\psi,\ovpn\rangle + L_\o\circ P_\o(\psi) \}
=(I-J)(\o,\o_{P\vp}). $
This concludes the proof, since $P\vp\le\vp$
and so the supremum must actually be attained at
$\vp\in\PSH(M,\o)$, as $Lu\le Lv$ if $u\le v\le0$,
and we can normalize $\psi$ so that $\sup\psi=0$.
\end{proof}

\begin{remark}
{\rm
As pointed out by Berman, one can, in fact, avoid using the result of 
Berman--Boucksom
to prove the preceeding Lemma,
since simpler convexity arguments already show that the supremum
must be attained at $\varphi$. However, that result
is useful to show
that a maximizer $\psi$ must satisfy $\o_{\psi}^n=\o_{\varphi}^n$
(and hence $\psi\in L^\infty$
\cite{Kolodziej} and so equal to $\varphi$
up to a constant \cite{Bl2003}). In addition, 
the result of Berman--Boucksom is essentially needed
if one replaces $\o_\vp^n$ by  a more general measure or even
volume form $\nu$; the
result can be seen as the starting point of the variational approach to constructing
a weak solution of the equation $\o_\vp^n=\nu$.
}
\end{remark}

\subsection{Equivalence of functionals}
\label{EquivalenceSubSec}

Motivated by variational calculus, one expects
that the KEE problem is solvable if and only if it
can be cast as a variational problem with a coercive functional.
One calls $E$ coercive if
there exist uniform positive constants $A,B$ such that the $c$-sublevel
set of $E$ is contained in the $A(c+B)$-sublevel set
of $J$, i.e., $E\ge \frac1A J-B$. In particular, $E$ is then proper,
thus bounded from below. We have seen a number of functionals
that have KEE metrics as critical points. The following
basic result says that they are all more or less
equivalent as far as boundedness, coercivity, and existence
of KEE metrics is concerned. To state it we introduce some notation.
Suppose that $\mu>0$, and define
\begin{equation}
\label{}
\calH_\o^{e,+}=\{\vp\in\calH_\o^e\,:\, \Ric\,\ovp-(1-\be)[D]
\hbox{\ is a positive
current}\}.
\end{equation}
This space is nonempty as a corollary of the  existence theorem for the
case $\mu=0$.
Let
$
l(\o)=\inf_{\vp\in\calH_\o^e} F^\be(\o,\ovp)
$
and
$$
l_k(\o)=\begin{cases}
\inf_{\vp\in\calH_\o^e} E_k(\o,\ovp),& \hbox{for\  $k=0,1$},\cr
\inf_{\vp\in\calH_\o^{e,+}} E_k(\o,\ovp),& \hbox{for\   $k=2,\ldots,n$}.\cr
\end{cases}
$$

\begin{thm}
\label{EquivLowerBdThm}
 Let $\mu>0$.
(i)
The lower bounds of $E^\be_k$ and that of $F^\be$ are related by
\beq
\label{EkLowerBoundEq}
\mu l(\o)+\frac1V\int f_\o\o^n= l_0(\o)
= \mu^{-k}l_k(\o)-\mu I_k(\o,\mu^{-1}\Ric\o-\mu^{-1}(1-\be)[D]).
\eeq
In particular, $F^\be$, $E^\be_0$ and $E^\be_1$ are simultaneously
bounded or unbounded from below on $\calH^e_\o$.
\hfill\break
(ii)
The coercivity of $E^\be_0$ is equivalent to that of $F^\be$.
\end{thm}

\bpf
(i)
First, by Proposition \ref{EkProp}
and the fact that $I_k$ is nonnegative on
$\calH^e_\o\times\calH^e_\o$, it follows that
if $E^\be_0$ is bounded below on $\calH^e_\o$
then $E^\be_k,\, k=1\ldots n$ is bounded below on $\calH^{e,+}_\o$.
In particular this is true for $E_n^\be$,
but then, by Proposition \eqref{EnFProp}
and the edge version of the Calabi--Yau theorem (Conjecture \ref{TianConj},
that guarantees that $\Ric_\be^{-1}:\calH^e_\o\ra
\calH^{e,+}_\o$ is an isomorphism)
it follows that $F^\be$ is bounded
below on $\calH^e_\o\times\calH^e_\o$. But
by a formula of Ding--Tian \cite{DingTian1993}
\beq
\label{DingTianEq}
(E^\be_0-\mu F^\be)(\o,\ovp)=\frac1V\int f_\o\on-f_{\ovp}\ovpn \ge\int f_\o\on,
\eeq
where we used the normalization \eqref{foEq} and
Jensen's inequality
$\frac1V\int f_{\ovp}\ovpn\le\log\frac1V\int e^{f_{\ovp}}\ovpn=0$.
This proves (i), since the precise lower  bounds
\eqref{EkLowerBoundEq} can be deduced from the
proof and a theorem of Ding--Tian---see \cite[Remark 4.5]{RJFA}.

We now give a second derivation, due to Berman,
of a special case of (i), namely the equivalence
of the lower bounds of $F^\be$ and $E^\be_0$.
First,
$
F^\be(\o,\,\cdot\,)\ge -C,
$
is equivalent to
$-\mu L_\o\ge \log\frac1V\int e^{-\mu\vp}\nu-\mu C$.
The Legendre transform,
in the sense of the previous
subsection,
 is order-reversing. Thus, according to
\eqref{EntropyDualityEq} and Lemma \ref{BermanLemma},
$
\mu(I-J)(\o,\ovp)\le \Ent(\nu,\ovpn)+\mu C
$
i.e.,
$
E^\be_0(\o,\,\cdot\,)\ge \frac1V\int f_\o\on-\mu C.
$
This concludes this derivation since the Legendre transform is an involution.

(ii) Suppose that $E^\be_0$ is coercive. Then, by definition,
$E^\be_0-\eps(I-J)$ is bounded from below. It follows
from (i) that so is $F^\be$, but now with $\mu$ replaced
by $\mu+\eps$, in other words
$$
-L_\o(\vp)\ge \frac1{\mu+\eps}\log\frac1V\int e^{-(\mu+\eps)\vp}\nu-C,
$$
for all $\vp\in\PSH(M,\o)\cap C^0$.
Normalize $\vp$ so that $\int \vp\on=0$, and
substitute $\frac{\mu}{\mu+\eps}\vp\in\PSH(M,\o)$
in this inequality to obtain
$$
J(\o,\o_{\mu\vp/(\mu+\eps)})=-L_\o\Big(\frac{\mu}{\mu+\eps}\vp\Big)
\ge \frac1{\mu+\eps}\log\frac1V\int e^{-\mu\vp}\nu-C,
$$
so
$$
\frac1{\mu}\log\frac1V\int e^{-\mu\vp}\nu-C'
\le
\frac{\mu+\eps}\mu J(\o,\o_{\mu\vp/(\mu+\eps)})
\le
\Big(\frac\mu{\mu+\eps}\Big)^{1/n} J(\o,\o_{\vp})
=:-(1-\eps')L(\o,\o_{\vp}),
$$
where the last inequality is due to Ding \cite[Remark 2]{Ding1988}.
Thus, $F^\be\ge \eps'J-C'$. The converse follows from
\eqref{DingTianEq}.
\epf

\begin{remark}
{\rm
Part (i) and its proof above is due to \cite{RJFA}.
The special case of the equivalence of
$F^\be$ and $E^\be_0$ 
being bounded below was obtained independently by H. Li \cite{HLi} using
results of Perelman \cite{SesumTian} on the Ricci flow,
and a third proof was later given by
Berman \cite{Berman2010} using the Legendre
transform, as presented above. Part (ii) is a special
case of a result of Berman (which allows to replace $\nu$ by
a rather general probability measure),
which generalized a result of Tian and
its subsequent refinement by Phong et al. (in the smooth case)
\cite{Tian1997,Tianbook,PSSW}. For a comparison of properness
and coercivity in the KE setting we refer to \cite[\S2]{Tian2012}.
}
\end{remark}

\section{The Ricci continuity method}
\label{RCMSec}

This section describes a new input that goes into the proof of
Conjecture \ref{TianConj} and Theorem \ref{KEEExistenceThm} in \S\ref{PositiveSubSec}.
It is a new continuity method that is essentially the only one that can be used to prove
existence of KEE metrics, since the classical continuity methods that were previously used
to construct KE metrics break down when $\be$ belongs to the more challenging
regime $(1/2,1)$ where, among other things, the curvature of the reference geometry
is no longer bounded (Lemma \ref{LiRLemma}). We describe all of this in detail in
\S\ref{MeetsSubSec}--\S\ref{ConvergenceSubSec}. Subsection
\ref{RISubSec} is an interlude about the Ricci iteration that origingally motivated
the Ricci continuity method. Finally, \S\ref{OtherApprSubSec} describes other approaches
to existence.

\subsection{Ricci flow meets the continuity method}
\label{MeetsSubSec}

The Ricci continuity method was introduced in \cite{R08}
and was further developed and first
used systematically to construct KE(E) metrics in \cite{MR2012,JMR}.
The idea is to prove existence of a continuity path, or, in other words,
a one-parameter family of \MA equations and solutions thereto, in a canonical
geometric manner. To that end, we start with the Ricci flow
\begin{equation}
\label{RFEq}
\frac{\del \o(t)}{\del t}
=
-\Ric\, \o(t)+(1-\be)[D]+\mu\o(t),\quad \o(0)=\o\in\calH,
\end{equation}
Fix $\tau\in(0,\infty)$.
The (time $\tau$) Ricci iteration is the sequence $\{\o_{k\tau}\}_{k\in\NN}\subset\calH_\o$,
satisfying the equations
$$
\o_{k\tau} =\o_{(k-1)\tau}+\tau\mu\o_{k\tau}-\tau\Ric\o_{k\tau}+\tau(1-\be)[D],
\qquad \o_{0\tau}  =\o,
$$
for each $k\in\NN$ for which a solution exists in $\calH_\o$.
This is the backwards Euler discretization of the flow \eqref{RFEq} \cite{R08}. Equivalently,
let $\o_{k\tau}=\o_{\psi_{k\tau}}$, with
$\psi_{k\tau}=\sum_{l=1}^k\vp_{l\tau}$. Then,
\begin{equation}
\label{RIEq}
\o_{\psi_{k\tau}}^n=\on e^{f_\o-\mu\psi_{k\tau}+\frac1\tau\vp_{k\tau}}.
\end{equation}
We now change slightly our point of view by fixing $k=1$ but instead varying
$\tau$ in $(0,\infty)$.
This yields $\o_{\vp_{\tau}}^n=\on e^{f_\o+(\frac1\tau-\mu)\vp_{\tau}},$
and setting $s:=\mu-\frac1\tau$, we obtain
\begin{equation}
\label{RCMEq}
\o_\vp^n=\on e^{f_\o-s\vp}, \quad s\in(-\infty,\mu],
\end{equation}
where $\vp(-\infty)=0$, and $\o_{\vp(-\infty)}=\o$. We call this the
{\it Ricci continuity path (Ricci CP).}

Aside from the formal derivation that relates
\eqref{RCMEq} to the Ricci flow,
it turns out that the two equations
share several key analytic and geometric properties:
\hfill\break
(i) Short-time existence: the Ricci CP exists for all
$0<\tau<<1$, i.e., for all $s<<-1$.
\hfill\break
(ii) Monotonicity: $\dot E_0(\o(0),\ovps)\le0$ with equality
iff $\o(0)$ is KEE.

\noindent
Moreover, the Ricci CP inherits one additional property
that is not satisfied by the Ricci flow, and which provides
an important advantage:
\hfill\break
(iii) Ricci lower bound: Along the Ricci CP $\Ric\,\o_{\vp(s)}
>s\o_{\vp(s)}$. Moreover, this holds even
if the initial metric has unbounded Ricci curvature!
Note that properties (i) and (iii) were first noticed by
Wu and Tian--Yau, respectively,
in their study of non-compact KE metrics with negative Ricci curvature
\cite{Wu,TianYau1986}, while (ii) goes back
to \cite{BandoMabuchi}.

In the rest of this section we will explain how to obtain a unified
proof of existence of both smooth and edge KE metrics
(the proof of Conjecture \ref{TianConj} and Theorem \ref{KEEExistenceThm})
 using the Ricci CP. In \S\ref{OtherApprSubSec}
we also review other approaches to existence.
But first, we make a comparison to some other CPs
and explain where each of them would break down in the edge setting.

\subsection{Other continuity methods}
\label{OtherCMSubSec}

The continuity path \eqref{RCMEq} has several useful properties, some already
noted above, which are necessary for the proof of Conjecture \ref{TianConj}
and Theorem \ref{KEEExistenceThm}
when $\be > 1/2$.
In other words, one could use various CPs (including the Ricci CP) when $\be$ is
in the ``orbifold regime" $\be\in(0,1/2]$, but it seems that
the CPs we discuss below break down in the regime $\be\in(1/2,1)$.
To illustrate this,  we now describe these CPs and where they
fail.

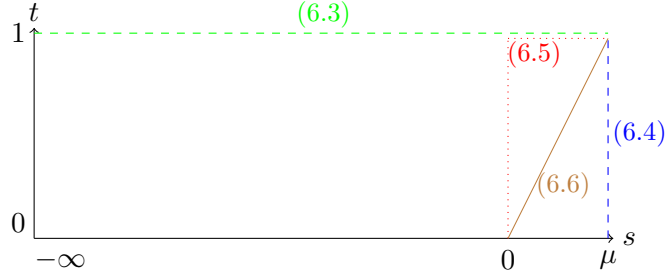
\begin{figure}
	\centering
\begin{subfigure}{.4\linewidth}
	\centering
\hglue-1cm
		\begin{tikzpicture}[scale=.7]
			\draw [->] (-9,0) -- (2,0);
			\draw [->] (-9,0) -- (-9,4);
			\draw [green, dashed] (-9,3.9) -- (1.9,3.9);
			\node [green] at (-3.5,4.25) {\small \eqref{RCMEq}};
			\draw [blue, dashed] (1.9,0) -- (1.9,3.9);
			\node [blue] at (2.5,2) {\small \eqref{CalabiCPEq}};
			\draw [red, dotted] (0,0) -- (0,3.8);
			\draw [red, dotted] (0,3.8) -- (1.9,3.8);
			\node [red] at (.5,3.5) {\small \eqref{AubinCPEq}};
			\draw [brown] (0,0) -- (1.9,3.8);
			\node [brown] at (1.05,1) {\small \eqref{DKCPEq}};
			\node  at (-8.5,-.4) {$-\infty$};
			\node  at (2.3,0) {$s$};
			\node  at (0,-.4) {$0$};
			\node  at (1.9,-.4) {$\mu$};
			\node  at (-9,4.3) {$t$};
			\node  at (-9.3,0.3) {$0$};
			\node  at (-9.3,3.9) {$1$};
		\end{tikzpicture}
\end{subfigure}
\caption{Continuity paths of the form \eqref{TwoParamCMEq} with $(s,t)\in A$ (assuming $\mu>0$).
The Ricci continuity path is \eqref{RCMEq}.}
		\label{FigCP}
\end{figure}

When $\mu\le0$, Calabi suggested the following path \cite[(11)]{Calabi1957}
that was later used by Aubin and Yau \cite{Aubin1978,Yau1978}
\begin{equation}
\label{CalabiCPEq}
\o_{\vp_t}^n=\on e^{tf_\o+c_t-\mu\vp_t},\q t\in[0,1].
\end{equation}
In the case $\mu>0$, Aubin suggested the following extension of Calabi's
path \cite{Aubin1984}
\begin{equation}
\label{AubinCPEq}
\o_{\vp_t}^n
=
\begin{cases}
\on e^{tf_\o+c_t},\q\qq\!\! t\in[0,1],\cr
\on e^{f_\o-(t-1)\vp_t},\q t\in[1,1+\mu].
\end{cases}
\end{equation}
Still when $\mu>0$, an alternative path was considered by Demailly and Koll\'ar
\cite[(6.2.3)]{DemaillyKollar}, given by
\begin{equation}
\label{DKCPEq}
\o_{\vp_t}^n
=
\on e^{t(f_\o/\mu-\vp_t)},\q t\in[0,\mu].
\end{equation}
All of these paths, as well as the Ricci CP, corresponds to different curves within the
two-parameter family of equations
\begin{equation}
\label{TwoParamCMEq}
\o_{\vp}^n=e^{tf_\o+c_t-s\vp}\on,
\quad
c_t:=-\log\frac1V\int_M e^{tf_\omega}\omega^n, \quad (s,t)\in A,
\end{equation}
where $A:=(-\infty,0]\times[0,1]\;\cup [0,\mu]\times\{1\}$
(see Figure \ref{FigCP}).

In the smooth setting, any one of these paths may be used to prove existence of
a KE
metric, assuming the K-energy is proper. Note that different paths
have been used to prove existence, depending on $\mu$: \eqref{CalabiCPEq}
when $\mu\le0$, and \eqref{AubinCPEq} or \eqref{DKCPEq} when $\mu>0$. Below, we will prove
existence in a unified
manner, i.e., regardless of the sign of $\mu$ or whether $\be=1$ or $\be\in(0,1)$. In fact, we show that when $\be\in(0,1/2]\cup\{1\}$
then \eqref{TwoParamCMEq} has a solution for each $(s,t)\in A$.
On the other hand, when $\be\in(0,1/2)$ the Ricci curvature of the reference metric $\o$
is unbounded from below as a corollary of Lemma \ref{LiRLemma}.
Thus, for each $(s,t)\in A$,
$$
\Ric\o_{\vp(s,t)}=(1-t)\Ric\o+s\o_{\vp(s,t)}+(\mu t-s)\o+(1-\be)[D],
$$
and this has a lower bound {\it only if} $t=1$.
Thus, on the one hand, the Chern--Lu inequality,
which requires such a lower bound is inapplicable; on the other hand,
the Aubin--Yau inequality which requires a lower bound on the bisectional curvature
of the reference geometry is inapplicable once again due to Lemma \ref{LiRLemma}.
This reasoning sifts out naturally the Ricci CP among all other possible curves in $A$.

\subsection{Short time existence}
\label{ShortTimeSubSec}

Intuitively, the Ricci continuity path \eqref{RCMEq} has the trivial solution
$\o(-\infty)=\o$ at $s=-\infty$.
Producing solutions for very negative $s$
can be considered as the continuity method analogue of
showing short-time existence for the Ricci flow.
However, it is not possible to apply the implicit function theorem
directly to obtain solutions for large negative finite values of $s$
(this observation is due to Wu, who noted that
the last displayed equation on \cite[p. 589]{TianYau1986}
is valid for $s_0>0$ but not for $s_0=0$).
 Indeed, reparametrizing \eqref{RCMEq} by
setting $\sigma = -1/s$, then the linearization of the Monge-Amp\`ere equation at $\sigma$ equals
$\sigma \Delta_{\vp(-1/\sigma)} - 1$, which degenerates at $\sigma = 0$.
More concretely, $L_\sigma:=\sigma \Delta_{\vp(-1/\sigma)} - 1$
has bounded operator norm when considered as acting from
$\Dw$ to $C^{0,\ga}_w$ only when $\sigma>0$: there is no constant $C>0$
such that $||L_0v||_{\Dw}=||v||_{\Dw}\le C||v||_{C^{0,\ga}_w}$, of course.

Thus a different method is needed to produce a solution of \eqref{RCMEq}
for sufficiently negative, but finite, values of $s$.
We present two arguments. The first, described in \S\ref{TwoTrickSubSec},
works only when $\be\in(0,1/2]\cup\{1\}$.
Wu's original argument also only works when $\be$ is in that range; in
\S\ref{NewtonSubSec} we present a generalization of Wu's argument that does not require
lower curvature bounds on $\o$ and thus is applicable for all $\be\in(0,1]$.
Finally, we remark that an interpretation of the $s\ra-\infty$ limit
in terms of thermodynamics together with a variational approach 
has been given by Berman \cite{Berman-zerotemp}.

\subsubsection{The two-parameter family trick}
\label{TwoTrickSubSec}

When $\be\in(0,1/2]\cup\{1\}$, the difficulty with applying the implicit function
theorem can be circumvented as follows \cite{JMR}.
Indeed, 
the original continuity
path \eqref{RCMEq} embeds into the two-parameter family \eqref{TwoParamCMEq}, and
it is trivial that solutions exist for the finite parameter values $(s,0)$.
Thus, while this does not show directly that our original equation has
solutions for all $(s,1)$ with $s$ sufficiently negative, it yields that result eventually,
{\it provided} we have a priori estimates for all values $(s,t)$.
This is somewhat reminiscent of adding variables or symmetries to a given equation in order
to solve it.

The reason this trick does not seem to work when $\be>1/2$ is that
it is not clear how to obtain the a priori estimates
needed to carry out the rest of the continuity argument for the two-parameter family,
unless $\be\in(0,1/2]\cup\{1\}$.
In essence then, the classical continuity path with parameter values
$\{(-\mu,t), \, t\in[0,1]\}$ may simply fail to exist within
the space of $\tildeDs$-regular \K edge potentials. It would be interesting
to understand the maximal set of values $(s,t)$ for which
a solution exists in \eqref{TwoParamCMEq},
\beq
\label{MSetEq}
M:=\{(s,t)\in A\,:\, \eqref{TwoParamCMEq} \h{ admits a
solution $\vp(s,t)\in\PSH(M,\o)\cap\Dw$} \},
\eeq
as well as analogues of this set for lower regularity classes.

\subsubsection{Newton iteration arguments}
\label{NewtonSubSec}

Thus, to handle the general case, another method must be used to obtain a solution
of \eqref{RCMEq} for some very negative value of $s$.  Wu used a
Newton iteration argument to obtain such a solution in a different setting \cite[Proposition 7.3]{Wu}.
However, his argument requires a Ricci curvature bound on the reference metric
(see \cite[p. 431]{Wu} where the expression $\Delta_\o f_\o$, that on $M\sm D$ equals
the scalar curvature up to a constant, enters), which we lack.
What follows is an adaptation of Wu's argument that requires
no curvature control on the reference metric, and thus requires more delicate estimates.
We compare our approach to Wu's in Remark \ref{WuRemark}.

Define
\[
N_{\sigma}:\Dw\ra C^{0,\gamma}_w, \ \  \  N_\s(\Phi):=\log (\o_{\sigma\Phi}^n/ e^{f_\o}\o^n)-\Phi.
\]
This is equivalent to the original \MA equation
\eqref{RCMEq} upon substituting $\sigma=-1/s$ and
$\Phi=-s\vp$.  Note that $DN_\s|_\Phi=\sigma\Delta_{\sigma\Phi}-\h{Id}$.
Now, suppose that $\s\Phi\in\calH_\o$, and, say, $s<-1$.
By the maximum principle (adapted to the edge setting by adding a barrier
function, see \S\ref{ConvergenceSubSec}), the nullspace of $DN_\s$ is trivial
provided $s < 0$. Thus, Theorem~\ref{DsThm} (i) implies that
$DN_\s|_\Phi:\Dw\ra C^{0,\gamma}_w$ is an isomorphism, with
\begin{equation}
\label{est1}
||u||_{\Dw}\le C||DN_{\s}u||_{C^{0,\gamma}_w}.
\end{equation}
Denote by $DN_\s|_{\Phi}^{-1}$ the inverse of this map on $C^{0,\gamma}_w$.

\begin{prop}
\label{NewtonProp}
Define,
$\Phi_0=0,\;
\Phi_k=(\h{Id}-DN_\s|_{\Phi_{k-1}}^{-1}\circ N_\s)(\Phi_{k-1}), \quad k\in\NN.
$
There exists $0<\sigma_0\ll 1$ and $\gamma'>0$,
such that if $\s\in(0,\s_0)$ then
 $\s\lim_{k\ra\infty}\Phi_k\in\calD_w^{0,\gamma'}\cap\PSH(M,\o)$
solves \eqref{RCMEq} with $s=-1/\s$.
\end{prop}

The proof appears in \cite[\S9]{JMR}. The crucial step is
showing that $\s\Phi_1$ has {\it small} $\,\calD_w^{0,\gamma'}$
norm, and therefore it is still a \K edge potential. We now
sketch some of the details.

When $k=1$,  $N_\s(0) = -f_\o$ and $\Phi_1=(\s\D_\o-\Id)^{-1}f_\o$.
To prove that $\sigma \Phi_1 \in \calH_\o$, it suffices to show that the pointwise norm
$|\ddbar \, \s \Phi_1|_\o$ is small, for then $\o + \sqrt{-1} \ddbar \, \s \Phi > 0$.
By Theorem~\ref{DsThm} (i), it is enough to prove that $\Delta_\o (\s \Phi_1)$ is small
in $C^{0,\gamma}_w$. However, by definition,
\[
\s \Delta_\o \Phi_1 = \s\Delta_\o (\s \Delta_\o - 1)^{-1} f_\o = \Delta_\o (\Delta_\o + s)^{-1} f_\o.
\]
The difficulty is that
$f_\o \in C^{0,\gamma}_w$ for $\gamma \in (0, 1/\beta - 1)$,
but not higher in the wedge H\"older scale.
To overcome this, we consider $f$ as
varying over a range of function spaces (some of which $f_\o$
does not belong to!), and estimate
the norm of the map
$
C^{\ell_1,\gamma_1}_w \ni f \mapsto \Delta_\o(\Delta_\o +s)^{-1} f \in C^{\ell_2, \gamma_2}_w,
$
for different values of $(\ell_j,\gamma_j)$,
and interpolate.
This eventually leads to the estimate
$||\Delta_\o(\Delta_\o + s)^{-1} f_\o||_{w; 0,\gamma''} \leq C |s|^{-\eta}$
for some $\eta > 0$ and $\gamma'' \in (0,\gamma)$,
proving that $\s\Phi_1$ is a \K edge potential for small enough $\s$.
The rest of the proof of Proposition \ref{NewtonProp} then
follows by induction.

\begin{remark}{\rm
\label{WuRemark}
It is worth comparing the above approach to Wu's original
argument. Appropriately translating
Wu's argument into our setting one would have considered the operator
$
N'_\s(\Phi):=\log (\tilde \o_{\sigma\Phi}^n/\o^n)-\Phi,
$
with $\tilde\o:=\o-\sigma\i\ddbar f_\o$.
This is clearly equivalent to the original complex \MA equation.
With this definition,
$$
\Phi_1=(\s\Delta_{\tilde\o}-1)^{-1}\log\frac{\tilde\o^n}{\on}.
$$
When a small multiple of $f_\o$ is a \K edge potential, equivalently
when the Ricci curvature of the reference $\o$ is uniformly bounded,
then $\o$ and $\tilde\o$ are uniformly equivalent for all
small $\s$. Then it is straightforward
to show $\s\Phi_1$ is small in $\Dw$, essentially from
its definition and by computing $D^2N'$.
The approach we described before was devised precisely to circumvent this lack of differentiability
of $f_\o$.

}
\end{remark}

\subsection{Convergence} 
\label{ConvergenceSubSec}

To show the convergence of the Ricci CP,  in particular implying the existence
of a KE(E) metric, one must prove, as usual, openness and closedness
of the set $M\cap(-\infty,\mu]\times\{1\}$ in
$(-\infty,\mu]\times\{1\}$ (recall \eqref{MSetEq});
indeed, since this set is nonempty by \S\ref{ShortTimeSubSec},
this implies that it is equal to $(-\infty,\mu]\times\{1\}$.
Openness follows as described in \S\ref{PositiveSubSec}.
Closedness follows from the following a priori estimate.

\begin{thm}
\label{APrioriEstTotalThm}
Along the Ricci continuity path \eqref{RCMEq},
\beq
\label{MainAPrioriEstEq}
||\vp(s)||_{\Dw}\le C,
\eeq
where $C=C(||\vp(s)||_{L^\infty(M)},M,\o,\be,n)$.
When $\mu\le0$ or the twisted Mabuchi energy is proper then
\beq
\label{C0APrioriEstEq}
||\vp(s)||_{L^\infty(M)}\le c,
\eeq
with $c$ depending only on $M,\o,\be,n$.
\end{thm}

For very negative values of $s$ this is a consequence of \S\ref{ShortTimeSubSec}.
Thus, \eqref{C0APrioriEstEq} follows from Proposition
\S\ref{UniformSubSec}, and \eqref{MainAPrioriEstEq} is a corollary
of the a priori estimates \eqref{ChernLuIneqPotentialEq}
for the Laplacian of $\vp(s)$ together with the maximum principle,
and \eqref{HolderTianEstEq} for the H\"older
semi-norm of the Laplacian of $\vp(s)$.
These estimates are described in detail in \S\ref{APrioriSubSec} below.
There is, however, one caveat in applying the maximum principle in the edge setting:
the maximum could be attained on $D$ and then, as the metric blows up there,
one cannot make sense of its Laplacian. A trick due to Jeffres \cite{Jeffres2}
is to add the barrier function $c|s|_h^\eps$ with $c,\eps>0$ small. This
is easily seen to ``push" the maximum away from $D$, while not changing the value
of the function being maximized by a whole lot:
the latter fact is obvious, while the maximum is pushed away from $D$
precisely because the gradient of the barrier function blows up near $D$. 
One can even let $c$ tend to
zero to see that the same exact estimates as in the smooth case hold.
An improved version of this maximum principle is proved in \cite[Lemma 5.1]{JMR}.

\subsection{The Ricci iteration}
\label{RISubSec}

As explained in \S\ref{MeetsSubSec},
the Ricci continuity method is motivated by the Ricci iteration,
introduced in \cite{R08} (cf. Keller \cite{Keller}). It
is then natural to go back and prove convergence of the Ricci iteration
$\{\o_{k\tau}\}_{k\in\NN}$. When $\tau=1$, this leads to a particularly natural
result:
$$
\lim_{k\ra\infty}\Ric_\be^{-k}\o=\oKE,
$$
where $\oKE$ is a KEE metric of angle $\be$, $\Ric^{-k}:=(\Ric^{-1})^k$, and $\Ric^{-1}$ is
the twisted inverse Ricci operator defined in \S\ref{StaircaseSubSec}.
It is interesting to note that when $\be=1$, results of Donaldson
\cite{Donaldson-numerical} show that $\Ric^{-1}$ can be approximated
by certain finite-dimensional approximations, and this was
further studied by Keller \cite{Keller} yielding Bergman type approximations
to KE metrics.

The convergence of the Ricci iteration when $\tau>1$, or when $\tau=1$ and
the $\alpha$-invariant is bigger than one was proved in \cite{R08} and
adapts to the edge setting once the a priori estimates
needed for the Ricci CP are established. The weak convergence when $\tau\le1$
in general was first established in \cite{BBEGZ}
(and the strong convergence then can be deduced from arguments of \cite{R08,JMR})
  using a new pluripotential
estimate from \cite{BBGZ,Berman2010} that can be stated as follows:

\begin{lem}
\label{BermanLemma}
Suppose $J(\o,\ovp)\le C$.
Then for each $t>0$ there exists $C'=C'(C,M,\o,t)$
such that $\int_M e^{-t(\vp-\sup\vp)}\on\le C'$.
\end{lem}

In particular, since the K-energy decreases along the Ricci iteration,
the properness assumption means that $J(\o,\o_{k\tau})$ is uniformly
bounded, independently of $k$. Thus, rewriting \eqref{RIEq} as
%\begin{equation}
%\label{RISecondEq}
$$
\o_{\psi_{k\tau}}^n=\on
e^{f_\o-(1-\frac1\tau)\psi_{k\tau}-\frac1\tau\psi_{(k-1)\tau}},
$$
%\end{equation}
 choosing $p$ sufficiently large
depending only on $\tau$, say
$p/3=\max\{1-\frac1\tau,\frac1\tau\}$.
Using \Kolodziej's estimate and the \Holder inequality this yields
the uniform estimate $\h{\rm osc}\, \psi_{k\tau}\le C$.
Unlike for solutions of \eqref{RCMEq}, the functions $\psi_{k\tau}$ need not be changing signs.
But an inductive argument shows that
$|(1-\frac1\tau)\psi_{k\tau}-\frac1\tau\psi_{(k-1)\tau}|\le C$
\cite[p. 1543]{R08}.
The higher derivative estimates follow as in the Ricci CP since the Ricci curvature is uniformly bounded from below along the iteration.

It is natural to also hope for similar results for a suitable K\"ahler--Ricci edge flow
\eqref{RFEq}.
For Riemann surfaces, a rather complete understanding is given by
\cite{MRSesum}, as described in \S\ref{RSSubSec}.
Different approaches to short time existence in higher dimensions
are developed by Chen and Wang \cite{ChenXWangY,ChenXWangY2,WangY},
Liu--Zhang \cite{LiuZhang}, as well as by Mazzeo and the author \cite{MR}.

\subsection{Other approaches to existence}
\label{OtherApprSubSec}

The Ricci continuity method gives a unified proof of
the classical results of Aubin, Tian, and Yau on existence of KE
metrics in the smooth setting, and naturally generalizes to give
new and optimal existence results for KEE metrics.
This was the main contribution of \cite{JMR,MR2012}, in addition
to the the linear theory and higher regularity.
Later, two other alternative approaches to existence were brought to fruition.

The first is a combination of a variational approach of Berman,
and an approximation technique of Campana, Guenancia, and P\v aun.
Guenancia--P\v aun, building on work of Campana--Guenancia--P\v aun,
developed a smooth approximation method \cite{GP,CGP} to prove
existence assuming a weak, say $C^0$, solution exists. In their scheme,
such a solution is obtained by using the variational approach
of Berman \cite{Berman2010}, under a properness assumption on
the K-energy (when $\mu\le0$ such a $C^0$ solution exists automatically by
\Kolodziej's estimate \cite{Kolodziej}).
One first approximates the reference form $\o$ \eqref{oEq} by
a particular sequence of smooth cohomologous \K forms
$\o_0+\i\ddbar\psi_\eps$ (with $\lim_{\eps\ra0}\psi_\eps=(|s|^2_h)^\be$),
and then solves the regularized \MA equation
\beq
\label{RegularizedMAEq}
(\o+\i\ddbar(\psi_\eps+\phi_\eps))^n=
\o_0^ne^{f+\mu(\psi_\eps+\phi_\eps)}(|s|^2+\eps^2)^{\be-1}.
\eeq
These equations can be solved by standard results in the smooth
setting, i.e., when $\eps>0$. It thus suffices to prove a Laplacian
estimate. When $\be\in(0,1/2]$ this follows directly from the
classical Aubin--Yau estimate \eqref{AYIneqPotentialSecondEq},
once a careful and tedious computation establishes that
the bisectional curvature of $\o_0+\i\ddbar\psi_\eps$ is
bounded from below independently of $\eps$ \cite{CGP}.
When $\be\in(0,1)$ this follows by additional clever and lengthy
computations showing that the negative contribution of the bisectional
curvature of $\o_0+\i\ddbar\psi_\eps$ in the right hand side of
\eqref{AYIneqPotentialSecondEq}, can be cancelled by adding
terms of the form $\D_\o \chi(\eps^2+|s|^2_h)$ on the left hand side,
where $\chi:\RR_+\ra\RR_+$ is a certain auxiliary function \cite{GP}.
As remarked in \cite{GP}, a somewhat similar trick appears in \cite{Brendle}
to deal with $C^3$ estimates, and is reviewed in \S\ref{Calabi3rdSubSec}.
Later, Datar--Song observed that relying on Lemma \ref{LiRLemma}
and the Chern--Lu inequality as developed in \cite[\S7]{JMR}
(Corollary \ref{CLCorProp}) one can avoid the aforementioned
lengthy computations. Either way,
the main advantage of this method over a continuity method is
that, as first observed by Berman \cite{Berman2010},
no openness argument is needed.

The second is an angle-deforming continuity method
that applies in the special case of a smooth plurianticanonical
divisor $D\in|-mK_M|$ in a Fano manifold. It was introduced
by Donaldson in lectures at Northwestern University
in 2009 and later published in \cite{Donaldson-linear} in
the case $D$ is anticanonical,
and the immediate, yet useful, extension
to the case of a plurianticanonical divisor, was noted by Li--Sun \cite{LiSun}.
Here, one constructs first a KEE of angle $\be$ along $D$
for some small $\be_0\in(0,1)$. Equation \eqref{EdgeKEEq} reads,
$\Ric{\oKE}_{\beta_0}=\be_0{\oKE}_{\beta_0}+(1-\be_0)[D]/m
$, and now we may
consider $\be_0$ as a {\it parameter} and try to deform it to a given $\be$. This was first achieved for
{\it all} small $\be_0$
by the combined results of Berman \cite{Berman2010} and \cite{JMR}:
the former shows that for all small $\be>0$ the twisted K-energy
is proper, while the latter shows that properness implies existence.
Alternatively, Berman also observes that when $\be_0=1/k$
for $k\in\NN$ sufficiently large, an orbifold KE metric can be constructed using Demailly--Koll\'ar's
orbifold
version of Tian's $\alpha$-invariant existence criterion \cite{DemaillyKollar,Tian1987}.
Next, Donaldson's openness result implies that the KEE metric of angle
$\be_0$ can be deformed to a KEE of slightly larger angle, as long as
the Lie algebra $\aut(M,D)$ is trivial. This always holds in this
Fano setting (but not in general \cite{CR})  as first observed by Berman \cite[p. 1291]{Berman2010}
(an algebraic proof of this was later given by Song--Wang \cite{SongWang}).
Finally, the recently announced results of \cite{CDSun,Tian2013}, together
with Berman's observation
that properness of the twisted K-energy
is an open property (in $\be$), can be combined to prove existence.

\section{A priori estimates for Monge--Amp\`ere equations}
\label{APrioriSubSec}

This section surveys the a priori estimates pertinent
for the study of the (possibly degenerate) complex \MA equations
\eqref{TwoParamCMEq}, both in the smooth setting ($\be=1$)
and the edge setting ($\be\in(0,1)$).
The $L^\infty$ estimate can be proved in at least three different ways,
as discussed in \S\ref{UniformSubSec} with little
dependence on the (possibly unbounded) curvature of the background geometry.
The Laplacian estimate on the other hand is quite sensitive to the latter, and
more care is needed here. We take the opportunity to give a rather self-contained
introduction to the Laplacian estimate for the complex \MA equation in
\S\ref{UniformityMetricSubSubSec}--\ref{UniformityMetricTwoSubSubSec}.
The Chern--Lu inequality was used first by Bando--Kobayashi in the 80's to
obtain a Laplacian estimate with bounded reference geometry, but 
fell into disuse since and was first systematically put to
use for a general class of \MA equations (more specifically, whenever
the solution metric has Ricci curvature uniformly bounded from below) in the author's work
on the Ricci iteration \cite{R08} and fully exploited in \cite{JMR}
to obtain estimates under a one-sided curvature bound on the reference geometry.
This is described in \S\ref{CLSubsubsec}--\ref{CLCorSubSec}.
Traditionally, except in those three articles, the Laplacian estimate
was essentially always derived using the Aubin--Yau estimate
(with the exception of \cite{Berman2010} that used the estimate from \cite{R08}).
The latter estimate depends on a lower bound on the bisectional curvature of
the reference metric, and is therefore not directly applicable for the Ricci
continuity method. However, it always seemed curious to the author that the Chern--Lu
inequality can be derived as a corollary of a general statement about holomorphic mappings,
while the Aubin--Yau estimate is classically derived using a lengthy and rather
un-enlightening computation. In \S\ref{ReverseCLSubSec} we describe a new inequality
on holomorphic {\it embeddings},
that we call the {\it reverse Chern--Lu inequality} that yields the Aubin--Yau estimate
as a corollary (\S\ref{AYAsCorSubSec}--\ref{UniformityMetricTwoSubSubSec}).
The rest of this Section describes approaches to \Holder continuity of the metric.
When $\be\in(0,1/2]\cup\{1\}$, the asymptotic expansion of \eqref{ExpansionBetaLessHalfEq}
proves that third mixed derivatives of the type $\vp_{i\b j k}$ are bounded.
Subsection \ref{Calabi3rdSubSec} indicates how to obtain a uniform estimate for such
derivatives by slightly modifying the original approach of Calabi in the smooth setting.
In general, the expansion \eqref{ExpansionAllBetaEq} shows that $\vp_{i\b j k}\not\in L^\infty(M)$
but that $\vp_{i\b j k}\in L^2$. Tian's approach to proving a uniform {\it local}
$W^{3,2}$ estimate on $\vp$ is the topic of \S\ref{TianW32SubSec}; Campanato's
characterization of \Holder spaces implies a uniform $\Dw$ estimate on $\vp$.
Finally, \S\ref{EKSubSec} describes three other approaches to $\Ds$ estimates.

\subsection{Uniformity of the potential}
\label{UniformSubSec}

\Kolodziej's estimate gives a uniform bound
on the oscillation of the solution $u$ of  $(\o_0+\i\ddbar u)^n=F\o_0^n$
in terms of $||F||_{L^{1+\eps}(M,\o_0^n)},\o_0$, and
$\eps>0$ \cite{Kolodziej}.
By \eqref{fomegaSecondEq} it suffices to take
any $\eps$ in the range $(0,\frac{\be}{1-\be})$.
Thus,
\Kolodziej's estimate (together with the normalization
along \eqref{RCMEq}) directly provides
the
$L^\infty$ bound on the
\K potential along the Ricci CP for all $s\le0$,
and this is of course enough when $\mu\le0$.
A different approach is to use Moser iteration,
as in Yau's work in the smooth setting, which
directly adapts to this setting without change.
The real challenge is
then to obtain the estimate when $s>0$.

The method used in \cite{JMR} is to prove uniform bounds on
the Sobolev and Poincar\'e constants along the Ricci CP. Then,
a standard Moser iteration argument \cite{Tianbook} gives
a uniform control on the $C^0$ norm of $\vp(s)$ in terms
of $I(\o,\o_{\vp(s)})$, which is in turn controlled
under the properness assumption.

The proof of the uniformity of the Poincar\'e inequality is
a quick consequence of the asymptotic expansion of the solutions
$\vp(s)$. This expansion precisely shows that $\vp(s)\in W^{3,2}$
and so the integration by parts in the Bochner--Weitzenb\"ock formula
is justified, and readily implies $\lambda_1(\o_{\vp(s)})\ge s$
with strict inequality for all $s<\mu$, and equality when $s=\mu$
iff there exist holomorphic vector fields on $M$ tangent to $D$
\cite[Proposition 8]{Donaldson-linear},\cite[Lemma 6.1]{JMR}.

The Sobolev inequality is trickier, but still the key is to use
the validity of the integrated form of the Bochner--Weitzenb\"ock formula.
More precisely, standard results of Bakry and others on diffusive
semigroups imply both the existence and the uniformity of a
Sobolev inequality under a general curvature-dimension condition (that precisely corresponds
to the Bochner--Weitzenb\"ock inequality holding on a class of functions) as well as some
assumptions on the algebra of functions on which the (uniform) curvature-dimension condition holds
\cite{Bakry}.
In our setting, the existence of a (possibly non-uniform in $s$) 
Sobolev inequality is easily verified by the change
of coordinate $z\mapsto\zeta$ and a covering argument.
Moreover, one can verify, using basic results on polyhomogeneity of solutions
to quasilinear elliptic equations, that the class $\Dw$ satisfies the conditions
necessary for Bakry's approach to be carried out \cite[\S6]{JMR}.
Thus, this approach furnishes a {\it uniform} in $s$ Sobolev inequality.
This approach also give a uniform diameter estimate along the Ricci continuity
method. It is interesting to note that one could also use classical 
Riemannian geometry arguments (e.g., Croke's approach
for the isoperimetric inequality \cite{Croke}, and Myers' approach for
the diameter bound \cite{Myers})
provided one knew that between every two points in $M\setminus D$
there exists a minimizing geodesic entirely contained in $M\setminus D$.
This was shown very recently by Datar \cite{Datar} building
on a result on Colding--Naber \cite{ColdingNaber}.

A completely different approach is to regularize the equation and prove
that solutions $\vp(s)$ can be approximated by smooth \K metrics whose Ricci
curvature is also bounded from below by $s$, and then use the standard results
on Sobolev bounds \cite{CDSun,Tian2013}.

Finally, as in \cite{Berman2010}, one may use
the pluripotential estimate of Lemma \ref{BermanLemma}
to obtain a $C^0$ estimate via \Kolodziej's result.
In fact, more recent and sophisticated methods yield
a H\"older estimate in this setting (see, e.g., \cite{Kolodziej-Holder,Demailly-etal-Holder}).
Furthermore, under stronger regularity assumptions on the right hand
side there is also a Lipschitz estimate due to \Blocki~\cite{BlockiGrad}.

\subsection{Uniformity of the metric, I}
\label{UniformityMetricSubSubSec}

We say that $\o,\ovp$ are uniformly equivalent
if $C_1\o\le \ovp \le C_2\o$, for some
(possibly non-constant) $C_2\ge C_1>0$.
This is implied by either
\beq
\label{FirstC2Eq}
\hbox{$n+\D_\o\vp=\tr_\o\ovp\le C_2$ and
$\det_\o\ovp\ge C_1C_2^{n-1}/(n-1)^{n-1}$},
\eeq
or,
\beq
\label{SecondC2Eq}
\hbox{$n-\D_\ovp\vp=\tr_{\ovp}\o\le 1/C_1$ and
$\det_\o\ovp\le C_1C_2^{n-1}(n-1)^{n-1}$};
\eeq
conversely, it implies
$\tr_\o\ovp\le nC_2$ and
$\det_\o\ovp\ge C_1^n$, as well as
$\tr_{\ovp}\o\le n/C_1$ and
$\det_\o\ovp\le C_2^{n}$.
Indeed, $\sum(1+\la_j)\le A$, and
$
\Pi(1+\la_j)\ge B
$
implies $1+\la_j\ge (n-1)^{n-1}B/A^{n-1}$;
conversely,
$
\Pi(1+\la_j)\ge \big(\frac1n\sum\frac 1{1+\la_j} \big)^{-n}
\ge C_1^n.
$

Let $\iota:(M,\ovp)\ra(M,\o)$ denote the identity map.
Consider $\del\iota^{-1}$
either as a map from $T^{1,0}M$ to itself, or
as a map from $\Lambda^n T^{1,0}M$ to itself.
Alternatively,
it is section of $T^{1,0\,\star}M\otimes T^{1,0}M$,
or of $K_M\otimes K_M^{-1}$, and we may
endow these product bundles with the product metric
induced by $\o$ on the first factor, and by $\ovp$ on the second factor.
Then,  \eqref{FirstC2Eq} means that the norm squared
of $\del\iota^{-1}$, in its two guises above,
is bounded from above by $C_2$, respectively
bounded from below by $C_1C_2^{n-1}/(n-1)^{n-1}$.
Similarly,  \eqref{SecondC2Eq} can be interpreted in
terms of $\del\iota$.

Now, $\det_\o\ovp=\ovp^n/\on=F(z,\vp)$ is given
solely in terms of $\vp$ (without derivatives),
it thus suffices to find an upper bound for either
$|\del\iota^{-1}|^2$ 
or $|\del\iota|^2$
(from now on we just consider maps on $T^{1,0}M$). 

The standard way to approach this is by using the maximum
principle, and thus involves computing the
Laplacian of either one of these two quantities.
The classical approach, due to Aubin \cite{Aubin1976,Aubin1978}
and Yau \cite{Yau1978}, is to estimate the first, while
a more recent approach is to estimate the second
\cite{R08,JMR,MR}, and this builds on using and finessing older
work of Lu \cite{Lu68} and Bando--Kobayashi \cite{BK88}.

We take the opportunity to explain here
both of these approaches in a unified manner since
such
a unified treatment seems to be missing in the literature.
In particular, essentially all known Laplacian estimates are seen to be a direct
corollary of the Chern--Lu inequality or its reverse form.
We now explain this in detail.

\subsection{Chern--Lu inequality}
\label{CLSubsubsec}

Let $f:(M,\o)\ra(N,\eta)$ be a holomorphic map between \K manifolds.
We choose two holomorphic coordinate charts
$(z_1,\ldots,z_n)$ and $(w_1,\ldots,w_n)$
centered at a point $z_0\in M$ and at a point $f(z_0)\in N$, respectively, such that the first
is normal for $\o$ while the second is
normal for $\eta$. In those coordinates, we consider
the map 
$$
f:z=(z^1,\ldots,z^n)\mapsto f(z)=(f^1(z),\ldots,f^n(z)),
$$
and write $\o=\i g_{i\b j}(z)dz^i\w \overline{dz^j},\;
\eta=\i h_{i\b j}(w)dw^i\w \overline{dw^j}$, and
$$
\del f|_{T^{1,0}M}
=
\frac{\del f^j(z)}{\del z^i}dz^i|_z\otimes\frac{\del}{\del w^j}\big|_{w(z)}
=
f^j_idz^i|_z\otimes\frac{\del}{\del w^j}\big|_{w(z)},$$
so
\beq
\label{fsqEq}
|\del f_{T^{1,0}M}|^2=g^{i\b l}(z)h_{j\b k}(f(z))f^j_i(z)
\overline{f^k_l(z)}.
\eeq
Thus, at $z_0$,
\beq
\label{CLFirstEq}
\begin{aligned}
\Delta_\o|\del f_{T^{1,0}M}|^2(z_0)
&=\sum_{p,q} g^{p\b q} \frac{\del^2(g^{i\b l}h_{j\b k}f^j_i
\overline{f^k_l})}{\del z^p\overline{\del z^p}}
\cr
&=\sum_p g^{p\b q}\Big[
g^{i\b l}h_{j\b k,d\b e}f^j_i\overline{f^k_l}f^d_p\overline{f^e_q}
-h_{j\b k}g^{i\b t}g^{s\b l}g_{s\b t,p\b q} f^j_i\overline{f^k_l}
+g^{i\b l}h_{j\b k}f^j_{ip}\overline{f^k_{lq}}
\Big]
\cr
&=
-\o^{\#}\otimes\o^{\#}\otimes R_\eta(\del f,\dbar f,\del f,\dbar f)
+(\Rico)^{\#}\otimes\eta(\del f,\dbar f)
+g^{p\b q}g^{i\b l}h_{j\b k}f^j_{ip}\overline{f^k_{lq}},
\end{aligned}
\eeq
where the last line is now coordinate independent.
Here $R_\eta$ denotes the curvature tensor of $\eta$ (of type $(0,4)$),
while $\o^{\#}$ denotes the metric $g^{-1}$
on $T^{1,0\,\star}M$ (i.e., of type $(2,0)$), and similarly  $(\Rico)^{\#}$
denotes the $(2,0)$-type tensor obtained from $\Rico$
by raising indices using $g$.
The last term in \eqref{CLFirstEq} is equal to $|\nabla\del f|^2$,
the covariant differential of $\del f$, a section
of $ T^{1,0}M\otimes T^{1,0}M\otimes f^\star T^{1,0}N$.
Note, finally, that $\eta(\del f,\dbar f)=f^\star\eta$,
and similarly other terms above can be expressed as pull-backs
from $N$ to $M$.

\begin{prop} {\rm (Chern--Lu inequality)}
\label{CLProp}
Let $f:(M,\o)\ra(N,\eta)$ be a holomorphic map between \K manifolds.
Then,
\beq
\label{CLInEq}
\begin{aligned}
|\del f|^2\D_\o\log|\del f|^2
&=
((\Rico)^{\#}\otimes\eta)(\del f,\dbar f)
-(\o^{\#}\otimes\o^{\#}\otimes R_\eta)(\del f,\dbar f,\del f,\dbar f)
+e(f)
\cr
&\ge
(\Rico)^{\#}\otimes\eta(\del f,\dbar f)
-\o^{\#}\otimes\o^{\#}\otimes R_\eta(\del f,\dbar f,\del f,\dbar f),
\end{aligned}
\eeq
where $e(f)=|\nabla\del f|^2-|\del f|^2|\del\log|\del f|^2|^2$.

\end{prop}

The proof follows from the previous paragraph, the
identity $u\Delta_\o\log u=\Delta_\o u-u|\del\log u|^2$,
and the Cauchy--Schwarz inequality:
indeed, since $ f$ is holomorphic,
\beq
\label{CauchyScEq}
|\del f|^2|\nabla\del f|^2\ge \langle\nabla\del f,\dbar f\rangle^2
=|\del|\del f|^2|^2,
\eeq
therefore $e(f)\ge0$. Note that $\del^2f$ and $\del f$ are
of course sections of different bundles, and so we are abusing
notation a bit when we write $\langle\del\del f,\dbar f\rangle$, however
the meaning should be clear from \eqref{fsqEq}.
Note that the right hand side of \eqref{CLInEq} can
be thought of as the Ricci curvature of the bundle
$T^{1,0\,\star}M\otimes f^\star T^{1,0}N$ equipped with the metric
$g^{\#}\otimes f^\star h$.

A few words about the history of inequality \eqref{CLInEq}.
It was shown by Lu (more generally on Hermitian manifolds) \cite{Lu68} (though
``$=$" should be replaced by ``$\ge$"
in \cite[(4.13)]{Lu68}). Chern \cite{Chern}
carried out a similar computation earlier for
$\Delta_\o |\del f|^2$
for $\del f$ considered as a map from ${\Lambda^nT^{1,0}M}$
to ${\Lambda^nT^{1,0}N}$.
As pointed out to us by Donaldson, \eqref{CLFirstEq} is
in fact a special case of the computation of Eells--Sampson
(with a slightly different sign convention)
\cite[(16)]{EellsSampson} on the Laplacian of the energy density
of a harmonic map (holomorphic maps are harmonic by
op. cit., pp. 116--118). In fact, as Eells--Sampson observe,
when $f$ is an immersion, \eqref{CLFirstEq} can also be proved using the
Gauss equations.

\subsection{A corollary of the Chern--Lu inequality}
\label{CLCorSubSec}

The next result shows how to estimate
$\D_\ovp\log|\del\iota|_{T^{1,0}M}|^2$ solely under
an upper bisectional curvature bound on the target,
and a generalized lower Ricci curvature bound on the domain.

\begin{prop}
\label{CLCorProp}
Let $f:(M,\o)\ra(N,\eta)$ be a holomorphic map between \K manifolds.
Assume that $\Ric\o\ge -C_1\o-C_2f^\star\eta$
and that $\h{\rm Bisec}_\eta\le C_3$, for some $C_1,C_2,C_3\in\RR$.
Then,
\begin{equation}
\label{FirstTraceChernLuIneq}
\Delta_{\o} \log |\del f|^2\ge -C_1-(C_2+2C_3)|\del f|^2.
\end{equation}
In particular, if $f=\iota:(M,\o)\ra (M,\eta)$ is
the identity map, and $\o=\eta+\i\ddbar\vp$ then
\begin{equation}
\label{ChernLuIneqPotentialEq}
\Delta_{\o} \big(
\log \tr_\o\eta
-(C_2+2C_3+1)\vp
\big)
 \ge -C_1-(C_2+2C_3+1)n + \tr_\o\eta.
\end{equation}

\end{prop}

Using the Chern--Lu inequality to prove a Laplacian estimate for
complex \MA equations seems to go back to Bando--Kobayashi \cite{BK88},
who considered the case $\Ric\o\ge -C_2 \eta$. Next, the case
$\Ric\o\ge -C_1\o$ first appeared in proving a priori Laplacian estimate
for the Ricci iteration \cite{R08}, where both the Bando--Kobayashi estimate
and the Aubin--Yau estimate do not work directly.
Proposition \ref{CLCorProp} combines both cases, and first appeared
in \cite[Proposition 7.1]{JMR}.

\subsection{The reverse
Chern--Lu inequality}

\label{ReverseCLSubSec}

The following is a new, reverse form of the Chern--Lu inequality.

\begin{prop}{\rm (Reverse Chern--Lu inequality)}
\label{ReverseCLProp}
Let $f:(M,\o)\ra(N,\eta)$ be a holomorphic map between \K manifolds
that is a biholomorphism onto its image.
Then,
\beq
\begin{aligned}
|\del f|^2\D_\eta\log|\del f|^2\circ f^{-1}
&=
-(\o^{\#}\otimes \Ric\,\!\eta)(\del f,\dbar f)
+((R_\o)^{\#}\otimes \eta\otimes \eta^{\#})(\del f,\dbar f,\del f^{-1},\dbar f^{-1})
+e(f)
\cr
&\ge
-\o^{\#}\otimes \Ric\,\!\eta(\del f,\dbar f)
+(R_\o)^{\#}\otimes \eta\otimes \eta^{\#}(\del f,\dbar f,\del f^{-1},\dbar f^{-1}),
\end{aligned}
\eeq
where
$e(f)=
|\nabla^{1,0}\del f|^2
-
|\del f|^2|\nabla^{1,0}_\eta\log|\del f|^2|^2
$.
\end{prop}

\bpf
Mostly keeping the notation of \S\ref{CLSubsubsec} we compute
$\D_\eta\tr_\o\eta$ with respect to two holomorphic coordinate charts,
but now we only assume $z=(z^1,\ldots,z^n)$ is normal for $g$, and let
$w=(w^1,\ldots,w^n)=f(z)$. By our assumption on $f$, $w$ is a holomorphic
coordinate on $f(M)$.
Then,
\beq
\label{ReverseCLIntermediate}
\begin{aligned}
\D_\eta|\del f_{T^{1,0}M}|^2
&=
\sum_{p,q}h^{p\b q}\frac{\del^2}{\del w^p\overline{\del w^q}}\Big[
 g^{i\b l}(z)h_{j\b k}(f(z))f^j_i(z)\overline{f^k_l(z)}
\Big]
\cr
&=
\sum_{p,q}
h^{p\b q}\Big[
 g^{i\b l}h_{j\b k,p\b q}f^j_i(z)\overline{f^k_l(z)}
-
g^{i\b t}g^{s\b l}g_{s\b t,d\b e}h_{j\b k}(f^{-1})^d_p\overline{(f^{-1})^e_q}
f^j_i\overline{f^k_l}
\Big]\cr
&=
-(\o)^{\#}\otimes \Ric\,\!\eta(\del f,\dbar f)
+(R_\o)^{\#}\otimes \eta(\del f,\dbar f,\del f^{-1},\dbar f^{-1})
+e_1(f).
\end{aligned}
\eeq
Here,
$$
h^{p\b q}h_{j\b k,p\b q}
=
h^{p\b q}(-R_{j\b kp\b q}
+h^{s\b t}h_{j\b t,p}h_{s\b k,\b q}),
$$
so
\beq\label{e1Eq}
e_1(f):=
g^{i\b l}(z)h^{p\b q}h^{s\b t}h_{j\b t,p}h_{s\b k,\b q}(w)
f^j_i(z)\overline{f^k_l(z)}=|\nabla^{1,0}\del f|^2.
\eeq
Also, by $(R_\o)^\#$ we denote the $(2,2)$-type version
of the curvature tensor.
Next,
$$
e_2(f):=
\D_\eta|\del f|^2
-
|\del f|^2\D_\eta\log|\del f|^2
=
|\del f|^2|\nabla^{1,0}_\eta\log|\del f|^2|^2.
$$
This time,
\beq
\label{CauchyScSecondEq}
e_1(f)|\del f|^2
=
|\nabla^{1,0}\del f|^2|\del f|^2\ge
\langle\nabla^{1,0}\del f,\dbar f\rangle^2
=|\nabla^{1,0}_\eta|\del f|^2|^2
=
e_2(f)|\del f|^2,
\eeq
using the K\"ahler condition.
Therefore,
$$
\begin{aligned}
|\del f|^2\D_\eta\log|\del f|^2 
\ge
-(\o)^{\#}\otimes \Ric\,\!\eta(\del f,\dbar f)
+(R_\o)^{\#}\otimes \eta(\del f,\dbar f,\del f^{-1},\dbar f^{-1}),
\end{aligned}
$$
as desired.
\epf
Note that \eqref{CauchyScSecondEq} seems to simplify, or at least cast invariantly,
Aubin--Yau's derivation, done in
coordinates, of a similar inequality, cf.
\cite[(2.15)]{Yau1978},\cite[p. 99]{Siu}.

\begin{remark}
\label{OtherCLRemark}
{\rm
The reverse Chern--Lu inequality is {\it not} the
same inequality one would obtain from the Chern--Lu
inequality 
by considering the inverse of
the identity map. In fact, the latter would yield
$$
|\del f^{-1}|^2\D_\eta\log|\del f^{-1}|^2
\ge
(\Ric\,\!\eta)^{\#}\otimes\o(\del f^{-1},\dbar f^{-1})
-\eta^{\#}\otimes R_\o(\del f^{-1},\dbar f^{-1},\del f^{-1},\dbar f^{-1}),
$$
or specifically, considering the map $\iota^{-1}:(M,\o)\ra(M,\ovp)$,
\beq\label{ChernLuInOtherDirectionIneq}
\D_\o\log(n+\D_\o\vp)\ge
(\Rico)^{\#}\otimes\ovp(\del\iota,\dbar\iota)
-\o^{\#}\otimes R_\ovp(\del\iota,\dbar\iota,\del\iota,\dbar\iota),
\eeq
which is less useful, since it is hard to estimate
the full bisectional curvature of a solution to
a complex \MA equation (which only controls the Ricci curvature).

}
\end{remark}

\subsection{The Aubin--Yau inequality as a corollary}
\label{AYAsCorSubSec}

We now demonstrate how the classical Aubin--Yau inequality
\cite{Aubin1976,Yau1978}
(see Siu \cite[p. 114]{Siu} for a comparison between the approaches
of Aubin and Yau) can
be deduced using the reverse Chern--Lu inequality. It allows
to work under somewhat complementary curvature assumptions
to those in Proposition \ref{CLCorProp}.

\begin{prop} {\rm (Aubin--Yau Laplacian estimate)}
In the above, let $f=\h{\rm id}:(M,\o)\ra (M,\eta)$ be
the identity map, and assume that
$\Ric\eta\le C_1\o+C_2\eta$
and that $\h{\rm Bisec}_\o\ge -C_3$, for some $C_1,C_2,C_3\in\RR$.
Then,
\beq\label{ReverseCLFirstEq}
\tr_\o\eta\D_\eta\log\tr_\o\eta
\ge
-n(C_1+C_3)-C_2\tr_\o\eta-C_3\tr_\o\eta\,\tr_\eta\o.
\eeq
In particular, if $\eta=\o+\i\ddbar\vp$ then
\begin{equation}
\label{AYIneqPotentialEq}
%\begin{aligned}
\Delta_{\ovp} \big(
\log \tr_\o\ovp
-(C_3+1)\vp
\big)
 \ge -n\frac{C_1+C_3}{\tr_\o\ovp}-(C_2+n(C_3+1)) + \tr_\ovp\o,
%\end{aligned}
\end{equation}
hence
\begin{equation}
\label{AYIneqPotentialSecondEq}
%\begin{aligned}
\Delta_{\ovp} \big(
\log \tr_\o\ovp
-(2C_3+C_1+1)\vp
\big)
 \ge -(C_2+n(2C_3+C_1+1)) + \tr_\ovp\o.
%\end{aligned}
\end{equation}

\end{prop}

\bpf
Proposition \ref{ReverseCLProp} implies
\eqref{ReverseCLFirstEq} by direct computation.
Since $\tr_{\ovp}\o=n-\D_{\ovp}\vp$,
this inequality is equivalent to \eqref{AYIneqPotentialEq}.
Since $\tr_\o\ovp\,\tr_\ovp\o\ge n$ (since $\tr A\, \tr A^{-1}\ge n$
for every positive matrix $A$), this last inequality
implies \eqref{AYIneqPotentialSecondEq}.
\epf

\begin{remark}
{\rm
An older inequality of Aubin \cite[p. 408]{Aubin1970}
also follows from \eqref{ReverseCLIntermediate}
when $\o$ has nonnegative bisectional curvature.
Aubin used this inequality to prove the Calabi conjecture under this
curvature assumption.
}
\end{remark}

\begin{remark}
{\rm
The reverse Chern--Lu inequality also yields, by considering
the inverse of the identity map $\iota^{-1}:(M,\ovp)\ra(M,\o)$,
\beq\label{NotUsefulReverseBackwardCL}
\D_\o\log(n-\D_\ovp\vp)\ge
-(\ovp)^{\#}\otimes \Ric\,\!\o(\del f,\dbar f)
+(R_\ovp)^{\#}\otimes \o(\del f,\dbar f,\del f^{-1},\dbar f^{-1}),
\eeq
which is not very useful for the \MA equations we consider
here for the same reasons as in Remark \ref{OtherCLRemark}.
In summary, there are four quantities one can estimate, and
the corresponding four inequalities are
\eqref{ChernLuIneqPotentialEq} (Chern--Lu),
\eqref{ChernLuInOtherDirectionIneq} (Chern--Lu backwards),
\eqref{AYIneqPotentialEq} (reverse Chern--Lu, i.e., Aubin--Yau),
and
\eqref{NotUsefulReverseBackwardCL} (reverse Chern--Lu backwards),
and it is the first and the third which are most useful.
One may also easily derive four corresponding parabolic versions of the above
inequalities: we leave the details to the reader.
}
\end{remark}

\subsection{Uniformity of the metric, II}
\label{UniformityMetricTwoSubSubSec}

The estimates of the previous paragraphs imply uniformity
of the metric under various curvature assumptions coupled with
uniform estimates on the potential and/or the volume form.
Let us now state these consequences carefully.
Part (i) in the next result is a corollary of the Chern--Lu inequality, and seems
to have first been formulated in this generality in \cite{JMR}.
Part (ii) is a corollary of the reverse Chern--Lu inequality,
and seems to be phrased in this generality for the first time here.
Part (iii) in the case
$\Ric\ovp=\Ric\o+\i\ddbar\psi_2-\i\ddbar\psi_1$
is due to P\v aun \cite{Paun2008} (whose more
quantitative formulation
appears in \cite{BermanG,BBEGZ}) 
together with Campana--Guenancia--P\v aun \cite{CGP} that allows to assume only
a lower bound on the bisectional curvature of the reference metric (without
an upper bound on the scalar curvature). 
Here we explain
how it (or rather its slight generalization to the case
$\Ric\ovp\le\Ric\o+\i\ddbar\psi_2-\i\ddbar\psi_1$),
too, follows from the reverse Chern--Lu inequality.

\begin{cor}
\label{AYLaplacianEstimateCor}
Let $\vp\in \Dw\cap C^4(M\sm D)\cap\PSH(M,\o)$.
\hfill\break
(i) Suppose that $\Ric\ovp\ge -C_1\o- C_2\ovp$ and $\h{\rm Bisec}_\o\le C_3$ on $\MsmD$.
Then
\beq\label{FirstLaplacianEstimate}
-n<\Delta_\o\vp\le
(C_1+n(C_2+2C_3+1))e^{(C_2+2C_3+1)\osc\vp}-n.
\eeq
\hfill\break
(ii) Suppose that $\Ric\ovp\le C_1\o+ C_2\ovp$ and $\h{\rm Bisec}_\o\ge -C_3$ on $\MsmD$.
Then
\beq\label{SecondLaplacianEstimate}
-n<\Delta_\o\vp\le
\frac{(C_2+n(2C_3+C_1+1))^{n-1}}{n-1}
\Big|\Big|e^{(2C_3+C_1+1)(\vp-\min\vp)}\frac{\ovp^n}{\o^n}\Big|\Big|_{L^\infty}-n.
\eeq
\hfill\break
(iii) Let $\psi_1,\psi_2\in \Dw$ with $\psi_2/C_4\in \PSH(M,\o)$.
Suppose that $\Ric\ovp\le\Ric\o+\i\ddbar\psi_2-\i\ddbar\psi_1$
and $\h{\rm Bisec}_\o\ge -C_3$ on $\MsmD$.
Then,
\beq\label{ThirdLaplacianEstimate}
-n<\Delta_\o\vp\le
\frac{(nA)^{n-1}}{n-1}
\Big|\Big|e^{(A\vp+\psi_2-\min(A\vp+\psi_2))}\frac{\ovp^n}{\o^n}\Big|\Big|_{L^\infty}-n,
\eeq
where $A:=1+C_4+C_3+\frac{|\inf\Delta_\o\psi_1|}n$.
\end{cor}

\bpf
(i)
By \eqref{ChernLuIneqPotentialEq},
$\tr_\o\ovp(p)\le C_1+n(C_2+2C_3+1)$,
where $\log \tr_\o\ovp-(C_2+2C_3+1)\vp$ is maximized
at $p\in M\sm D$, proving \eqref{FirstLaplacianEstimate}.
Notice that here we used that the maximum is always attained on $M\sm D$,
as explained in \S\ref{ConvergenceSubSec} (see \cite[\S7]{JMR} for details) 
by using barrier functions of the form $\eps_1|s|^{\eps_2}$.

\noindent (ii)
By \eqref{AYIneqPotentialSecondEq},
$\tr_\ovp\o(p)\le C_2+n(2C_3+C_1+1)$,
where $\log \tr_\o\ovp-(2C_3+C_1+1)\vp$ is maximized
at $p\in M\sm D$.
But,
$\tr_\o\ovp\frac{\o^n}{\ovp^n}\le\frac1{n-1}(\tr_\ovp\o)^{n-1}$.
Thus,
$$
\baeq
\max_M\,\tr_\o\ovp
& \le \tr_\o\ovp(p)e^{(2C_3+C_1+1)(\vp(p)-\min\vp)}
\cr
& \le
e^{(2C_3+C_1+1)(\vp(p)-\min\vp)}\frac{\ovp^n}{\o^n}(p)
\frac1{n-1}(C_2+n(2C_3+C_1+1))^{n-1}.
\eaeq
$$

\noindent (iii) Proposition \ref{ReverseCLProp} implies
$$
\tr_\o\ovp\Delta_{\ovp}\log\tr_\o\ovp\ge
-nC_3-C_3\tr_\o\ovp\tr_\ovp\o-s_\o+\D_\o\psi_1-\D_\o\psi_2,
$$
where $s_\o$ is the scalar curvature of $\o$.
However, since we do not want to assume an upper bound for the scalar 
curvature (e.g., if $\h{\rm Bisec}_\o\le C_5$
then $-s_\o\ge -n(n+1)C_5$ \cite[pp. 168--169]{KobNom2}), we observe that Proposition \ref{ReverseCLProp} 
also implies
\beq\label{RCLCorCGPEstEq}
\tr_\o\ovp\Delta_{\ovp}\log\tr_\o\ovp\ge
R_\o^{\#}\otimes \ovp\otimes\ovp^{\#}-s_\o+\D_\o\psi_1-\D_\o\psi_2.
\eeq
Recall that here the first term on the right hand side is defined 
as the contraction of the curvature tensor of $\o$ that appears in 
\eqref{ReverseCLIntermediate}.
As in \cite[Lemma 2.2]{CGP}, one can then combine the terms depending on 
the curvature of $\o$, as we now explain. Indeed, this is most easily
seen upon diagonalization (choosing normal holomorphic coordinates so that at
a given point $p$, $\o$ is represented by the identity matrix,
while $\ovp$ is represented by a diagonal matrix $\diag\{\l_1,\ldots,\l_n\}$
while $\ovp^{\#}$ is represented by $\diag\{1/\l_1,\ldots,1/\l_n\}$).
Then,
$$
R_\o^{\#}\otimes \ovp\otimes\ovp^{\#}-
s_\o
=\sum_{j,k}\Big(\frac{\l_j}{\l_k}-1\Big)R_{j\b j,k\b k}
=\sum_{j\le k}\Big(\frac{\l_j}{\l_k}+\frac{\l_k}{\l_j}-2\Big)R_{j\b j,k\b k}
\ge -C_3\sum_{j\le k}\Big(\frac{\l_j}{\l_k}+\frac{\l_k}{\l_j}-2\Big),
$$
where, by assumption, the numbers $R_{j\b j,k\b k}$ are bounded from
below $-C_3$ and 
$R_{j\b j,k\b k}=R_{k\b k,j\b j}$ by the symmetries of the curvature tensor.
But now,
$$
\sum_{j\le k}\Big(\frac{\l_j}{\l_k}+\frac{\l_k}{\l_j}-2\Big)
=
\sum_{j,k}\Big(\frac{\l_j}{\l_k}-1\Big)
=
-n(n+1)+\sum_{j,k}\frac{\l_j}{\l_k}
\le 
-n(n+1)+\sum_{j}\frac{1}{\l_j}\sum_k\l_k,
$$
and this last expression equals 
$-n(n+1)+\tr_\ovp\o\tr_\o\ovp(p)$.
Thus, returning to \eqref{RCLCorCGPEstEq}, we now deduce, 
$$
\Delta_{\ovp}\log\tr_\o\ovp\ge
(n(n+1)C_3+\D_\o\psi_1-\D_\o\psi_2)/\tr_\o\ovp
-C_3\tr_\ovp\o.
$$
Since $\tr_\o\ovp\,\tr_\ovp\o\ge n$,
$
\Delta_{\ovp}\log\tr_\o\ovp\ge
-[C_3+\frac1n|\inf\Delta_\o\psi_1|]\tr_\ovp\o-\D_\o\psi_2/\tr_\o\ovp,
$
hence,
$$
\Delta_{\ovp}\Big(\log\tr_\o\ovp
-\Big[1+C_3+\frac{|\inf\Delta_\o\psi_1|}n\Big]\vp
\Big)
\ge
\tr_\ovp\o-n\Big[1+C_3+\frac{|\inf\Delta_\o\psi_1|}n\Big]
-\frac{\D_\o\psi_2}{\tr_\o\ovp}.
$$
So far the proof followed that of (ii). P\v aun's additional observation
is that the term $\D_\o\psi_2/\tr_\o\ovp$ can be controlled
even though it has the `wrong' sign, under the plurisubharmonicity assumption.
Indeed, since $0\le C_4\o+\i\ddbar\psi_2$ then each eigenvalue of this
nonnegative form (that we take with respect to $\ovp$)
is controlled by the sum of the eigenvalues, i.e.,
$C_4\o+\i\ddbar\psi_2\le (C_4\tr_\ovp\o+\Delta_\ovp\psi_2)\ovp$.
Taking the trace of this inequality with respect to $\o$,
$$
-\frac{\Delta_\o\psi_2}{\tr_\o\ovp}\ge \frac{nC_4}{\tr_\o\ovp}-C_4\tr_\ovp\o-\Delta_\ovp\psi_2
\ge -C_4\tr_\ovp\o-\Delta_\ovp\psi_2.
$$
Thus,
\beq
\begin{aligned}
\Delta_{\ovp}\Big(\log\tr_\o\ovp
-\Big[1+C_4+C_3+\frac{|\inf\Delta_\o\psi_1|}n\Big]\vp
+\psi_2\Big)
&\ge
\tr_\ovp\o
\cr
&\mskip-150mu-n\Big[1+C_4+C_3+\frac{|\inf\Delta_\o\psi_1|}n\Big].
\end{aligned}
\eeq
Arguing as in (ii) proves \eqref{ThirdLaplacianEstimate}.
\epf

In all three cases, it follows that
$\frac1C \o\le \o_{\vp(s)}\le C\o$, with $C$ depending only on
the constants appearing in
\eqref{FirstLaplacianEstimate}--\eqref{ThirdLaplacianEstimate} and
in \eqref{FirstC2Eq}--\eqref{SecondC2Eq}, with the precise
dependence computed in \S\ref{UniformityMetricSubSubSec}. Of course,
if one assumes equalities in the expressions for $\Ric\ovp$ instead
of inequalities, then one can express the ratio of the volume
forms more explicitly to make the estimates
\eqref{SecondLaplacianEstimate}--\eqref{ThirdLaplacianEstimate}
more explicit.

\subsection{Uniformity of the connection: Calabi's third derivative estimate}
\label{Calabi3rdSubSec}

The term $e_1(f)=|\nabla\del f|^2$ of \eqref{e1Eq}
is the norm squared of the
connection associated to $g^{\#}\otimes f^\star h$.
on $T^{1,0\,\star}M\otimes f^\star T^{1,0}N$.
When $f=\h{id}$ it equals $|\nabla^{1,0}\ddbar\vp|^2$,
with the norm taken with respect to $\ovp\otimes\ovp\otimes\o$.
But assuming that the metrics $\ovp$ and $\o$ are uniformly equivalent,
this term is uniformly equivalent to a term obtained by using the norm
associated to $\ovp$ alone, namely
$S:=|\nabla^{1,0}\del f|_{\ovp}^2$.
Under appropriate bounds on $R_\o$ and $\Ric\ovp$ it follows
from \eqref{ReverseCLIntermediate}
that $\D_\ovp \tr_\o\ovp\ge C_1S-C_2$, with $C_i$ depending
also on the equivalence between $\o$ and $\ovp$. Standard computations
going back to Calabi \cite{Calabi1958}
also show that $\D_\ovp S\ge -C_3'|R_\o|S-C_4'|DR_\o|S^{1/2}
\ge -C_3S-C_4$, with constants
depending on bounds on $R_\o$ and its covariant derivative.
Thus, $\D_\ovp(S+\frac{C_3+1}{C_1}\tr_\o\ovp)\ge S-\frac{C_3+1}{C_1}C_2-C_4$.
Using the maximum principle now yields a bound on $S$.
This summarizes the proof in the smooth setting
\cite[p. 410]{Aubin1970},\cite[\S3]{Yau1978},\cite{PhongSesumSturm}.

In the edge setting, $R_\o$ and its covariant derivative are no longer
bounded, and so this approach has somewhat limited applicability.
For instance, when $\be\in(0,1/2)$, Brendle observed that
$|\nabla R_\o|\le \frac C{|s|_h^{\eps-\be}}$ for some $\eps>0$. Thus an inequality
from the previous paragraph implies
$\D_\ovp S\ge -C_3S-\frac{C_4'}2(|DR_\o|^2+S)\ge
-(C_3+\frac{C_4'}2)S-\frac{C_4'}2\frac {C^2}{|s|_h^{2\eps-2\be}}$.
But since $\D_\ovp|s|_h^{2\eps}\ge C_5|s|_h^{2\eps-2\be}-C_6$, the maximum
principle can be applied, this time to $S+\frac{C_3+1}{C_1}\tr_\o\ovp+|s|_h^{2\eps}$,
to conclude \cite{Brendle}.

\subsection{Tian's $W^{3,2}$-estimate}
\label{TianW32SubSec}

A general result due to Tian  \cite{tian83}, proved in his M.Sc. thesis, gives a local a priori estimate in $W^{3,2}$
for solutions of both real and complex Monge--Amp\`ere equations under the assumption that the solution
has bounded real or complex Hessian and the right hand side is at least H\"older. By Campanato's
 classical integral characterization
of H\"older spaces
\cite{Campanato,HanLin}
this implies a uniform H\"older estimate on the Laplacian. This result can be seen as an alternative to
the Evans--Krylov theorem (and in fact appeared independently around the same time).
Let $B_1\subset\CC^n$ be the unit ball. Consider the equation
\beq
\label{MAlocalEq}
\det[u_{i\b j}]
=
e^{F-cu}, \q \h{on $B_1$}.
\eeq
The following is due to Tian. The proof in  \cite{tian83} is written for the real \MAno,
but applies equally to the complex \MAno.

\begin{thm}
\label{TianEstThm}
Suppose that $u\in C^4\cap \PSH(B_1)$ satisfies \eqref{MAlocalEq}.
For any $\gamma \in (0,1)$, there exists and $C> 0$ such that
for any $0 < a < 1/2$,
\begin{equation}
\label{HolderTianEstEq}
\int_{B_a} |u_{i\b j k}|^2 \le C a^{2n-2+2\gamma}.
\end{equation}
Thus, $||\vp_{i\bar j}||_{C^{0,\gamma}(B_{1/4})}\le C'$.
The constants $C,C'$ depend only on $\gamma,\beta,\o,n$, $||u_{i\b j}||_{L^\infty(B_1)},$
$||F||_{C^{0,\gamma}(B_1)}$,
and $||u||_{L^\infty(B_1)}$.
\end{thm}

The main ideas are as follows.
First, an easy computation shows
that $(\det[u_{k\b l}]u^{i\b j})_i=0$ \cite{tian83,WJ}.
Consider the \MA equation
$\det[u_{i\b j}]=h$. Taking the logarithm and
differentiating twice, multiplying by $h$,
and using the previous identity, yields
\begin{equation}
\label{Defineh2Eq}
-hu^{i\b s}u^{t\b j}u_{t\b s,\b l}u_{i\b j,k}+(hu^{i\b j}u_{k\b l,i})_{\b j}=h_{k\b l}-h_kh_{\b l}/h,
\q \h{for each $k,l$}.
\end{equation}
Thus, the \MA equation roughly becomes a second order {\it system} of
equations in divergence form for the Hessian, with a quadratic nonlinearity, resembling the harmonic map equation.
In the harmonic map setting, a result of Giaquinta--Giusti shows that a bounded weak
$W^{1,2}$-solution is necessarily H\"older \cite{GiaquintaGiusti}.
Since we already have bounds on the real/complex Hessian, the situation is quite analogous.
As shown by Tian, the \MA equation can be treated in just the same way, proving
Theorem \ref{TianEstThm}. A key difference between the two settings is
dealing with a system, so some extra algebra facts are needed.
This method makes clear that no additional curvature assumptions on the reference geometry are needed
for this estimate. 

The discussion so far was in the absence of edge singularities.
However, a very nice feature of the above result is that 
since its proof involves {\it integral} quantities, they carry over
verbatim to the edge situation $C_\be\times\CC^{n-1}$ (recall \eqref{CbetaEq}), modulo
only one very minor difference \cite[Theorem B.1]{JMR}: since the vector field
$V_1=z_1^{1-\be}\del_{z_1}=\del_\ze$ is multivalued, then if $f$ is a smooth
function in $(z_1,\ldots,z_n)=(z_1,Z)$ then choosing any branch and considering
the {\it function} $f_1=V_1f$ on the model wedge, $f_1$ satisfies the boundary condition
\begin{equation}
\label{eq:boundary-1}
f_1( r e^{\sqrt{-1} 2\pi \beta}, Z)\, =\, e^{\sqrt{-1} 2\pi (1-\beta)} f_1(r,Z).
\end{equation}
Thus, the Sobolev inequality satisfied by such functions deteriorates as $\be$
approaches 1 (but this is harmless if $\be$ is fixed, for instance). Moreover,
harmonic functions with such boundary conditions satisfy the estimate
\begin{equation}\label{eq:har-5}
||d h||^2_{L^2(C_\beta(a),\o_\be)} \le C a^{2n-4 + 2\beta^{-1}}
||d h||^2_{L^2(C_\beta(1),\o_\be)},
\end{equation}
instead of the usual one with $a^{2n}$. Again, this is harmless, and the only
effect it has is to restrict the range of possible H\"older exponents
in the following immediate corollary of Theorem \ref{TianEstThm}  \cite{tian83,JMR}.
Consider the singular equation
\beq
\label{MAlocalSingEq}
\det[u_{i\b j}]
=
e^{F-cu}|z_1|^{2\be-2}, \q \h{on $B_1\setminus\{z_1=0\}$}.
\eeq

\begin{cor}
\label{TianEstConicCor}
Suppose that $u\in C^4(\{z_1\not=0\})\cap \PSH(B_1)$
satisfies \eqref{MAlocalSingEq}.
For any $\gamma \in (0,\frac1\be-1)\cap(0,1)$, there exists a $C> 0$ such that
for any $0 < a < 1/2$,
\begin{equation}
\label{eq:main}
\int_{B_a} |u_{i\b j k}|^2 \le C a^{2n-2+2\gamma}.
\end{equation}
Here $f_i:=V_if$, where
$V_1=z_1^{1-\be}\del_{z_1}=\del_\ze$ (a choice of one branch), $V_2=\del_{z_2},\ldots, V_n=\del_{z_n}$.
Thus, $||u_{i\bar j}||_{C^{0,\gamma}}\le C'$.
The constants $C,C'$ depend only on $\gamma,\beta,\o,n,||u_{i\b j}||_{L^\infty(B_1)},
||F||_{C^{0,\gamma}(B_1)}$,
and $||u||_{L^\infty(B_1)}$.
\end{cor}

\subsection{Other H\"older estimates for the Laplacian}
\label{EKSubSec}

Next, we describe several other
H\"older estimates on the Laplacian of
a solution of a complex \MA equation.
The first is a weaker $\De$ estimate that follows
the Evans--Krylov method. The other estimates are
alternative approaches to the $\Dw$ estimate descried in \S\ref{TianW32SubSec}
(we also mention the approach in \cite{CalamaiZheng}
under the restriction $\be\in(0,2/3)$).

\subsubsection{The complex Evans--Krylov edge estimate}

Compared 
with the approach presented in \S\ref{TianW32SubSec},
perhaps a more well-known approach to the
H\"older estimate on the Laplacian of
a solution of a complex \MA equation is an adaptation of
the Evans--Krylov estimate to the complex setting.
This has been carried out for the original formulation \cite{Siu}
and in divergence form \cite{WJ}; the latter requires slightly
less control on the right hand side of the equation than the former
(see also \cite{Bl}). Let us
then concentrate on the relevant adaptation to the singular setting.

The standard complex Evans--Krylov estimate adapts rather
 easily to the edge setting to give  a uniform $\De$
estimate in terms of a Laplacian estimate.
The key point is to use properties of the edge H\"older spaces under rescaling,
and a Lipschitz estimate on the right hand side that is
valid for the Ricci continuity path (and fails for some other
paths). We present the result, closely following \cite[\S8]{JMR}.

The Evans--Krylov technique is local; we may thus
concentrate on the neighborhood of $D$
where $r\le 1$. We cover this region with {\it Whitney cubes}
\[
W = W_R(y_0):=\{ p\in M\sm D\,:\, |y(p)-y_0|<R,\quad \th(p)\in I_{\pi \be} ,\quad r(p)\in(R,2R)\} \subset M \sm D,
\]
where $R>0$ and $y_0\in D$. Here $I_{\pi\be}$ denotes any interval in $S^1_{2\pi \be}$ of length $\pi \be$, so
each $W_R(y_0)$ is simply connected.  Clearly $\{r \leq 1\} \setminus D$ is covered by the union of such cubes.
Our goal is to show that there exists a fixed $\gamma'\in(0,1)$ and a uniform $C>0$ such that
$||P_{ij}\vp||_{\calC^{0,\gamma'}_e(W)}\le C$ for all such Whitney cubes $W$ (here $P_{ij}$ are the special
second order operators discussed in \S 2.1, cf.\ the discussion in \S 3.4). Taking the supremum
over $W$ gives the uniform estimate $[\vp]_{\calD^{0,\gamma'}_e}\le C'$.

The key point here is that 
$\log r$ and $y_i/r$ are distance functions
for the complete metrics $(d\log r)^2$ and $|dZ|^2/r^2$,
respectively, and the model metric is the product of these two, i.e.,
$\o_\be/r^2=(d\log r)^2+d\th^2+|dZ|^2/r^2$. Thus, we may in fact restrict to such cubes
when computing the \Holder norm. Indeed, if the supremum in the definition
of the \Holder seminorm was nearly attained for two points $p,q\in\MsmD$ not contained in one such cube then either $r(p)/r(q)>2$ or $r(q)/r(p)>2$. But then
the distance between $p$ and $q$ with respect to $\o_\be/r^2$ would be
bigger than some fixed constant, and the \Holder seminorm would then
be uniformly controlled by the $C^0$ norm, which is, by assumption, already bounded.
Similarly we see that
$|y(p)/r(p)-y(q)/r(q)|$ must be uniformly bounded.

This property of $\o_\be/r^2$ manifests itself in another way in
the proof. Namely, denoting $\Sla(r,\th,y)=(\la r,\th,\la y)$
the dilation map then
$$
||S_\lambda^\star f||_{C^{0,\gamma}_e(W_1(y_0))}
=
||f||_{C^{0,\gamma}_e(W_\lambda(y_0))}.
$$
Thus, we can perform all estimates on a cube of fixed size.

Finally, the last crucial observation 
is that for any $\vp\in\calH^e_\o$, the metric
$\ovp$ when viewed in a ``microscope" looks
(up to perhaps
a linear transformation of the original coordinates)
very close to the model (product) metric $\o_\be$.
In particular, given $\eps>0$, there exists $\lambda_0=\lambda_0(\eps,\vp)$ such
that for all $\lambda>\lambda_0$,
\begin{equation}
\label{RescaledVolEq}
(1-\eps) |d\vec{Z}|^2 \leq \la^{-2}S_\lambda^\star\o(\vec Z) \leq (1+\eps) |d\vec{Z}|^2,
\end{equation}
where
$\vec{Z} = (Z_1, \ldots, Z_n)$
are holomorphic coordinates
on $W_1$, where $\la \vec Z=\vec\zeta=(\zeta,z_2,\ldots,z_n)$,
are the original coordinates on $W_\la$.

To summarize, when we work in the rescaled cubes 
(that we may assume are of fixed
size) and use the coordinates $\vec Z$, the  
``rescaled pulled-back" metric is essentially equivalent to the model Euclidean metric
\eqref{RescaledVolEq}.
Fortunately, also the complex \MA equation we are trying to solve
transforms very nicely under this same rescaling coupled with pull-back
under the dilation map. So, at the end of the day, one may
simply apply the standard Evans--Krylov argument on this cube of fixed
size (that is disjoint from $D$). There is one small caveat, however.
In the (divergence form of the) Evans--Krylov argument one differentiates the \MA equation twice
to obtain a differential inequality and one must control the right
hand side in $C^{0,1}$.
More precisely, provided then that we can estimate the Lipschitz norm of the right hand side
$e^{f_\o-s\vp}\on$, we can
carefully put all these observations together to prove
a uniform $\De$ estimate by directly applying
the divergence form complex Evans--Krylov estimate
to the (uniformly elliptic) metric $\la^{-2}S_\lambda^\star\o$
on the (fixed size) cube $W_1(y_0)$, and to the
local \K potential $\la^{-2}(\psi+\vp)\circ S_\la$
(here $\psi$ is a local potential for $\o$).
The required aforementioned Lipschitz estimate is proved in
\cite[Lemmas 4.4,8.3]{JMR} and can be summarized as follows.
\begin{lem}
\label{hLipschitzLemma}
Let $\log h(s,t):= \log F+\log\det[\psi_{i\bar j}]=tf_\o+c_t-s\vp+\log\det[\psi_{i\bar j}]$, with $s>S$. Then
the following estimates hold with constants independent of $t,s$:\hfill\break
(i) For $\be \leq 2/3$, $||h(s,t)||_{w; 0,1} \le C=C(S,M,\o,\be,||\vp(s,t)||_{\calC_w^{0,1}})$.
\hfill\break
(ii) For $\be \le 1$,
$||h(s,1)||_{w; 0,1} \le C=C(S,M,\o,\be,||\vp(s,1)||_{\calC_w^{0,1}}).$
\end{lem}

The point of the proof of this lemma is that $f_\o$
and $\det[\psi_{i\b j}]$ are in $C^{0,\frac2\be-2}_w$,
and therefore so is $h$. This is however in $C^{0,1}_w$
only when $\be\le2/3$. Fortunately, the combination
$f_\o+\log\det[\psi_{i\b j}]$ is nevertheless always
in $C^{0,1}_w$. On the other hand, this crucial cancellation
is false
for $tf_\o+\log\det[\psi_{i\b j}]$ when $t<1$!
Thus, we must use the Ricci continuity path  in
this juncture.

To summarize, we have:

\begin{thm}
\label{HolderEstimateThm}
Let $\vp(s)\in \De\cap \calC^4(M\sm D)\cap \PSH(M,\o)$ be a solution to \eqref{RCMEq} with $s>S$
and $0 < \beta \leq 1$. Then
\begin{equation}
\label{HolderTwoEstimateEq2}
||\vp(s)||_{\De}\le C,
\end{equation}
where $\gamma > 0$ and $C$ depend only on $M,\o,\beta,S$ and $||\Delta_\o\vp(s)||_{C^0},||\vp(s)||_{C^0}$.
\end{thm}

\subsubsection{A harmonic map type argument vs. a Schauder type argument}

Equation \eqref{RescaledVolEq} was simply a consequence of the definition of a
\K edge metric. However, it turns out that under some geometric assumptions
it is possible to obtain a priori control on $\lambda(\eps,\vp)$.
Intuitively, of course, if along a family $\{\o_{\vp_j}\}$ of \K edge
metrics $\sup_j\lambda(\eps,\vp_j)=\infty$ then there exists $p\in D$ such that
 the family of metrics, while being bounded, fails to be uniformly continuous up to
the boundary at $p\in D$. This kind of behavior can be easily ruled out by the existence of
a unique tangent cone at $p$. In \cite[II]{CDSun} this is proved
when $\o_{\vp_j}$ are KEE metrics of angle $\be_j<\be_\infty<1$.
This can be generalized to \K edge metrics with a uniform lower bound
on the Ricci curvature \cite{Tian-unpublished}. One idea is to
notice that a rescaled limit will be Ricci flat and then
use results of \cite{CCTian,CDSun,Tian2013}.

Be it as it may, this improved control on the metric allows to
upgrade to \Holder bounds. We explain two approaches.

The first is similar in spirit to the proof of Theorem \ref{TianEstThm}
and is due to Tian \cite{Tian-unpublished}.
In the following paragraph all norms, covariant derivatives and
Laplacian are with respect to $\o_\be$.
Thus, let $v$ be the unique $\o_\be$-harmonic (1,1)-form
equal to $\ovp$ on $\del B_a(y)\subset U$, where $U$ is a neighborhood
of $D$ (this neighborhood is the same for the whole family
of potentials $\vp$ we are considering)
on which $(1-\eps)\o_\be<\ovp<(1+\eps)\o_\be$.
Since $\o_\be$ is $\o_\be$-harmonic,
then $\Delta |v-\o_\be|^2=|\nabla(v-\o_\be)|^2\ge0$. The maximum principle
then gives that $-\eps\o_\be<v-\o_\be<\eps\o_\be$ on $B_a(y)$, and thus
also $-\eps\o_\be<v-\ovp<\eps\o_\be$ on $B_a(y)$.
Multiplying \eqref{Defineh2Eq} by $\hat\o:=v-\ovp$ and integrating directly implies
that $||\nabla \hat\o||^2_{L^2(B_a(y))}$ is controlled from above
by $C\eps||\nabla\ovp||^2_{L^2(B_a(y))}+Cr^{2n}$.
Now, since
$
||\nabla v||^2_{L^2(B_\sigma(y))}
\le
C\left(\frac\sigma a\right)^{2n-4+\frac2\be}
||\nabla v||^2_{L^2(B_a(y))}
$,
and using Dirichlet's principle $||\nabla v||_{L^2(B_a(y))}\le
||\nabla \ovp||_{L^2(B_a(y))}$, we conclude that
$$
||\nabla \ovp||^2_{L^2(B_\sigma(y))}
\le
C\left(\eps+\left(\frac\sigma a\right)^{2n-4+\frac2\be}\right)
||\nabla \ovp||^2_{L^2(B_a(y))}+Ca^{2n}.
$$
In fact, these arguments are very similar to the ones that go into
the proof of Theorem \ref{TianEstThm}; the only difference is
in showing the smallness of $||\nabla \hat\o||^2_{L^2(B_a(y))}$. The latter
proof uses completely elementary tools (Moser iteration, essentially).
At any rate, given the inequality above, it is standard to
show the estimate \eqref{eq:main}. The arguments above assume that
$F_{i\b j}\in L^\infty$, but a close examination shows that had
$F$ only been assumed \Holder than by elementary arguments one
must replace the term $Ca^{2n}$ by a slightly more singular term
$Ca^{2n-\delta}$, however still sub-critical. We omit the elementary details.

The second approach is to cleverly use the Schauder estimate for $\Delta$.
A key point is that given any $\delta\in(0,1)$ there exists a sufficiently small
ball $B_a(y)$, with $a$ uniformly positive, such that
$$
||\det[u_{i\b j}]-\Delta u||_{C^{0,\gamma}_w}\le \delta[u]_{\Dw}.
$$
This can be proved by elementary properties of the function $\det(\,\cdot\,)\!-\tr(\,\cdot\,)$
near the identity
on the space of positive Hermitian matrices \cite[Lemma 2.2]{tian83},\cite[II]{CDSun}.
By the \MA equation we know that $\det[u_{i\b j}]$ is uniformly controlled. Thus,
by the triangle inequality and Schauder estimate of Theorem \ref{DsThm} (i),
$$
||\Delta u||_{C^{0,\gamma}_w}\le \delta[u]_{\Dw}+C\le C\delta
(||\Delta u||_{C^{0,\gamma}_w}+||u||_{C^{0,\gamma}_w}+1).
$$
Choosing $\delta$ small enough gives a uniform bound on $||\Delta u||_{C^{0,\gamma}_w}$.
In the argument above we have been particularly sloppy in
keeping track of the scaling; we refer to \cite[II]{CDSun} for details.

\subsubsection{An approximation by orbifolds argument}

In the very recent revision of the paper of Guenancia--P\v aun
(that appeared during the final revision of the present article)
one finds a yet different approach to the $\Dw$ estimate, that
we only attempt to sketch briefly, restricting for simplicity
to the case $D$ is smooth. First, the authors 
assume that $\be$ is rational, namely $\be=p/q$ for $p,q\in\NN$
with no common nontrivial divisors. Then, working on the ramified
$q$ cover essentially reduces one to the situation of edge metrics
of angles $2\pi p$. In fact, under this cover, that amounts to the
map $z\mapsto z^{1/q}$, the metric $|z|^{2\be-2}|dz|^2$ pulls-back
to $q^2|w|^{2p-2}|dw|^2$ (substitute $z=w^q$). When $p=1$ this
solves the problem, of course (the orbifold case). When $p\in\NN$,
the authors show that all the tools from the standard complex Evans--Krylov
theorem have appropriate analogues in this large angle regime.
Namely, a Sobolev inequality with an appropriate constant and
a Harnack inequality (for such degenerate (vanishing along $D$
to possibly high order) metrics). 
The key point here is that these inequalities
do not break down when $p$ tends to infinity. Then, the authors
express the \MA equation in terms possibly singular vector
fields $\frac1q w^{1-p}\del_w, \del_{z_2},\ldots,\del_{z_n}$,
to obtain the desired conclusion. Finally, the authors
approximate an arbitrary $\be\in(0,1)$ by rational numbers,
and approximate the initial \MA equation by \MA equations 
with $\be$ replaced by those rational numbers. By stability
results for the complex \MA operator the solutions of these
approximate equations will converge to the solution of the original
equation. Thus, taking a limit, the $\Dw$ estimate carries over
to the solution of the original equation.

\section{The asymptotically logarithmic world}
\label{AsympSection}

As discussed in \S\ref{KEEEqSubSec},
a basic obstruction to existence of KEE metrics
is the cohomological requirement that the $\RR$-divisor
\beq\label{KbetaDefEq}
K^\be=K^\be_M:=K_M+D_\be:=K_M+\sum_{i=1}^r(1-\be_i){D_i}
\eeq
satisfy
\beq\label{KMEq}
\h{$-K^\be_M$ equals $\mu$ times an ample class, with $\mu\in\RR$.}
\eeq
Here, $\be:=(\be_1,\ldots,\be_r)\in(0,1]^r$, $M$ is smooth,
$D\ne0$ has simple normal crossings (as we will assume throughtout
this whole section),
and $K_M^\be$ is sometimes referred to as
the twisted canonical bundle associated
to the triple $(M,D,\be)$.
In this section we will be interested in
{\it classification questions}, perhaps the most
basic of which is:

\begin{question}\label{ClassifyQ}
What are all triples
$(M,D,\be)$ for which \eqref{KMEq} holds?
\end{question}

When $\mu\le0$, according to Theorem \ref{KEEExistenceThm}, such a classification is tantamount
to a classification of KEE manifolds of nonpositive Ricci curvature.
When $\mu>0$, such a classification would yield
a class of manifolds containing the KEE manifolds.
Narrowing this class down then of course depends
on notions of stability that are a further
challenging obstacle---more on that in \S\ref{logCalabiSec}.

Question \ref{ClassifyQ} is, of course,
too ambitious
in the sense that even when all $\be_i=1$
and $M$ is smooth there is no complete classification,
or list, of projective manifolds satisfying
\eqref{KMEq}, unless $\mu>0$ and $n$ is small.
In particular, when $\mu<0$, which is a subset of the world of
``general type" varieties, a classification is quite hopeless;
in \S\ref{LogNegSubSec} we will review what can still be said
when $n=2$.
Thus, we will largely concentrate on the case $\mu>0$
and further restrict to the small angle regime where
some classification can be achieved,
that furthermore has
interesting geometric consequences.
The small angle, or asymptotically logarithmic, regime,
can be thought of as the other extreme from the
smooth regime ($\be_i=1$). In the next few subsections
we discuss some of its interesting properties.

\subsection{Warm-up: classification of del Pezzo surfaces}
\label{WarmUpSubSec}

Which compact complex surfaces admit a \K metric
of positive Ricci curvature?
The following basic classification result, together
with the Calabi--Yau theorem gives a complete list.

\begin{thm}
\label{delPezzoThm}
Let $S$ be a compact complex surface. Then $c_1(S)$ is ample
if and only if $S$ is either $\PP^1\times\PP^1$, or otherwise $\PP^2$ blown-up
at at most 8 distinct points, of which no three are collinear, no
six lie on a conic, and no eight lie on a cubic with one
of the points being a double point.
\end{thm}

Del Pezzo first described some of the eponymous surfaces in the ninteenth century
\cite{DelPezzo}; more precisely, he described the ones with
up to six blown-up points, or, in his language, surfaces of
degree $d=c_1^2$ embedded in $\PP^d$. In other words,
the ones for which $-K_M$ is very ample, i.e., for which
the linear series $|-K_M|$ gives a projective embedding;
indeed, since $-K_M>0$,
by Riemann--Roch
$\dim H^0(M,\calO_M(-K_M))
=
\chi(-K_M)
=
\chi(\calO_M)+c_1^2
=1+c_1^2.
$
What are now known as del Pezzo surfaces are the surfaces
for which $-K_M$ is ample, i.e., those in the statement
of Theorem \ref{delPezzoThm}.
By the Kodaira Embedding Theorem, those are the surfaces
for which $|-mK_M|$ gives a projective embedding for some $m\in\NN$.
It is hard to trace precisely the original discoverers of those
remaining del Pezzo surfaces, let alone the first time
Theorem \ref{delPezzoThm} was stated in this form in the literature, but the contributions
of Clebsch, Segre, Enriques, Nagata, among others, played a crucial role. 
We refer to \cite{ConfortoEnriques,Dolgachev,CalabriCiliberto,
BabbitGoodstein,CilibertoSallent} for more references and historical notes.

This result is used, and generalized, in Theorem \ref{theorem:classification}.
We describe a proof, closely following Hitchin \cite{Hitchin1975}
(see also, e.g., \cite{Yau1974,Dolgachev,Friedman}). 
The detailed proof serves to motivate later classification results,
as well as to establish notation and basic results that are
useful later.

\begin{proof}
To start, one checks that indeed the surfaces in the statement
are del Pezzo. We concentrate on the converse.

{\it Step 1.}
Let $\O\in H^2(M,\RR)\cap H^{1,1}_{\dbar}$.
By Nakai's criterion \cite{Nakai,Buchdahl,Lazarsfeld} 
\beq\label{NMEq}
\h{$\O>0$ if and only if
$\O^2>0$ and $\O.C>0$ for every curve $C$ in $S$.}
\eeq
Thus, any blow-down
$\pi:\tilde S\ra S$ of $\tilde S$ with $c_1(\tilde S)>0$
will also satisfy $c_1(S)>0$:
indeed, if $E$ is the exceptional divisor of $\pi$
then $E^2=\deg L_E|_E=-1$ and \cite[p. 185,187]{GH}
\beq\label{KSEq}
K_{\tilde S}=\pi^\star K_{ S}+E.
\eeq
Thus, $c^2_1(S)=c^2_1(\tilde S)+1>0$.
Additionally, if $\Sigma$ is any holomorphic curve
in $S$ then the associated cohomology
class $[\Sigma]$ (represented by the current of integration along it),
that by abuse of notation we still denote by $\Sigma$,
satisfies $\tilde \Sigma+mE=\pi^\star\Sigma$, where
$\tilde \Sigma:=\overline{\pi^{-1}(\Sigma\sm E)}$ denotes
the proper transform of $\Sigma$, and where $m$
is the multiplicity of $\tilde \Sigma$ at the blow-up point.
Since cup product is preserved under birational transformations,
$
K_S.\Sigma=\pi^\star K_S.\pi^\star\Sigma
=
(K_{\tilde S}-E).\pi^\star\Sigma
=
K_{\tilde S}.(\tilde\Sigma+mE)\le 0$. Here,
we used the fact that any pulled-back class
has zero intersection number with the exceptional
divisor.

It thus remains to classify all del Pezzos with no $-1$-curves,
since by a classical theorem of Castelnuovo--Enriques \cite[p. 476]{GH}
there always exists such a birational blow-down $\pi$ contracting any
given $-1$-curve. The goal is to show that these ``minimal
del Pezzos" are precisely $\PP^2$ and $\PP^1\times\PP^1$.
To that end, one first observes that $S$ is rational,
i.e., birational to $\PP^2$. Indeed, by the Kodaira
Vanishing Theorem \cite[p. 154]{GH}, $H^k(S,\calO(mK_S))=0, k=0,1, m\in\NN$.
In particular, $H^0(S,\calO_S(2K_S))=0$, and
by Dolbeault's theorem and
Kodaira--Serre duality
$H^{0,1}_{\dbar}\cong H^1(S,\calO_S)\cong H^1(S,\calO_S(K_S))$.
Consequently, by the Castelnuovo--Enriques characterization
$S$ is rational \cite[p. 536]{GH}. The classification of
minimal rational surfaces implies that these are
precisely $\PP^2$ and the Hirzebruch surfaces,
\beq\label{HirzEq}
\FF_m:=\PP(\calO\oplus\calO(m)),
\eeq
with $m\in\NN_0\sm\{1\}$,
the projectivization
of the rank 2 bundle over $\PP^1$ that is obtained
by the direct sum of the trivial bundle and the degree
$m$ bundle. (Note that $\FF_0=\PP^1\times\PP^1$,
while $\FF_1$ is the blow-up of $\PP^2$ at one point, hence
is not minimal.)
Finally, observe that $\FF_m, m\ge2$ are not del Pezzo:
indeed they contain a rational holomorphic curve of self-intersection $-m$,
while by adjunction, any curve $C$ on a del Pezzo surface $S$ satisfies
$C^2=2g_C-2-K_S.C>-2$.

{\it Step 2.} Next, we classify the admissible blow-ups
of $\PP^2$ and $\PP^1\times\PP^1$. Since the two-point
blow-up of the former equals the one-point blow-up of
the latter, we may concentrate on $\PP^2$.
Since $c_1^2(\PP^2)=9$,
by \eqref{KSEq} and the relation following it,
at most 8 blow-ups are allowed, according to \eqref{NMEq}.
Next, it remains to determine
the allowable configurations of blow-ups. Suppose that
$k\le 8$ points have been blown up and that the resulting
surface is not del Pezzo.  Then by \eqref{NMEq} this means
that there exists a
curve $C\subset S$
with $C.K_S\ge 0$. We may assume $C$ is irreducible,
since at least one of its components will have nonnegative
intersection number with $K_S$.
Denote by $\pi:S\ra\PP^2$
the blow-down map, by $E\subset S$ the exceptional divisor,
and by $\{p_1,\ldots,p_k\}\subset \PP^2 $ the blown-up
points. We denote by $q_1,\ldots,q_l$ the singular
points of $C$, and let $p_{k+i}:=\pi(q_i),\, i=1,\ldots,l$.
 We also let $\Sigma=\pi(C)$ and write $E=\sum_{i=1}^kE_i$,
with $\pi^{-1}(p_i)=E_i$.
Since $\pi$ is an isomorphism outside $E$, $\Sigma$
will have multiplicity one everywhere except, possibly,
at the points $\{p_i\}\cap\Sigma$, and we denote each of these multiplicities
by $m_i$.
Denote by $d$ the degree of $\Sigma$.
By the genus formula for planar irreducible
projective curves \cite[p. 220, 505]{GH}
the genus of $\Sigma$ equals $g_\Sigma = (d-1)(d-2)/2$
precisely when $\Sigma$ is smooth, and in general
\beq\label{GenusInEq}
2g_\Sigma\le (d-1)(d-2)-\sum_{i=1}^{k+l}m_i(m_i-1)
\le (d-1)(d-2)-\sum_{i=1}^km_i(m_i-1).
\eeq
Here the genus of a possibly singular irreducible curve is
defined either as $\dim H^1(\Sigma,\calO_\Sigma)$ \cite[p. 494]{GH},
or as the genus of its unique desingularization \cite[p. 500]{GH};
in particular,
it is nonnegative.
Letting $m:=\sum_{i=1}^km_i$, convexity of $f(x)=x^2$
gives $\sum m_i^2/k\ge (m/k)^2$. Combined with $g_\Sigma\ge0$
this yields $m^2/k-m\le (d-1)(d-2)$, hence
\beq\label{mEq}
2m\le k+\sqrt{k^2+4k(d-1)(d-2)}.
\eeq
Additionally, $C=\pi^\star\Sigma-\sum_{i=1}^km_iE_i$,
and so by \eqref{KSEq},
\beq\label{miEq}
0\le K_S.C=(-\pi^\star 3H+E).(\pi^\star\Sigma-\sum_{i=1}^km_iE_i)
=-3d+m,
\eeq
since $\Sigma\in|d H|$, where $H\ra\PP^2$ is the
hyperplane bundle. Here we are implicitly assuming
that the blown-up points are distinct, or in other words
that we have not blown up any point on the exceptional divisor
of a previous blow-up. This is justified by the observation
that had we blown-up a point on a $-1$-curve, we would
obtain a $-2$-curve, contradicting the Fano assumption (recall
the end of the previous paragraph).

Thus,
$
3d\le m.
$
Plugging this back into \eqref{mEq} and expanding
the resulting inequality yields
\beq
\label{kdInEq}
(9-k)d^2\le 2k.
\eeq
It follows that $(k,d)\in\{(3,1),(4,1),(5,1),(6,1),(6,2),(7,1),
(7,2),(8,1),(8,2),(8,3),(8,4)\}$.
When $d\in\{1,2\}$ by \eqref{GenusInEq} all $m_i$ must equal 1,
thus equality holds in \eqref{GenusInEq}. Thus, the cases
$\{(k,1)\,:\, k=3,\ldots,8\}$ correspond to $\Sigma$
being a line passing through three or more of the $\{p_i\}$,
and the cases $\{(k,2)\,:\, k=6,\ldots,8\}$ correspond to $\Sigma$
being a smooth conic passing through six or more of the points.
If $d=3$ then by \eqref{miEq} at least one $m_i$
equals 2, thus by \eqref{GenusInEq} exactly one such $m_i$ exists.
Thus, $(k,d)=(8,3)$ and $\Sigma$ is a singular rational cubic with
a double point passing through one of the blow-up points.
Finally, the case $(8,4)$ is the only one that is excluded,
since it forces $m\ge 12$, while the equation following \eqref{GenusInEq},
namely, $m^2/8-m\le 6$ implies $m\le 12$, thus $m=12$. Now equality
in the latter means precisely that equality also occurs in
$\sum m_i^2/k\ge (m/k)^2$. By strict convexity of $f(x)=x^2$
this means that all the $m_i$ are equal; but since they are also
all at least 2, this implies $m\ge 16$, a contradiction.
\end{proof}

\subsection{Log Fano manifolds}
\label{LogFanoSubsec}

The definition of log Fano manifolds goes back to work of Maeda \cite{Maeda}.
\begin{definition}
\label{definition:log-del-Pezzo-pair} We say that the
pair $(M,D=\sum D_i)$ is \emph{ log Fano} if
$-K_M-D$ is ample.
\end{definition}
\noindent

In dimension 2, these are also called
log del Pezzo surfaces
(to avoid confusion, we remark that some
authors use this terminology to refer to rather different
objects).
The motivation for the adjective
``logarithmic", according to Maeda, is from the work of
Iitaka on classification of open algebraic varieties
where logarithmic differential forms are used
to define invariants of the pair. The open variety associated
to $(M,D)$ is the Zariski open set $M\sm D$.

Maeda posed the following problem.

\begin{problem}\label{MaedaProblem}
Classify log Fano manifolds.
\end{problem}

This problem has a beautiful inductive structure.
Indeed, by the adjunction formula \cite[p. 147]{GH}, any component $D_i$
of $D$, or more precisely the pair $(D_i,\sum_{j\ne i}D_i\cap D_j)$,
is itself a log Fano manifold of one dimension
lower, to wit
$$
K_{D_i}+\sum_{j\ne i}D_j|_{D_i}=(K_M+D)|_{D_i}.
$$
When $n=1$, log Fanos consist precisely of $(\PP^1,\{\h{point}\})$
(we always omit the case of empty boundary,
that in this dimension corresponds to $(\PP^1, \emptyset)$).
Thus, the first step in Problem
\ref{MaedaProblem} should be a classification for $n=2$.
This was provided by Maeda \cite[\S3.4]{Maeda}.

\begin{thm}\label{MaedaThm}
Log del Pezzo surfaces $(S,C)$ are classified as follows:
\hfill\break
(i) $S\cong\mathbb{P}^2$, and $C$ is a line in $S$,%
\hfill\break
(ii)
$S\cong\mathbb{P}^2$, and $C=C_1+C_2$, where each $C_i$ is a line in $S$.%
\hfill\break
(iii)
$S\cong\mathbb{P}^2$, and $C$ is a smooth conic in $S$.%
\hfill\break
(iv)
$S\cong\mathbb{F}_{n}$ for any $n\geqslant 0$, and $C$ is a smooth  rational curve in $S$ such that $C^2=-n$ (such curve is unique if $n\geqslant 1$).%
\hfill\break
(v)
$S\cong\mathbb{F}_n$ for any $n\geqslant 0$, and $C=C_1+C_2$
where $C_1$ is as in (iv) and $C_2$ is a smooth fiber (i.e.,
a smooth rational curve such that $C_2^2=0, C_2.C_1=1$).
\hfill\break
(vi)
$S\cong\mathbb{F}_1$, and $C$ is a smooth  rational curve such that
$C\in |C_1+C_2|$, with $C_1,C_2$ as in (v).
\hfill\break
(vii)
$S\cong\mathbb{P}^1\times\mathbb{P}^1$, and $C$ is a smooth rational
curve in $|H_1+H_2|$ where $H_1,H_2$ are lines in each copy of $\mathbb{P}^1$.

\end{thm}

Building on this result and considerable more work, Maeda then
tackles the case $n=3$. Much more recently, Fujita was able
to obtain some results in higher dimensions, especially
for pairs with high log Fano index \cite{Fujita}.

\subsection{Asymptotically log Fano manifolds}

In this section we finally get to
the case $0<\be_i\ll 1$, that generalizes both
extremal cases $\be_i=1$ and $\be_i=0$ studied
in the last two sections.

\begin{definition}
\label{definition:log-Fano} We say that a pair $(M,D)$
is {\it (strongly) asymptotically log Fano} if
the divisor $-K_M^\be=-K_M-\sum_{i=1}^r(1-\beta_i)D_i$ is ample for (all)
sufficiently small $(\beta_1,\ldots,\beta_r)\in(0,1]^r$.
\end{definition}

Both of these classes (strongly asymptotically log Fano and asymptoticaly log Fano)
generalize the class of log Fanos since ampleness (of $-K_X-D$) is an open property.
They also generalize the class of Fanos (at least in small dimensions, see the discussion
surrounding Problem \ref{FanoProb}), since if  $D$ is a smooth anticanonical divisor
in a Fano $M$, then $(M,D)$ is strongly asymptotically log Fano.

The notion of
strongly asymptotically log Fano coincides with that of asymptoticaly log Fano
in the case $r=1$, i.e.,
when $D$ consists of a single smooth component.
However, they differ in general, as the following example
demonstrates.

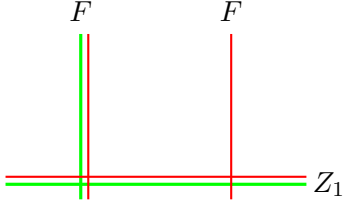
\begin{figure}
	\centering
	%\begin{subfigure}{0.4\linewidth}
	\centering
		\begin{tikzpicture}
			\draw [green, very thick] (-2,0) -- (2,0);
			\draw [green, very thick] (-1,-.2) -- (-1,2);
			\draw [red, thick] (1,-.2) -- (1,2);
			\draw [red, thick] (-2,.1) -- (2,.1);
			\draw [red, thick] (-0.9,-.2) -- (-0.9,2);
			\node at (-1,2.3) {$F$};
			\node at (1,2.3) {$F$};
			\node at (2.3,0) {$Z_{1}$};
		\end{tikzpicture}
		\caption{A non-strongly asymptotically log del Pezzo pair
(see Example \ref{AsympExample}): $\FF_1$
with boundary consisting of two fibers and the $-1$-curve
(the fibration structure of the
surface $\FF_1$ is indicated in green, the boundary $C$ (consisting
of those three curves) in red).}
	%\end{subfigure}%
\label{Figure1}
\end{figure}

\begin{example}\label{AsympExample}
{\rm
(See Figure \ref{Figure1}.)
Let $S=\mathbb{F}_{1}$ (recall \eqref{HirzEq}) and $C=C_1+C_2+C_3$ where
$C_1,C_2\in |F|$ are both fibers and $C_3=Z$ is
the $-1$-curve. Note that $-K_S=2Z+3F$. Then
$-K_S^\be.Z=\be_1+\be_2-\be_3$ and so $(S,C)$
is not strongly asymptotically log del Pezzo.
However, one may verify that it is asymptotically log del Pezzo:
when $\be_1+\be_2>\be_3$, the class $-K_S-(1-\be_1)C_1-(1-\be_2)C_2-(1-\be_3)C_3$
is positive. Note this pair is the
blow-up of the pair $\PP^2$ with two lines at their intersection (the pair
II.1B in Figure \ref{Figure3} below).
}
\end{example}

We pose the following problem.

\begin{problem}\label{ClassifyProblem}
	Classify the (strongly) asymptotically log del Pezzo surfaces
and Fano 3-folds.

\end{problem}

In \S\ref{SALDPSubSec} below we explain the solution to this problem
in the case of a smooth divisor \cite{CR}, where
Problem \ref{ClassifyProblem} is solved more
generally for strongly asymptotically log del Pezzos. This
generalizes Maeda's result and the classical classification
of del Pezzo surfaces (Theorems \ref{MaedaThm} and \ref{delPezzoThm}).
The general (i.e., not necessarily strongly) case of surfaces, as well as the three dimensional case
are open and challenging.

\subsubsection{Comparison between the asymptotic and classical
logarithmic regimes}

We end this subsection with a few comparisons between the log Fano
and asymptotically log Fano regimes, emphasizing the
flexibility in the asymptotic classes as opposed to the
rigidity of the class of log Fanos.

One point of similarity between both classes is that unlike Fano
manifolds, for neither logarithmic classes is the degree
of the logarithmic anticanonical bundle bounded uniformly for
fixed dimension. E.g., when $n=2$, $K_M^2\le 9$ for
del Pezzos, while following the notation of
Example \ref{AsympExample}, $(\FF_m,Z)$
(see Theorem \ref{MaedaThm} (iv)) satisfies
$(-K_M-Z)^2=((m+2)F+Z)^2=m+4$.
Thus, already in the log Fano class there are infinitely
many non-diffeomorphic pairs. Another property shared by both classes is that if
$M$ is Fano then in fact
$-K_M-\sum_{i=1}^r(1-\beta_i)D_i$ is ample for {\it all}
$(\beta_1,\ldots,\beta_r)\in(0,1]^r$.

Aside from these properties though, these classes are quite different.

First, the asymptotic notion is no longer inductive, in the
sense that $(D_i,D_i\cap\cup_{j\ne i}D_j)$,
is not necessarily itself asymptotically log Fano.
In fact,
$$
K_{D_i}+\sum_{j\ne i}(1-\be_j)D_j|_{D_i}
=K^\be_M|_{D_i}+\be_iD_i|_{D_i},
$$
and the right hand side may fail to be negative.
Perhaps the simplest example is the pair $(\PP^2,\h{smooth cubic curve})$,
where the boundary is an elliptic curve, hence not Fano.
Thus, in every dimension one may encounter boundaries that were
absent from the classification in lower dimensions.

Second, while $D$ is always connected in the classical setting
\cite[Lemma 2.4]{Maeda},
this is certainly not
so in the asymptotic regime.
As an example, consider $M=\PP^1\times\PP^1$ and $D=D_1+D_2$
a union of two disjoint lines in the same linear series,
say, $D_i\in|H_1|$. However, there is an upper bound on the
number of disjoint components $D$ may have.
The reason $D$ is always connected in the non-asymptotic regime
is the standard logarithmic short exact sequence
$$
0\ra \calO_M(-D)\ra \calO_M\ra \calO_D\ra 0.
$$
Note that $H^0(M,\calO_M(-D))=\{0\}$ since holomorphic functions
on $M$ vanishing on $D$ must be identically zero, as $M$
is connected. Also, $H^1(M,\calO_M(-D))=\{0\}$ since
by Serre duality this vector space is isomorphic to
$H^1(M,\calO_M(K_M+D))=\{0\}$, by Kodaira Vanishing.
Therefore, $H^0(M,\calO_M)\cong H^0(D,\calO_D)$, and thus
the connectivity of $D$ is `inherited' from that of $M$.

Other more refined connectivity properties are also
interesting to compare.
According to Maeda (op. cit.), when $(M,D)$ is log Fano,
$D$ is always ``strongly connected," meaning that
any two components of $D$ intersect. This follows
immediately from the inductive structure already
mentioned. Indeed, this certainly holds for $n=1$.
Suppose now that $D_1$ intersects both $D_2$
and $D_3$. Then since $(D_1,\sum_{j=2}^rD_1\cap D_j)$,
is itself log Fano, then by induction
$(D_2\cap D_1)\cap(D_3\cap D_1)\not=\emptyset$,
therefore also $D_2\cap D_3\not=\emptyset$, as desired.
Such strong connectivity again fails in the asymptotic
world. In fact, Example \ref{AsympExample},
or even simpler, the disconnected
example or the previous paragraph, or even $(\PP^1,\h{2 distinct points})$
(which are both strongly asymptotically log Fano) provide instances of that.

Moreover, the number of components in the boundary of log Fanos
is bounded from above by the dimension (op. cit.): by strong connectivity
any two components intersect and thus all components have
a common point. But the components cross normally!
There is no analogue for this property in the asymptotic regime.
As we will see, the number of boundary components can be arbitrary.

Finally, the class of asymptotically log Fano manifolds seems like
a more natural generalization of the class of Fano manifolds than
the class of log Fanos.
Indeed the latter do contain the Fano manifolds as a subclass
if one allows the case of empty boundary. However
the class of Fanos actually can be considered as a subset
of the asymptotically log Fanos, if one considers pairs $(M,D)$ with
$M$ Fano and $D\in|-K_M|$ a snc divisor, when such a divisor
exists. As an aside, we mention that this last existence
problem is known to hold for all smooth Fano up to dimension three (then even
a smooth anticanonical divisor exists by the
classification
and a theorem Shokurov \cite{Fano,Iskovskih,Shokurov,Shokurov2,MukaiUmeumura,MoriMukai,
IskovskihProkhorov}),
and in general it falls under the world of the Elephant
conjectures going back to Iskovskih \cite{Shokurov60}. In fact, even the
existence
of such a divisor is open, although examples show that, in general, 
one needs to allow for worse singularities than snc \cite{HoringVoisin}.

\begin{problem}\label{FanoProb}
Determine whether an anticanonical divisor
exists on a smooth Fano manifold, and whether it has
some regularity, at least in sufficiently low dimensions.
\end{problem}

One approach to this problem that does not seem to have been tried so far
would be to use Geometric Measure Theory. Indeed, any (holomorphic)
divisor is automatically a minimal submanifold, in fact area
minimizing in its homology class by Wirtinger's inequality
\cite[\S5.4.19]{Federer}. In other words, the \K form provides
for a calibration in the sense of Harvey--Lawson \cite{HarveyLawson1984}.
The question is
then whether an area minimizing representative of the homology class
$[-K_M]$ can be found that is also a complex subvariety,
and if so whether it has some regularity beyond that provided
by general results of GMT.
In view of \cite{MicallefWolfson} this seems to be a delicate question.
 In the real setting, a
famous result says that hypersurfaces can have singularities
only in codimension 7 or higher \cite[Theorem 37.7]{SimonBook}.
Perhaps one approach to Problem \ref{FanoProb} would be
to develop a regularity theory for complex hypersurfaces.
The rigidity of the holomorphic setting might just be
enough for such a theory, which in the general real codimension
greater than one setting breaks down, of course
aside from Almgren's fundamental result saying that singularities
then occur in real codimension two or higher \cite{Almgren2000}.

\subsection{Classification of strongly asymptotically log del Pezzo surfaces}
\label{SALDPSubSec}

The following result gives a complete classification of
strongly asymptotically log del Pezzo surfaces with
smooth connected boundary.

\begin{thm}
\label{theorem:classification} Let $S$ be a smooth surface (the surface),
and let $C$ be an irreducible smooth curve on $S$ (the boundary
curve). Then $-K_S-(1-\beta)C$ is ample for all sufficiently
small $\beta>0$ if and only if $S$ and $C$ can be described as
follows:
\begin{itemize}\setlength{\itemsep}{-1\parsep}
\item [$\mathrm{(I.1A})$] $S\cong\mathbb{P}^2$, and $C$ is a smooth cubic elliptic curve,%
\item [$\mathrm{(I.1B})$] $S\cong\mathbb{P}^2$, and $C$ is a smooth conic,%
\item [$\mathrm{(I.1C})$] $S\cong\mathbb{P}^2$, and $C$ is a line,%
\item [$\mathrm{(I.2.n})$] $S\cong\mathbb{F}_{n}$ for any $n\ge 0$, and $C=Z_n$,%
\item [$\mathrm{(I.3A})$] $S\cong\mathbb{F}_1$, and $C\in |2(Z_1+F)|$,%
\item [$\mathrm{(I.3B})$] $S\cong\mathbb{F}_1$, and $C\in |Z_1+F|$,%
\item [$\mathrm{(I.4A})$] $S\cong\mathbb{P}^1\times\mathbb{P}^1$, and $C$ is a smooth elliptic curve of bi-degree $(2,2)$,%
\item [$\mathrm{(I.4B})$] $S\cong\mathbb{P}^1\times\mathbb{P}^1$, and $C$ is a smooth  rational curve of bi-degree $(2,1)$,%
\item [$\mathrm{(I.4C})$] $S\cong\mathbb{P}^1\times\mathbb{P}^1$, and $C$ is a smooth  rational curve of bi-degree $(1,1)$,%
\item [$\mathrm{(I.5.m})$] $S$ is a blow-up of the surface in $\mathrm{(I.1A)}$ at $m\le 8$ distinct points on the boundary curve such that $-K_{S}$ is ample, i.e., $S$ is a del Pezzo surface, and $C$ is the proper transform of the boundary curve in $\mathrm{(I.1A)}$, i.e., $C\in|-K_{S}|$,%
\item [$\mathrm{(I.6B.m})$] $S$ is a blow-up of the surface in $\mathrm{(I.1B)}$ at $m\ge 1$ distinct points on the boundary curve, and $C$ is the proper transform of the boundary curve in $\mathrm{(I.1B)}$,%
\item [$\mathrm{(I.6C.m})$] $S$ is a blow-up of the surface in $\mathrm{(I.1C)}$ at $m\ge 1$ distinct points on the boundary curve, and $C$ is the proper transform of the boundary curve in $\mathrm{(I.1C)}$,%
\item [$\mathrm{(I.7.n.m})$] $S$ is a blow-up of the surface in $\mathrm{(I.2.n)}$ at $m\ge 1$ distinct points on the boundary curve, and $C$ is the proper transform of the boundary curve in $\mathrm{(I.2)}$,%
\item [$\mathrm{(I.8B.m})$] $S$ is a blow-up of the surface in
$\mathrm{(I.3B)}$ at $m\ge 1$ distinct points on the boundary curve, and $C$ is the proper transform of the boundary curve in $\mathrm{(I.3B)}$,%
\item [$\mathrm{(I.9B.m})$] $S$ is a blow-up of the surface in $\mathrm{(I.4B)}$ at $m\ge 1$ distinct points on the boundary curve with no two
of them on a single curve of bi-degree $(0,1)$, and $C$ is the proper transform of the boundary curve in $\mathrm{(I.4B)}$,%
\item [$\mathrm{(I.9C.m})$] $S$ is a blow-up of the surface in $\mathrm{(I.4C)}$ at $m\ge 1$ distinct points on the boundary curve, and $C$ is the proper transform of the boundary curve in $\mathrm{(I.4C)}$.%
\end{itemize}
\end{thm}

\begin{figure}
	\centering
	\begin{subfigure}{0.2\linewidth}
	\centering
	\begin{tikzpicture}[scale=.3]
		\draw [green, very thick] (0,0) circle (3);
		\draw [red,thick] plot [smooth, tension=1.2] coordinates {(2,0) (0,-1) (1,1) (-1,.5) } ;
	\end{tikzpicture}
	\caption{I.1A}
	\end{subfigure}%
	\begin{subfigure}{0.2\linewidth}
	\centering
	\begin{tikzpicture}[scale=.3]
		\draw [green, very thick] (0,0) circle (3);
		\draw [red,thick] plot [smooth,tension=1.3] coordinates {(-.7,-1.5) (1.3,.5) (-1.4,1) } ;
	\end{tikzpicture}
	\caption{I.1B}
	\end{subfigure}%
	\begin{subfigure}{0.2\linewidth}
	\centering
	\begin{tikzpicture}[scale=.3]
		\draw [green, very thick] (0,0) circle (3);
		\draw [red,thick] plot [smooth, tension=1.2] coordinates {(-1.2,-1) (.2,-.4) (1.6,.5) } ;
	\end{tikzpicture}
	\caption{I.1C}
	\end{subfigure}%
%\end{figure}

%\begin{figure}[t]
	\centering
	\begin{subfigure}{0.2\linewidth}
	\centering
	\begin{tikzpicture}[scale=.4]
			\draw [green, dashed] (-2,0) -- (2,0);
			\draw [green, dashed] (0,-.2) -- (0,2);
			\draw [red, thick] (-2,.1) -- (2,.1);
			\node at (0,2.3) {$F$};
			\node at (2.3,0) {$Z_{n}$};
		\end{tikzpicture}
		\caption{I.2.n}
	\end{subfigure}%
	\begin{subfigure}{0.2\linewidth}
	\centering
		\begin{tikzpicture}[scale=.4]
			\draw [green, dashed] (-2,0) -- (2,0);
			\draw [green, dashed] (0,-.2) -- (0,2);
			\draw [red, thick] plot [smooth, tension =2] coordinates {(-.5,.5) (1,1) (-.5,1.5)} ;
			\node at (0,2.3) {$F$};
			\node at (2.3,0) {$Z_{1}$};
		\end{tikzpicture}
		\caption{I.3A}
	\end{subfigure}%
	\begin{subfigure}{0.2\linewidth}
	\centering
		\begin{tikzpicture}[scale=.4]
			\draw [green, dashed] (-2,0) -- (2,0);
			\draw [green, dashed] (0,-.2) -- (0,2);
			\draw [red, thick] (-2,1.5) -- (2,1.5);
			\node at (0,2.3) {$F$};
			\node at (2.3,0) {$Z_{1}$};
		\end{tikzpicture}
		\caption{I.3B}
	\end{subfigure}%
%\end{figure}
%
%\begin{figure}
	\centering
	\begin{subfigure}{.2\linewidth}
	\centering
	\begin{tikzpicture}[scale=.4]
			\draw [green, very thick] (-2,0) rectangle (2,4);
			\draw [green, dashed] (0,0) -- (0,4);
			\draw [green, dashed] (-2,2.1) -- (2,2.1);
			\draw [red, thick] (0,2) circle (1.4);
		\end{tikzpicture}
		\caption{I.4A}
	\end{subfigure}
	\begin{subfigure}{.2\linewidth}
	\centering
		\begin{tikzpicture}[scale=.4]
			\draw [green, very thick] (-2,0) rectangle (2,4);
			\draw [green, dashed] (0,0) -- (0,4);
			\draw [green, dashed] (-2,2.1) -- (2,2.1);
			\draw [red, thick] plot [smooth, tension=2] coordinates {(-1,1) (.2,3.4) (1,1.2)} ;
		\end{tikzpicture}
		\caption{I.4B}
	\end{subfigure}
	\begin{subfigure}{.2\linewidth}
	\centering
		\begin{tikzpicture}[scale=.4]
			\draw [green, very thick] (-2,0) rectangle (2,4);
			\draw [green, dashed] (0,0) -- (0,4);
			\draw [green, dashed] (-2,2.1) -- (2,2.1);
			\draw [red, thick] (-1,.5) -- (1.5,3) ;
		\end{tikzpicture}
		\caption{I.4C}
	\end{subfigure}
\caption{Strongly asymptotically log del Pezzo surfaces with
smooth connected boundary: the minimal pairs (the
surface is indicated in green, the boundary in red,
the fibration structure, when one exists, is indicated by the dashed green lines). 
The remaining pairs
listed in Theorem \ref{theorem:classification} are obtained by blowing-up
along the boundary curves as follows: I.1A as described in Theorem \ref{delPezzoThm};
I.1B, I.1C, I.2.n, I.3B, and I.4C at any number of distinct points; I.4B at any  number
of distinct point with no two on a single $(0,1)$-fiber. Note that I.4A and I.3A
may also be blown-up but these cases are covered by blow-ups of I.1A and I.4B,
respectively.
}
\label{Figure2}
\end{figure}
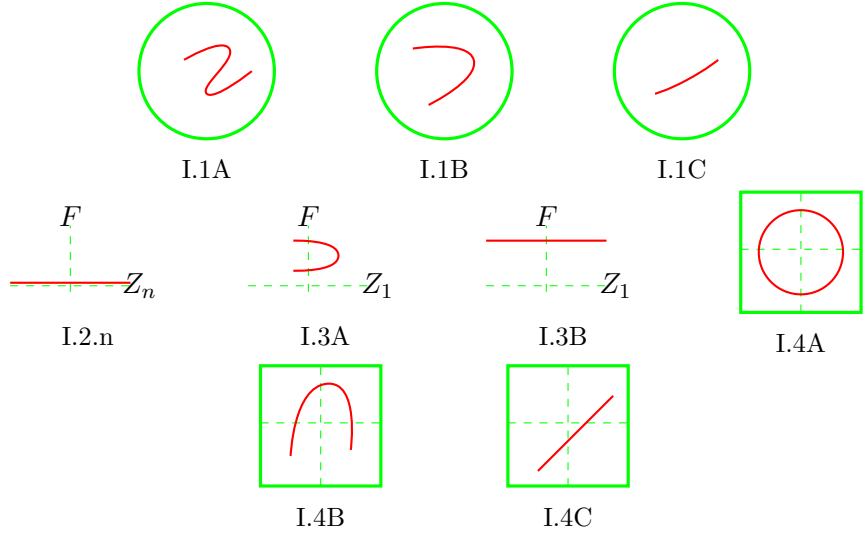

%\vfill\eject

The proof appears in \cite{CR}. We sketch the main
steps. Figure 1 illustrates the pairs graphically.

One starts by checking directly that indeed all the pairs
in the list above are asymptotically log Fano. Let us concentrate on the reverse
implication.

First, it follows from the asymptotic assumption that
$-K_S-C$ is nef. Also, $-K_S$ is big and nef, since it is a linear
combination of an ample class $-K_S-(1-\be)C$ and an effective
class $(1-\be)C$.
The former implies that the genus of $C$ is at most one. The latter,
together with a version of Nadel Vanishing Theorem and a theorem of
Castelnuovo, imply that $S$ is rational.
This further implies that if $C$ is elliptic then $C\in|-K_S|$,
and $S$ is del Pezzo, i.e., is one of (1A),(4A), or $(5_m)$.
On the other hand, if $C$ is rational, then it must ``trap" all
the negative curvature of $-K_S$. More precisely, the only
curve that can intersect $-K_S$ nonpositively is $C$,
and that happens if and only if $C^2\le-2$ (compare
to the end of Step 1 in the proof of Theorem \ref{delPezzoThm}). Thus, all other
negative self-intersection curves must be $-1$-curves. Furthermore,
these curves must be either disjoint from $C$, or intersect it
transversally at exactly one point. This motivates the following definition.

\begin{figure}
	\centering
	\begin{subfigure}{0.2\linewidth}
	\centering
	\begin{tikzpicture}[scale=.4]
		\draw [green, very thick] (0,0) circle (3);
		\draw [red,thick] (-2,-.8)--(2,-.4);
		\draw [red,thick] plot [smooth] coordinates {(-1.4,-1.5) (-1,1) (0, 1.3)  (1,-1)} ;
	\end{tikzpicture}
	\caption{II.1A}
	\end{subfigure}%
	\begin{subfigure}{0.2\linewidth}
	\centering
	\begin{tikzpicture}[scale=.4]
		\draw [green, very thick] (0,0) circle (3);
		\draw [red,thick] (-1.4,-1.4)--(2,.7);
		\draw [red,thick] (-.7,1)--(2,-.7);
	\end{tikzpicture}
	\caption{II.1B}
	\end{subfigure}
	\begin{subfigure}{0.2\linewidth}
	\centering
		\begin{tikzpicture}[scale=.4]
			\draw [green, very thick] (-2,0) -- (2,0);
			\draw [green, very thick] (0,-.2) -- (0,2);
			\draw [red, thick] (-2,.1) -- (2,.1);
			\draw [red, thick] (-2,1.5) -- (2,1.5);
			\node at (0,2.3) {$F$};
			\node at (2.3,0) {$Z_{n}$};
		\end{tikzpicture}
		\caption{II.2A.n}
	\end{subfigure}
	\begin{subfigure}{0.2\linewidth}
	\centering
		\begin{tikzpicture}[scale=.4]
			\draw [green, very thick] (-2,0) -- (2,0);
			\draw [green, very thick] (0,-.2) -- (0,2);
			\draw [red, thick] (-2,.1) -- (2,.1);
			\draw [red,thick] plot [smooth] coordinates {(-1,-0.3) (-0.4,.3) (1.3,1.5)};
			\node at (0,2.3) {$F$};
			\node at (2.3,0) {$Z_{n}$};
		\end{tikzpicture}
		\caption{II.2B.n}
	\end{subfigure}
	\begin{subfigure}{0.2\linewidth}
	\centering
		\begin{tikzpicture}[scale=.4]
			\draw [green, very thick] (-2,0) -- (2,0);
			\draw [green, very thick] (0,-.2) -- (0,2);
			\draw [red, thick] (-2,.1) -- (2,.1);
			\draw [red, thick] (0.1,-.2) -- (0.1,2);
			\node at (0,2.3) {$F$};
			\node at (2.3,0) {$Z_{n}$};
		\end{tikzpicture}
		\caption{II.2C.n}
	\end{subfigure}%
	\begin{subfigure}{.2\linewidth}
	\centering
		\begin{tikzpicture}[scale=.4]
			\draw [green, very thick] (-2,0) -- (2,0);
			\draw [green, very thick] (0,-.2) -- (0,2);
			\draw [red, thick] (-1.5,1.7) -- (2,1);
			\draw [red, thick] (-2,.3) -- (2,1.5);
			\node at (0,2.3) {$F$};
			\node at (2.3,0) {$Z_{n}$};
		\end{tikzpicture}
		\caption{II.3}
	\end{subfigure}
	\begin{subfigure}{.2\linewidth}
	\centering
		\begin{tikzpicture}[scale=.4]
			\draw [green, very thick] (-2,0) rectangle (2,4);
			\draw [green, dashed] (0,0) -- (0,4);
			\draw [green, dashed] (-2,2) -- (2,2);
			\draw [red, thick] (-1,.4) -- (1,2.4);
			\draw [red, thick] (1,.3) -- (-.4,3);
		\end{tikzpicture}
		\caption{II.4A}
	\end{subfigure}
	\begin{subfigure}{.2\linewidth}
	\centering
		\begin{tikzpicture}[scale=.4]
			\draw [green, very thick] (-2,0) rectangle (2,4);
			\draw [green, dashed] (0,0) -- (0,4);
			\draw [green, dashed] (-2,2) -- (2,2);
			\draw [red, thick] (.1,.1) -- (.1,3.9);
			\draw [red,thick] plot [smooth, tension =2] coordinates {(-1,1) (1,2) (-1,3)} ;
		\end{tikzpicture}
		\caption{II.4B}
	\end{subfigure}
%\end{figure}
%
%\begin{figure}
	\centering
	\begin{subfigure}{0.2\linewidth}
	\centering
	\begin{tikzpicture}[scale=.4]
		\draw [green, very thick] (0,0) circle (3);
		\draw [red,thick] (-1.4,-1.4)--(2,.7);
		\draw [red,thick] (-.7,1)--(2,-.7);
		\draw [red,thick] (-.5,-1.4) -- (-.2,1.7);
	\end{tikzpicture}
	\caption{III}
	\end{subfigure}
	\begin{subfigure}{.2\linewidth}
	\centering
		\begin{tikzpicture}[scale=.4]
			\draw [green, very thick] (-2,0) rectangle (2,4);
			\draw [green, dashed] (0,0) -- (0,4);
			\draw [green, dashed] (-2,2.1) -- (2,2.1);
			\draw [red, thick] (.1,0) -- (.1,4);
			\draw [red, thick] (-2,1.5) -- (2,1.5);
			\draw [red, thick] (-1.5,1) -- (1,3.5);
		\end{tikzpicture}
		\caption{III.2A}
	\end{subfigure}
	\begin{subfigure}{0.2\linewidth}
	\centering
		\begin{tikzpicture}[scale=.4]
			\draw [green, very thick] (-2,0) -- (2,0);
			\draw [green, very thick] (0,-.2) -- (0,2);
			\draw [red, thick] (-2,.1) -- (2,.1);
			\draw [red, thick] (-2,1.5) -- (2,1.5);
			\draw [red, thick] (.1,-.2) -- (.1,2);
			\node at (0,2.3) {$F$};
			\node at (2.3,0) {$Z_{n}$};
		\end{tikzpicture}
		\caption{III.3n}
	\end{subfigure}
	\begin{subfigure}{.2\linewidth}
	\centering
		\begin{tikzpicture}[scale=.4]
			\draw [green, very thick] (-2,0) rectangle (2,4);
			\draw [green, dashed] (-.5,0) -- (-.5,4);
			\draw [green, dashed] (-2,1.5) -- (2,1.5);
			\draw [red, thick] (-.4,0) -- (-.4,4);
			\draw [red, thick] (-2,1.6) -- (2,1.6);
			\draw [red, thick] (1,0) -- (1,4);
			\draw [red, thick] (-2,3) -- (2,3);
		\end{tikzpicture}
		\caption{IV}
	\end{subfigure}
\caption{Strongly asymptotically log del Pezzo surfaces with
general snc boundary: the minimal pairs. Essentially, these pairs
are obtained by ``degenerating" the elliptic/rational boundaries of Figure \ref{Figure2}
into cycles/chains of rational curves.
}
\label{Figure3}
\end{figure}
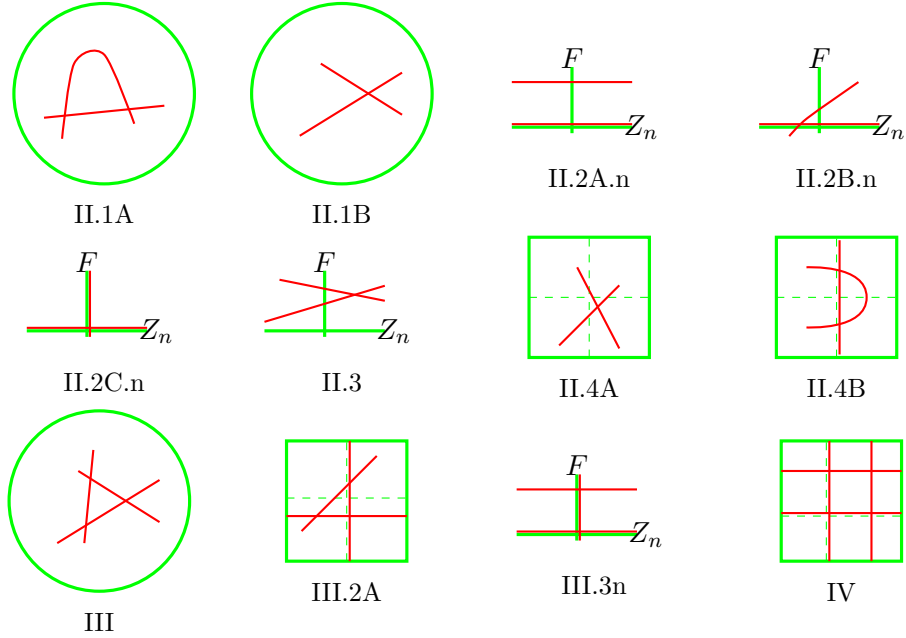

\begin{definition}
\label{definiton:log-del-Pezzo-minimally-ALF} We say that the
pair $(S,C)$ is \emph{minimal }   if there
exist no smooth irreducible rational $-1$-curve $E\ne C$ on the surface
$S$ such that 
$E\cap C\not=\emptyset$.
\end{definition}

The importance of this definition is in the following.

\begin{lem}
\label{corollary:log-del-Pezzo-blow-down} Suppose that $(S,C)$ is
non-minimal asymptotically log del Pezzo and let $E$ be as in Definition \ref{definiton:log-del-Pezzo-minimally-ALF}.
Then there exists a birational morphism $\pi\colon S\to s$
such that $s$ is a smooth surface, $\pi(E)$ is a point, the
morphism $\pi$ induces an isomorphism $S\setminus E\cong
s\setminus\pi(E)$, the curve $\pi(C)$ is smooth, and
$(s, \pi(C))$ is asymptotically log del Pezzo.
\end{lem}

Thus, it remains to classify all minimal pairs. First, one
proves that minimality implies the rank of the Picard group
of $S$ is at most two. Second, one shows, using the classical
theory of rational surfaces,
  that a minimal
pair with this rank restriction must be
$\mathrm{(I.1B)}$, $\mathrm{(I.1C)}$, $\mathrm{(I.2.n)}$,
$\mathrm{(I.3A)}$, $\mathrm{(I.3B)}$, $\mathrm{(I.4B)}$, or
$\mathrm{(I.4C)}$, and a non-minimal one must equal
$\mathrm{(I.6B.1)}$ or $\mathrm{(I.6C.1)}$. This concludes the proof
since  $\mathrm{(I.6B.m),(I.6C.m)}$, with $m\ge2$,
and (I.7.n.m),(I.8B.m),(I.9B.m), (I.9C.m), with $m\ge1$,
are precisely the only blow-ups of minimal pairs that are
still asymptotically log del Pezzo.

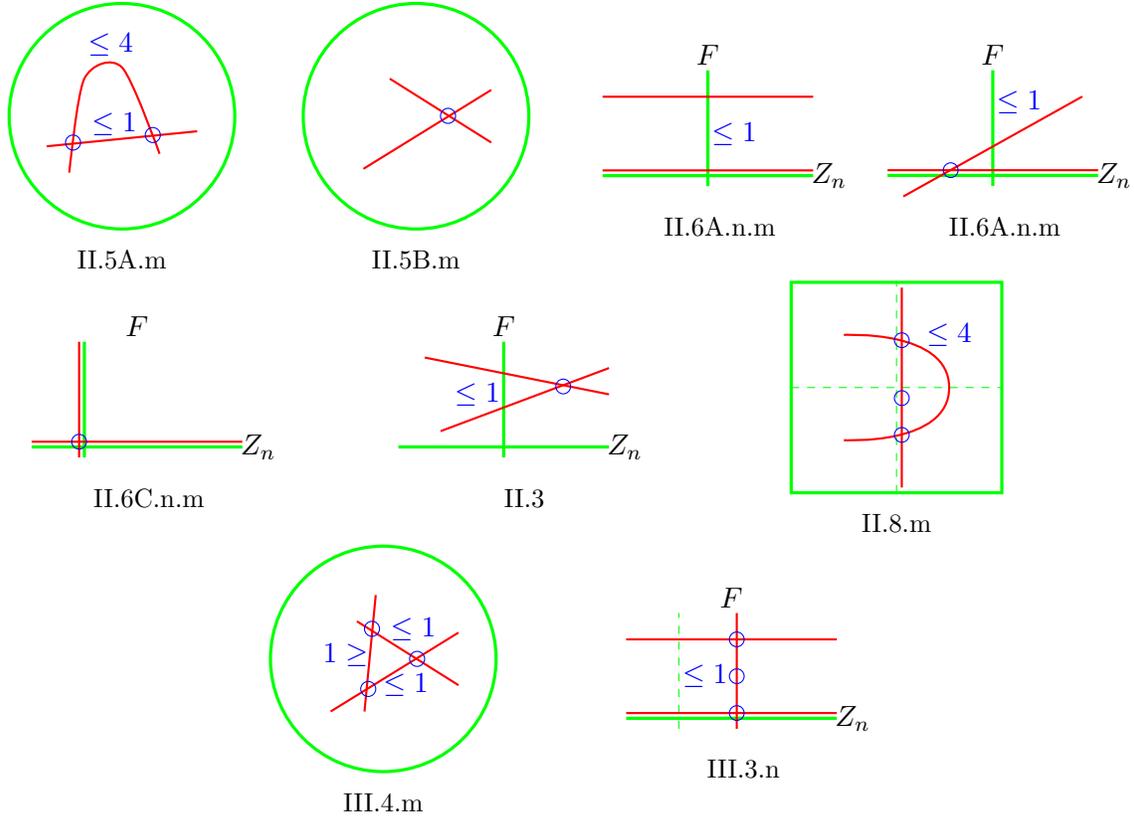
\begin{figure}
	\centering
	\begin{subfigure}{0.2\linewidth}
	\centering
	\begin{tikzpicture}[scale=.5]
		\draw [green, very thick] (0,0) circle (3);
		\draw [red,thick] (-2,-.8)--(2,-.4);
		\draw [red,thick] plot [smooth] coordinates {(-1.4,-1.5) (-1,1) (0, 1.3)  (1,-1)} ;
		\draw [blue] (.82,-.5) circle (.2);
		\draw [blue] (-1.3,-.7) circle (.2);
		\node [blue] at (-.2,-.2) {$\leq 1$};
		\node [blue] at (-.3,1.9) {$\leq 4$};
	\end{tikzpicture}
	\caption{II.5A.m}
	\end{subfigure}%
\hglue0.75cm
	\begin{subfigure}{0.2\linewidth}
	\centering
	\begin{tikzpicture}[scale=.5]
		\draw [green, very thick] (0,0) circle (3);
		\draw [red,thick] (-1.4,-1.4)--(2,.7);
		\draw [red,thick] (-.7,1)--(2,-.7);
		\draw [blue] (.85,.01) circle (.2);
	\end{tikzpicture}
	\caption{II.5B.m}
	\end{subfigure}
\hglue0.75cm
	\begin{subfigure}{0.2\linewidth}
	\centering
		\begin{tikzpicture}[scale=.7]
			\draw [green, very thick] (-2,0) -- (2,0);
			\draw [green, very thick] (0,-.2) -- (0,2);
			\draw [red, thick] (-2,.1) -- (2,.1);
			\draw [red, thick] (-2,1.5) -- (2,1.5);
			\node at (0,2.3) {$F$};
			\node at (2.3,0) {$Z_{n}$};
			\node [blue] at (.5,.8) {$\leq 1$};
		\end{tikzpicture}
		\caption{II.6A.n.m}
	\end{subfigure}
\hglue0.5cm
	\begin{subfigure}{0.2\linewidth}
	\centering
		\begin{tikzpicture}[scale=.7]
			\draw [green, very thick] (-2,0) -- (2,0);
			\draw [green, very thick] (0,-.2) -- (0,2);
			\draw [red, thick] (-2,.1) -- (2,.1);
			\draw [red, thick] (-1.7,-.4) -- (1.7,1.5);
			\node at (0,2.3) {$F$};
			\node at (2.3,0) {$Z_{n}$};
			\node [blue] at (.5,1.4) {$\leq 1$};
			\draw [blue] (-.8,.1) circle (.14);
		\end{tikzpicture}
		\caption{II.6A.n.m}
	\end{subfigure}
\hglue0.5cm
	\begin{subfigure}{0.2\linewidth}
	\centering
		\begin{tikzpicture}[scale=.7]
			\draw [green, very thick] (-2,0) -- (2,0);
			\draw [green, very thick] (-1,-.2) -- (-1,2);
			\draw [red, thick] (-2,.1) -- (2,.1);
			\draw [red, thick] (-1.1,-.2) -- (-1.1,2);
			\node at (0,2.3) {$F$};
			\node at (2.3,0) {$Z_{n}$};
			\draw [blue] (-1.1,.1) circle (.14);
		\end{tikzpicture}
		\caption{II.6C.n.m}
	\end{subfigure}%
\hglue0.25cm
	\begin{subfigure}{.4\linewidth}
	\centering
		\begin{tikzpicture}[scale=.7]
			\draw [green, very thick] (-2,0) -- (2,0);
			\draw [green, very thick] (0,-.2) -- (0,2);
			\draw [red, thick] (-1.5,1.7) -- (2,1);
			\draw [red, thick] (-1.2,.3) -- (2,1.5);
			\node at (0,2.3) {$F$};
			\node at (2.3,0) {$Z_{n}$};
			\draw [blue] (1.13,1.15) circle (.14);
			\node [blue] at (-.5,1) {$\leq 1$};
		\end{tikzpicture}
		\caption{II.3}
	\end{subfigure}
\hglue-1.5cm
	\begin{subfigure}{.4\linewidth}
	\centering
		\begin{tikzpicture}[scale=.7]
			\draw [green, very thick] (-2,0) rectangle (2,4);
			\draw [green, dashed] (0,0) -- (0,4);
			\draw [green, dashed] (-2,2) -- (2,2);
			\draw [red, thick] (.1,.1) -- (.1,3.9);
			\draw [red,thick] plot [smooth, tension =2] coordinates {(-1,1) (1,2) (-1,3)} ;
			\draw [blue] (.1,1.1) circle (.14);
			\draw [blue] (.1,2.9) circle (.14);
			\draw [blue] (.1,1.8) circle (.14);
			\node [blue] at (1,3) {$\leq 4$};
		\end{tikzpicture}
		\caption{II.8.m}
	\end{subfigure}
	\begin{subfigure}{0.2\linewidth}
	\centering
	\begin{tikzpicture}[scale=.5]
		\draw [green, very thick] (0,0) circle (3);
		\draw [red,thick] (-1.4,-1.4)--(2,.7);
		\draw [red,thick] (-.7,1)--(2,-.7);
		\draw [red,thick] (-.5,-1.4) -- (-.2,1.7);
		\draw [blue] (-.3,.8) circle (.2);
		\draw [blue] (-.4,-.8) circle (.2);
		\draw [blue] (.9,0) circle (.2);
		\node [blue] at (.8,.8) {$\leq 1$};
		\node [blue] at (-1,.1) {$1\geq$};
		\node [blue] at (.6,-.7) {$\leq 1$};
	\end{tikzpicture}
	\caption{III.4.m}
	\end{subfigure}
\hglue1.5cm
	\begin{subfigure}{0.2\linewidth}
	\centering
		\begin{tikzpicture}[scale=.7]
			\draw [green, very thick] (-2,0) -- (2,0);
			\draw [green, dashed] (-1,-.2) -- (-1,2);
			\draw [red, thick] (-2,.1) -- (2,.1);
			\draw [red, thick] (-2,1.5) -- (2,1.5);
			\draw [red, thick] (.1,-.2) -- (.1,2);
			\node at (0,2.3) {$F$};
			\node at (2.3,0) {$Z_{n}$};
			\node [blue] at (-.5,.8) {$\leq 1$};
			\draw [blue] (.1,1.5) circle (.14);
			\draw [blue] (.1,.1) circle (.14);
			\draw [blue] (.1,.8) circle (.14);

		\end{tikzpicture}
		\caption{III.3.n}
	\end{subfigure}
\caption{
\label{FigureSAlDPgeneral-snc}
Strongly asymptotically log del Pezzo surfaces with
general snc boundary: the remaining cases obtained by blowing-up
the minimal pairs in Figure \ref{Figure3}. A circle corresponds to a point that may not
be blown-up. An indication ``$\le k$" next to a curve means that
no more than $k$ distinct points may be blown-up on that curve.
In III.3.n no more than one point may be blown up on any single fiber
and none on the fiber belonging to the boundary.
}
\end{figure}

Building on this, the case of a snc boundary is handled in
\cite{CR}. Essentially, aside from a few cases of
disconnected boundary, the only new boundaries allowed
beyond the smooth connected boundary case are boundaries
that can be considered as ``degenerations" of smooth ones.
For instance, the smooth elliptic boundary of $\mathrm{(I.1A)}$
can be replaced by a triangle of lines, or a conic and a line;
the elliptic boundary of $\mathrm{(I.4A)}$ can similarly
break up to no more than 4 components.
However, an additional complication
in the snc case is that a
$-1$-curve in a non-minimal pair could be a component of the boundary.
Luckily, one can show that
such a curve must be at the `tail': it cannot
intersect two boundary components. Thanks to this,
paying attention to the combinatorical structure of the boundary,
the main idea from the proof of Theorem \ref{theorem:classification}
carries over to give a classification of strongly asymptotically log
del Pezzos \cite{CR}. We list these pairs in Figures \ref{Figure3}--\ref{FigureSAlDPgeneral-snc}, and refer to \cite{CR} for their precise construction.

\subsection{The negative case}
\label{LogNegSubSec}

In this subsection we assume that $\mu<0$. We fix
$M$ and seek
necessary and sufficient restrictions on a snc divisor $D\subset M$
in order to be an admissible boundary for all small enough
$\be$.

\begin{definition} 
\label{definition:log-del-Pezzo-good-pair}
 We say that a pair $(M,D)$
is {\it (strongly) asymptotically log general type} if
the divisor $K_M^\be=K_M+\sum_{i=1}^r(1-\beta_i)D_i$ is ample for (all)
sufficiently small $(\beta_1,\ldots,\beta_r)\in(0,1]^r$.
\end{definition}

The following theorem comes close to describing the (strongly) asymptotically
log general type surfaces as a subclass of the class
of log general type minimal surfaces.
The proof we describe is due to Di Cerbo \cite{DiCerbo2012}
(who showed $(ii)\;\Leftrightarrow\;(iv)$), with
slight modifications to include also the asymptotic classes
defined above.

\begin{prop}
\label{NegativeCaseProp}
Let $S$ be projective surface and $C\subset S$ a snc curve such
that $K_S+C$ is big and nef.
Consider the following statements:
\hfill\break
(i) $(S,C)$ is strongly asymptotically log
general type.
\hfill\break
(ii) $K^\be_S>0$ for $\be=\be_1(1,\ldots,1)$ for $0<\be_1\ll 1$.
\hfill\break
(iii) $(S,C)$ is asymptotically log general type.
\hfill\break
(iv) Every rational $-1$-curve not contained in $C$ 
intersects $C$ at least in two points,
and any rational $-2$-curve $F$ satisfies
$F\not\subset S\sm C$. 
Every  rational component $C_i$ of $C$ intersects
$\cup_{j\ne i}C_j$ at least in two points, and, if 
$C_i^2\in\{-1,-2\}$
then at least in three points.
\hfill\break
(v) Every $-1$-curve not contained in $C$ intersects $C$
at least in two points,
and any rational $-2$-curve $F$ satisfies
$F\not\subset S\sm C$. 
Every rational component $C_i$ of $C$ intersects
$\cup_{j\ne i}C_j$ at least in two points.
\hfill\break
Then $(i)\;\Rightarrow\;(ii)\;\Leftrightarrow\;(iv)\Rightarrow\;(iii)
\;\Rightarrow(v)$.

\end{prop}

Note that, similarly to the last section, the assumption
that $K_S+C$ is big and nef is of course a consequence
of $(S,C)$ being asymptotically log general type.

\begin{proof}
Suppose first that $(S,C)$ is (strongly) asymptotically
log canonical.
Let $F$ be a holomorphic rational curve. Then, $K_S^\be.F>0$, i.e.,
\beq\label{IntersectionEq}
C.F>-K_S.F+\sum\be_i C_i.F=2-2g_F+F^2+\sum\be_i C_i.F
=2+F^2+\sum\be_i C_i.F.
\eeq
Thus $(C-F).F\ge2$ which proves the second part of (v) by letting $F=C_i$.
To prove the second part of (iv), suppose now that
we are in the
strong regime. By putting $\be_i=\be_1\ll 1$, then
$(1-\be_1)C.F>2+F^2$, i.e.,
$(1-\be_1)F.(C-F)>2+\be_1F^2$.
If $F^2\in\{-1,-2\}$ this then majorizes
$2-2\be_1$, so
$F.(C-F)>2$ proving the second part of (iv).
Finally, suppose now that $F^2=-2$ (but not necessarily
in the strong regime). Then, $F^2<-2+F.\sum(1-\be_i)C_i$. Thus, if $F\subset S\sm C$
so that $F.C_i=0$ for each $i$, then necessarily $F^2<-2$, contradicting $F^2=-2$.
To prove the remainder of the first part of (iv) and (v),
let $F^2=-1$. If $F\not=C_i$ for each $i$,
then $F.C_i\ge0$, thus \eqref{IntersectionEq} implies
$C.F>1$. In the strong regime the same inequality
gives $C.F>(1-\be_1)^{-1}>1$ without further assumption on $F$.

Suppose now that (iv) hold. As $(K_S+C)^2>0$
also $(K_S^\be)^2>0$ for all small enough $|\be|$.
By Nakai's criterion \eqref{NMEq}, 
it remains
to show that $K_S^\be$ intersects positively with
every irreducible curve in $S$. By taking $|\be|$
sufficiently small, this is certainly the case
for every curve $Z$ such that $(K_S+C).Z>0$.
Thus, suppose that $(K_S+C).Z=0$.
Note that the cup product $Q$ on $H^{1,1}_{\dbar}$
has exactly one positive eigenvalue \cite[p. 126]{GH}.
Thus if $Q(x,x)>0$, $Q(x,y)=0$ for
some $x,y\in H^{1,1}_{\dbar}$ then
$Q(y,y)=Q(x\pm y,x\pm y)-Q(x,x)\mp2Q(x,y)=Q(x+y,x+y)-Q(x,x)$,
so necessarily $Q(y,y)<0$, otherwise $Q(ax+by,ax+by)>0$,
for all $a,b\in\RR$,
and $Q$ would have at least two positive eigenvalues.
Thus, $Z^2<0$. Now, by our assumption on $Z$,
$(C-Z).Z=2-2g_Z$ (here $g_Z$ denotes the genus
of the desingularization of $Z$). 
Therefore, $g_Z\ge1$ implies $C.Z\le Z^2<0$. But then
$
K_S^\be.Z
=
(K_S+C).Z
-
\sum\be_iC_i.Z=-\sum\be_iC_i.Z,
$
which is necessarily positive for certain $\be$ in any neighborhood
of $0\in\RR_+^r\sm\{0\}$ (for instance, if $\be_i=\be_1$).
On the other hand, suppose $g_Z=0$. 
As we just saw, we may suppose that $C.Z\ge 0$; since
$C.Z=2+Z^2$, this implies $Z^2\in\{-1,-2\}$.
If $Z^2=-2$, so $C.Z=0$, then either 
$Z\subset S\setminus C$---but this is precluded by the 
first part of (iv)---or, $Z\not\subset S\setminus C$,
so necessarily $Z=C_i$ for some $i$, but
then 
$C.Z=(Z+\sum_{j\not=i}C_j).Z\ge-2+3\ge1$,
contradicting the second part of (iv). 
If $Z^2=-1$ then $C.Z=1$. Thus, $Z\not\subset S\setminus C$.
If $Z\not=C_i$ for each $i$ we obtain a contradiction
to the first part of (iv). 
If $Z=C_i$ for some $i$, then 
$C.Z=(Z+\sum_{j\not=i}C_j).Z\ge-1+3\ge2$,
a contradiction.
\end{proof}

\subsection{Uniform bounds}

A natural question is whether there exist uniform bounds
on the asymptotic range of $\be$; and if so, what do
they depend on? This was first addressed by
Di Cerbo--Di Cerbo \cite{DiCerbo2013} in the case $\be=\be_1(1,\ldots,1)$,
and this subsection is mostly a review of these
results. As can be expected, the results are more complete
in the negative regime.

\subsubsection{Strongly asymptotically log general type regime}

Perhaps the simplest example of an asymptotically log general type
pair is $(\PP^n,D)$ with $D\in|\calO(n+2)|$. Then for every
$\be$ in the range $(0,\frac1{n+2})$
$K^\be_M$ is still positive. As shown
by Di Cerbo--Di Cerbo \cite{DiCerbo2013}, this is always the case
when restricting to the ray $\be=\be_1(1,\ldots,1)$.

\begin{prop}
\label{DiCerboProp}
Suppose that $(M,D)$ is such that
$K^\be_M>0$ for $\be=\be_1(1,\ldots,1)$ for some $0<\be_1\ll 1$
(recall (\ref{KbetaDefEq})).
Then the same is true for
$0<\be_1<\frac1{n+2}$.

\end{prop}

\bpf
First, recall the following fact:
\beq
\label{step1}
\h{if $C$ is an irreducible curve such that $(K_M+D).C=0$ then $K_M.C>0$}.
\eeq
In fact, for some $t<1$, $(K_M+tD).C>0$, and since $D.C=-K_M.C$
we conclude $0<(K_M+tD).C=(1-t)K_M.C$.

Second, $K_M+tD$ is nef for every $t\in[\frac{n+1}{n+2},1]$. This implies
the Proposition, since then for any $t\in(\frac{n+1}{n+2},1]$,
$K_M+tD$ is a convex combination of an ample divisor and a nef divisor,
hence positive by Kleiman's criterion. We now prove the nefness claim.
It suffices to show that if $C$ is an irreducible curve with $(K_M+D).C>0$
then $(K_M+tD).C\ge0$; indeed, this is already true
by the first paragraph if $(K_M+D).C=0$ (and since $K_M+D$ is nef as
a limit of ample divisors, it is always true that
$(K_M+D).C\ge0$). Now, we decompose $C$ according
to the cone theorem (see, e.g., \cite{DiCerbo2013})
\beq
\label{coneThm}
C\sim_\ZZ \sum_{i=1}^ra_iC_i+F, \qq a_i>0, \q F.K_X\ge0,\q C_i.K_X\in(-n-1,0).
\eeq
Thus, once again using that $K_M+D$ is nef,
$$
\baeq
(K_M+(n+1)(K_M+D)).C
&=
(K_M+(n+1)(K_M+D)).\Big(\sum_{i=1}^ra_iC_i+F\Big)
\cr
&\ge
-(n+1)\sum a_i+(n+1)\sum a_i C_i.(K_M+D)\ge0,
\eaeq
$$
since $C_i.(K_M+D)\ge1$ as otherwise,
by nefness of $K_M+D$, $C_i.(K_M+D)=0$ which
would imply $K_M.C_i=0$ by \eqref{step1}, contradicting
\eqref{coneThm}. Thus,
$K_M+tD$ is nef for every $t\in[\frac{n+1}{n+2},1]$.
\epf

\begin{remark}{\rm
In \cite{DiCerbo2013} it is shown
that \eqref{step1} also implies that
$K_M+(1-\be_1)D$ is ample for some $0<\be_1\ll 1$ in the
case $K_M+D$ is big.
}
\end{remark}

\subsubsection{Log Fano regime}

Theorem \ref{theorem:classification} implies
that an analogue of Proposition \ref{DiCerboProp}
in the asymptotically log Fano regime is false, and
the correct analogue remains to be found.
In the more restrictive log Fano regime (see \S\ref{LogFanoSubsec}),
Di Cerbo--Di Cerbo prove
an interesting first result in this direction \cite{DiCerbo2013},
based on deep results from algebraic geometry.
It is an a priori bound on the asymptotic regime depending
only on the degree of $-K_M-D$ and $n$.

\begin{prop}
Suppose $(M,D)$ is log Fano. Then $-K_M-(1-\be_1)D$ is positive
for every $\be_1\in[0,\be_\max)$ with
$\be_\max$ depending only on $n$ and $(-K_M-D)^n$.
\end{prop}

\section{
The logarithmic Calabi problem
}
\label{logCalabiSec}

For simplicity, in what follows we always suppose the boundary is smooth and connected
and that the dimension is two. We refer to \cite{CR} for more general considerations.

The preceding section sets the stage for the asymptotic logarithmic Calabi problem:

\begin{problem}\label{logCalabiProb}
Determine which (strongly) asymptotically logarithmic Fano manifolds admit KEE metrics
for sufficiently small $\be$.
\end{problem}

In dimension two, the smooth version of Calabi's problem
was solved by Tian in 1990 who showed that
among the list of Theorem \ref{delPezzoThm}, only
$\PP^2$ blown-up at one or two distinct points do not admit KE metrics \cite{Tian1990}.
In light of Theorem \ref{theorem:classification} it
is very natural and tempting to hope for
a counterpart for strongly asymptotically log
del Pezzo surfaces. The formulation conjectured in \cite{CR} is the following:

\begin{conj}
\label{UniformizationLogDelPezzoConj}
Suppose that $(S,C)$ is strongly asymptotically log
del Pezzo with $C$ smooth and irreducible.
Then $S$ admits KEE metrics with
angle $\beta$ along $C$ for
all sufficiently small $\beta$ if and
only if $(K_{S}+C)^2=0$.
\end{conj}

In Tian's solution of the smooth case the vanishing
of the Futaki invariant provided a necessary and sufficient
condition for existence.
More generally, since the work of Hitchin, Kobayashi, and many others,
a standard condition for the existence of canonical metrics
that can be described as zeros of an infinite-dimensional
moment map is some sort of `stability' condition.
How, then,
does Conjecture \ref{UniformizationLogDelPezzoConj}
fit into this scheme?

\subsection{First motivation: positivity classification and Calabi--Yau fibrations}

It turns out to be quite useful to re-classify the pairs
appearing in Theorem \ref{theorem:classification} according
to the positivity of the logarithmic anticanonical bundle
$-K_S-C$. We distinguish between four mutually exclusive classes.
Class $\mathrm{(\aleph)}$: $S$ is del Pezzo and $C\sim-K_S$;
class $\mathrm{(\beth)}$: $C\not\sim-K_S$ and $(K_S+C)^2=0$;
class $\mathrm{(\gimel)}$: $-K_S-C$ is big but not ample;
class $\mathrm{(\daleth)}$: $-K_S-C$ is ample.

\begin{thm}
\label{theorem:4-cases}
The asymptotically log del Pezzo pairs appearing in
Theorem \ref{theorem:classification}
are classified according to the positivity properties
$\mathrm{(\aleph)}$, $\mathrm{(\beth)}$, $\mathrm{(\gimel)}$,
and $\mathrm{(\daleth)}$ as follows:
\begin{itemize}\setlength{\itemsep}{-1\parsep}
\item[$\mathrm{(\aleph)}$] $(S,C)$ is one of
$\mathrm{(I.1A})$, $\mathrm{(I.4A})$, or $\mathrm{(I.5.m})$.

\item[$\mathrm{(\beth)}$] $(S,C)$ is one of
$\mathrm{(I.3A})$, $\mathrm{(I.4B})$, or $\mathrm{(I.9B.m})$.

\item[$\mathrm{(\gimel)}$] $(S,C)$ is one of
$\mathrm{(I.6B.m})$, $\mathrm{(I.6C.m})$, $\mathrm{(I.7.n.m})$,
$\mathrm{(I.8B.m})$, or
$\mathrm{(I.9C.m})$.

\item[$\mathrm{(\daleth)}$] $(S,C)$ is one of
$\mathrm{(I.1B})$, $\mathrm{(I.1C})$, $\mathrm{(I.3B})$,
$\mathrm{(I.2.n})$, or $\mathrm{(I.4C})$.
\end{itemize}
\end{thm}

This list nicely puts the discussion of \S\ref{WarmUpSubSec}--\ref{LogFanoSubsec} in perspective.
Class $(\daleth)$ is Maeda's class of log del Pezzo surfaces
\cite{Maeda},
while class $\mathrm{(\aleph)}$
is the classical class of del Pezzo surfaces
together with the information of a simple normal crossing
anticanonical curve.
The classes $\mathrm{(\beth)}$ and $\mathrm{(\gimel)}$
are new.

The next result is a structure result for surfaces of class $(\beth)$
\cite{CR}. It is slightly stronger than what Kawamata--Shokurov
basepoint freeness would give: there the relevant linear system
giving a morphism is $|-lK_{S}-lC|$, for some $l\in\NN$.

\begin{prop}
\label{proposition:conic-bundle} If
$(S,C)$ is of class $(\beth)$, then the linear system
$|-K_{S}-C|$ is free from base points and
gives a morphism $S\to\mathbb{P}^1$ whose general fiber
is $\mathbb{P}^1$, and every reducible fiber consists of exactly
two components, each a $\PP^1$.
\end{prop}

Thus, these surfaces  are conic bundles,
and the boundary $C$ intersects each generic fiber at
two points, whose fiber complement is a cylinder!
It is therefore tempting to conjecture:

\begin{conj}
\label{LimitConj}
Let $(S,C,\o_\be)$ be KEE pairs of class $(\aleph)$
or $(\beth)$.
Then $(S,C,\o_\be)$ converges in an appropriate sense to a
a generalized KE metric $\o_\infty$ on $S\setminus C$ as $\beta$ tends to zero.
In particular, $\o_\infty$ is a Calabi--Yau metric
in case $(\aleph)$, and a cylinder along each generic fiber
in case $(\beth)$.
\end{conj}

This conjecture is itself a generalization
of a folklore conjecture in K\"ahler geometry
saying that $S\setminus C$ equipped with the Tian--Yau
metric \cite{TY1990} should be a limit of KEE metrics
on $(S,C)$ when $S$ is of class $(\aleph)$
(see, e.g., \cite[p. 9]{Mazzeo1999}, \cite[p. 76]{Donaldson-linear}).

This gives strong motivation for
the `if' part of
Conjecture \ref{UniformizationLogDelPezzoConj} because it suggests
what the small-angle KEE metrics could be considered as a perturbation
of the complete Calabi--Yau metrics on the complement of $C$.
It also motivates the `only if' part: then there is no good limit,
as the limit class is `too' positive, which should morally preclude
the existence of a smooth non-compact complete metric on it
(having Myers' theorem in mind).

\subsection{Second motivation: asymptotic log canonical thresholds}

Perhaps further evidence for Conjecture \ref{UniformizationLogDelPezzoConj}
is given by the following result.

\begin{thm}
\label{theorem:log-del-Pezzo-alpha}
Assume $(S,C)$ is asymptotically log del Pezzo with $C$ smooth and irreducible.
Then
$$
\lim_{\beta\to 0^+}\alpha(S,(1-\beta)C)
=
\left\{\aligned%
&1\qquad\h{class $(\aleph)$},\\
&1/2\quad\h{class $(\beth)$},\\
&0\qquad\h{class $(\gimel)$ or $(\daleth)$}
\endaligned
\right.
$$
\end{thm}

The result for class $(\aleph)$ is shown by Berman \cite{Berman2010}, and
the remaining cases are shown in \cite{CR}.
Note that $0,1/2$ and $1$ are the Tian
invariants of $\PP^n, n\ra\infty,
\PP^1$, and $\PP^0$, respectively. It is then tempting to think
of $1/2$ as the Tian invariant of the generic rational fiber
of Proposition \ref{proposition:conic-bundle}, thus
suggesting existence of approximate conic metrics on the
football fibers, who should tend to cylinders in the limit.
On the other hand, the smallness of the log canonical threshold
for classes $(\gimel)$ and $(\daleth)$ suggests non-existence.

\subsection{Third motivation: explicit computations}

Tian's 1990 result in the smooth regime
mentioned earlier can also be phrased
equivalently by saying that a del Pezzo surface
admits a KE metric if and only if its automorphism
group is reductive (a simplification of Tian's original
proof has been obtained by the work of Cheltsov \cite{Cheltsov2001}
and Shi \cite{ShiYalong2010}, see also Odaka--Spotti--Sun 
\cite{OdakaSpottiSun}, and the expository article \cite{Tosatti-delPezzo}).
Given the logarithmic version of
Matsushima's criterion (Theorem \ref{reductiveThm}), it
is tempting to check how far reductivity gets us in the asymptotic
regime. Some explicit computations give \cite{CR}:

\begin{prop}
\label{NoKEEThm} The automorphism groups of all pairs
of class $(\aleph)$ and $(\beth)$ are reductive. The
pairs of classes $(\gimel)$ and $(\daleth)$
that have non-reductive automorphism groups, and hence admit
no KEE metrics, are: $\mathrm{(I.1C})$, $\mathrm{(I.2.n})$
with any $n\ge 0$, $\mathrm{(I.6C.m})$ with any $m\ge 1$,
$\mathrm{(I.7.n.m})$ with any $n\ge 0$ and $m\ge 1$,
$\mathrm{(I.6B.1})$, $\mathrm{(I.8B.1})$ and $\mathrm{(I.9C.1})$.
\end{prop}

Thus, Matsushima's criterion supports
Conjecture \ref{UniformizationLogDelPezzoConj} but
does not solve the problem in the singular setting.

Another tool is Tian's criterion for existence of KEE metrics,
which involves calculating log canonical threshold, and is especially
useful in the presence of large finite symmetries. Using such
tools, the following KEE metrics are constructed \cite{CR} on surfaces
of class $(\beth)$:

\begin{thm}
\label{KEEBethThm}
There exist strongly asymptotically log del
Pezzo pairs of type $\mathrm{(I.3A), (I.4B),}$ and $\mathrm{(I.9B.5)}$
that admit KEE metrics for all sufficiently small $\be$.
\end{thm}

An in-depth study of log canonical thresholds on pairs of class $(\aleph)$
was carried out by Mart\'\inodot nez-Garc\'\inodot a \cite[\S4]{JesusThesis}
and Cheltsov--Mart\'\inodot nez-Garc\'\inodot a \cite{VanyaJesus}
whose results give lower bounds on how large $\be$ can be taken
(see also \cite{CR} for some weaker bounds that hold in all dimensions).
In fact, their results show that the threshold depends on the representative
chosen in \h{$|-K_S|$}, and, similarly, examples of Sz\'ekelyhidi
\cite{Szekelyhidi2012} show the maximal allowed $\be$ might as well.
Conjecture \ref{UniformizationLogDelPezzoConj} predicts that such dependence
does not appear in the asymptotic regime.

In work in progress \cite{CR2}, the `only if' part of
Conjecture \ref{UniformizationLogDelPezzoConj} is verified by
using techniques adapted from work of Ross--Thomas on slope stability,
motivated in part by work of Li--Sun \cite{LiSun} that proved
non-existence in the small angle regime for $\h{(I.1B)},
\h{(I.3B)},$ and  $\h{(I.4C)}$.

\subsubsection*{Acknowledgements}

This article is an expanded version of a talk
delivered at the CRM, Montr\'eal, in July 2012.
I
am grateful to D. Jakobson for the invitation, the hospitality,
and the encouragement to write this article,
and
to the city of Montr\'eal for its inspiring Jazz Festival.
I have benefited from many discussions on these topics with many people over
the years, and in particular R. Berman,
S. Donaldson, and C. Li.
I thank
L. Di Cerbo and J. Mart\'\inodot nez-Garc\'\inodot a
for corrections, 
C. Ciliberto and I. Dolgachev for historical references on classical
algebraic geometry,
R. Conlon, 
T. Darvas, S. Dinew, H. Guenancia, 
T. Murphy, S. Zelditch, and the referees for comments improving the exposition,
and O. Chodosh and D. Ramos for creating many of the figures.
Sections \ref{KahlerSection}, \ref{RCMSec}, and \ref{APrioriSubSec}
are largely based on joint work with R. Mazzeo
and I am grateful to him for everything he has taught me
during our collaboration. 
Similarly, Sections \ref{AsympSection}--\ref{logCalabiSec} survey joint
work with I. Cheltsov and I am grateful to him for
an inspiring collaboration.
I am indebted to Gang Tian for introducing me to this subject,
for countless discussions and lectures over many years, encouragement,
and inspiration.
Finally, this survey is dedicated
to Eugenio Calabi on the occasion of his ninetieth birthday.
I feel privileged to have known him since my graduate
school days. He continues to be an inspiring r\^ole model
exemplifying the spirit of free and open exchange of ideas in
a (mathematical) world that has become increasingly competitive.
This research was supported by NSF grants DMS-0802923,1206284,
and a Sloan Research Fellowship.

\begin{spacing}{0}

\def\bi{\bibitem}

\end{spacing}

\bigskip

\bigskip

{\sc University of Maryland} 

{\tt yanir@umd.edu}

\end{document}